\documentclass{article}

\usepackage{microtype}
\usepackage{graphicx}
\usepackage{subfigure}
\usepackage{booktabs}

\usepackage{hyperref}

\usepackage[accepted]{icml2024}

\usepackage{amsmath}
\usepackage{amssymb}
\usepackage{mathtools}
\usepackage{amsthm}
\usepackage{array}

\usepackage{url}
\usepackage{amsfonts}
\usepackage{nicefrac}
\usepackage{microtype}
\usepackage{xcolor}
\usepackage{wrapfig}
\usepackage{graphicx}
\usepackage{multirow}
\usepackage{enumerate}
\usepackage{arydshln}
\usepackage{multicol}

\usepackage{makecell}
\usepackage{graphicx}

\usepackage[capitalize,noabbrev]{cleveref}

\theoremstyle{plain}
\newtheorem{theorem}{Theorem}[section]

\newtheorem{lemma}[theorem]{Lemma}
\newtheorem{corollary}[theorem]{Corollary}
\theoremstyle{definition}

\newtheorem{assumption}[theorem]{Assumption}
\theoremstyle{remark}
\newtheorem{remark}[theorem]{Remark}

\usepackage[textsize=tiny]{todonotes}

\icmltitlerunning{
SPABA: A Single-Loop and Probabilistic Stochastic Bilevel Algorithm Achieving Optimal Sample Complexity
}

\begin{document}

\twocolumn[
\icmltitle{SPABA: A Single-Loop and Probabilistic Stochastic Bilevel Algorithm Achieving Optimal Sample Complexity}

\icmlsetsymbol{equal}{*}

\begin{icmlauthorlist}
	\icmlauthor{Tianshu Chu}{bjut}
 	\icmlauthor{Dachuan Xu}{bjut}
	\icmlauthor{Wei Yao}{sustech,sustech2}
	\icmlauthor{Jin Zhang}{sustech2,sustech}
\end{icmlauthorlist}

\icmlaffiliation{sustech}{National Center for Applied Mathematics Shenzhen, Southern University of Science and Technology, Shenzhen, China}

\icmlaffiliation{bjut}{Institute of Operations Research and Information Engineering, Beijing University of Technology, Beijing 100124, P.R. China}

\icmlaffiliation{sustech2}{Department of Mathematics, Southern University of Science and Technology, Shenzhen, China }

\icmlcorrespondingauthor{Jin Zhang}{zhangj9@sustech.edu.cn}

\icmlkeywords{Machine Learning, ICML}

\vskip 0.3in
]

\printAffiliationsAndNotice{}

\begin{abstract}
While stochastic bilevel optimization methods have been extensively studied for addressing large-scale nested optimization problems in machine learning, it remains an open question whether the optimal complexity bounds for solving bilevel optimization are the same as those in single-level optimization. Our main result resolves this question: SPABA, an adaptation of the PAGE method for nonconvex optimization in \cite{li2021page} to the bilevel setting, can achieve optimal sample complexity in both the finite-sum and expectation settings. We show the optimality of SPABA by proving that there is no gap in complexity analysis between stochastic bilevel and single-level optimization when implementing PAGE. Notably, as indicated by the results of \cite{dagreou2022framework}, there might exist a gap in complexity analysis when implementing other stochastic gradient estimators, like SGD and SAGA. In addition to SPABA, we propose several other single-loop stochastic bilevel algorithms, that either match or improve the state-of-the-art sample complexity results, leveraging our convergence rate and complexity analysis. Numerical experiments demonstrate the superior practical performance of the proposed methods.
\end{abstract}

\section{Introduction}
\label{introduction}

Bilevel optimization, where one optimization problem is nested within the constraints of another, has extensive applications in fields such as transportation \cite{yang2001transport} and game theory \cite{von1952theory}.
In recent years, bilevel optimization has gained popularity in the machine learning community due to its broad range of applications, including hyperparameter optimization \cite{pedregosa2016hyperparameter,mackay2019self,lorraine2020optimizing}, meta-learning \cite{franceschi2018bilevel,ji2020convergence}, and neural architecture search \cite{liu2018darts,liang2019darts+}. Refer to recent survey papers \cite{liu2021investigating,zhang2023introduction} for more applications of bilevel optimization in machine learning, computer vision and signal processing.

Bilevel optimization tackles challenges arising from hierarchical optimization, where decision variables in the upper level are also involved in the lower level.
Typically, the bilevel optimization problems are formulated as
\begin{align}
\min_{x\in \mathbb{R}^{d_x}}\
& H(x):=f(x,y^*(x)) \label{UL}\\
\text{s.t.}\ \ \,  &
y^*(x):=\arg \min_ {y\in \mathbb{R}^{d_y}}g(x,y) , \label{LL}
\end{align}
where the upper-level (UL) objective $f(x,y)$ and the lower-level (LL) objective $g(x,y)$ are two smooth real valued functions defined on $\mathbb{R}^{d_x}\times \mathbb{R}^{d_y} $.
In this work, we focus on the setting where the LL objective $g(x,y)$ is strongly convex with respect to (w.r.t.) $y$ for any $x$, and the UL objective $f(x,y)$ is possibly nonconvex.

A commonly employed strategy for solving bilevel problems involves utilizing implicit differentiation, which yields the following expression for the hypergradient:
\begin{equation}\label{hypergradient-1}
\nabla H(x) =\nabla_1f(x,y^*(x))-\nabla_{12}^2g(x,y^*(x))z^*(x),
\end{equation}
where $z^*(x)=\left[\nabla_{22}^2g(x,y^*(x))\right]^{-1}\nabla_2 f(x,y^*(x))$.
The practical implementation of the gradient descent method encounters several challenges, including: the computation of the exact solution $y^*(x)$ at the lower level, and the inversion of the Hessian $\nabla_{22}^2 g$ at the point $(x, y^*(x))$.
Utilizing a warm start strategy in the LL updates, results by \cite{ji2021bilevel,Averaged2023} demonstrate that deterministic bilevel algorithms based on approximate implicit differentiation (AID) can achieve a convergence rate of $\mathcal{O}(\epsilon^{-1})$.
The convergence rate matches that of the gradient descent method for nonconvex single-level optimization.

However, deterministic approaches necessitate the evaluation of the full gradient at every iteration, demanding substantial computational resources. This drawback renders these methods unsuitable for large-scale machine learning tasks.
In many applications of interest, the objective functions $f$ and $g$ have the \textbf{finite-sum form}:
\begin{gather*}\label{finipro}
    f(x,y)=\frac{1}{n} \sum_{i=1}^n F_i\left(x, y \right),
    g(x,y)=\frac{1}{m} \sum_{j=1}^m G_j\left(x, y \right),
\end{gather*}
which captures the standard empirical risk minimization problems in machine learning.
Additionally, when dealing with a substantial or potentially infinite number of data samples, such as in online or streaming scenarios, $f$ and $g$ are commonly represented using the \textbf{expectation form}:
\begin{gather*}\label{generalpro}
f(x,y)=\mathbb{E}_{\xi}\left[F\left(x, y ; \xi\right)\right],
g(x,y)=\mathbb{E}_{\zeta}\left[G\left(x, y ; \zeta\right)\right].
\end{gather*}

To improve sample efficiency compared to full-batch methods, it makes sense to apply stochastic techniques from single-level optimization to the bilevel optimization context.
Unfortunately, the practical implementation of stochastic algorithms faces various challenges, such as: computing exactly the solution $y^*(x)$ at the lower level; inverting the inversion of the Hessian $\nabla_{22}^2 g$; and addressing the nonlinear characteristics of $\nabla H$ within functions $f$ and $g$.
Therefore, a natural question follows: {\it Whether the optimal complexity bounds for solving bilevel optimization are the same as those in single-level optimization?} In fact, even more basic questions are open: \\
\textbf{Question 1.}
{\it Is there a gap in complexity analysis between stochastic bilevel
and single-level optimization when implementing the same stochastic gradient estimator?}

In the literature, various strategies have been proposed to tackle these challenges. For instance, the existing methods \cite{ghadimi2018approximation,ji2021bilevel,yang2021provably,chen2021closing,guo2021stochastic,khanduri2021near,hong2023two} employ one or multiple iterations of stochastic gradient descent (SGD) for the LL problem while incorporating truncated stochastic Neumann series to approximate the Hessian inversion. However, the mentioned methods suffer from an additional factor of $\log(\epsilon^{-1})$ in both sample complexity and batch size. Hence, there exists a gap in complexity analysis between stochastic bilevel and single-level optimization when employing the same stochastic gradient estimator in the aforementioned methods.

To address the nonlinear characteristics of $\nabla H$ within functions $f$ and $g$ and avoid relying on the stochastic Neumann approximation, recent works \cite{arbel2022amortized, dagreou2022framework} have employed the decoupling approach, see Section~\ref{framwork1} for more details. This approach breaks down the hypergradient computation into three gradient estimates, as outlined in (\ref{DX})-(\ref{DZ}).
When utilizing the framework presented in \cite{arbel2022amortized, dagreou2022framework}, several studies have indicated that stochastic bilevel algorithms exhibit similar sample complexity when compared to their single-level counterparts.
For instance, in the context of \textbf{general expectation setting}, by implementing multiple SGD iterations in subroutines and using large batchsizes of order $\mathcal{O}(\epsilon^{-1})$, AmIGO from \cite{arbel2022amortized} achieves the same sample complexity $\mathcal{O}(\epsilon^{-2})$ to SGD for smooth nonconvex single-level optimization, that is required to get an $\epsilon$-stationary point, defined as $\mathbb{E}\|\nabla H(x)\|^2\leq \epsilon$ \cite{ghadimi2013stochastic}.
When a single SGD step is used, it results in SOBA \cite{dagreou2022framework}.
However, by the result in Appendix D of \cite{dagreou2022framework}, SOBA achieves a sample complexity of $\mathcal{O}(\epsilon^{-2.5})$ under standard smoothness assumptions.
Thus, the gap in complexity analysis between stochastic bilevel and single-level optimization, using the SGD gradient estimator, is on the order of $\mathcal{O}(\epsilon^{-0.5})$. Recently, this gap has been effectively addressed by MA-SOBA \cite{chen2023optimal}, which builds upon the SOBA algorithm by incorporating an additional standard momentum (also referred to as moving average) into the update of the UL variable.

Despite the simplicity and power of MA-SOBA, we lack a comprehensive understanding of Question 1 regarding other stochastic gradient estimators.  Specifically, if we additionally assume that the stochastic gradient satisfies a mean-squared smoothness property, as commonly assumed in the existing literature \cite{yang2021provably,khanduri2021momentum},
the lower bound for nonconvex stochastic single-level optimization can be improved to $\mathcal{O}(\epsilon^{-1.5})$ \cite{arjevani2023lower}.
It is natural to ask the following question: \\
\textbf{Question 2.}
{\it How to develop a fully single-loop algorithm for solving stochastic
bilevel optimization problems that achieves an optimal sample complexity $\mathcal{O}(\epsilon^{-1.5})$ under the bounded variance and mean-squared smoothness?}

In the \textbf{finite-sum setting}, as indicated by the result of SABA in Appendix D of \cite{dagreou2022framework}, there also exists a gap of order $\mathcal{O}((n+m)^{1/3})$ between stochastic bilevel and single-level optimization in complexity analysis, when using the SAGA gradient estimator.
In a recent work \cite{dagreou2023lower}, the authors introduce SRBA, which is a bilevel extension of the well-known SARAH algorithm \cite{nguyen2017sarah}. They demonstrate that SRBA achieves a better sample complexity of $\mathcal{O}((n+m)^{1/2}\epsilon^{-1})$, matching the lower bound they established for bilevel optimization.
Unfortunately, the current analysis of SRBA relies on the assumption of higher-order smoothness for both the UL and LL functions to achieve optimality. It is also worth noting that SRBA in \cite{dagreou2023lower} utilizes a double-loop structure. Consequently, natural questions arise:\\
\textbf{Question 3.}
{\it Is it possible to fill the gap between stochastic bilevel and single-level optimization when using the SAGA? How to develop a fully single-loop algorithm for solving stochastic
bilevel optimization problems that achieves an optimal sample complexity $\mathcal{O}((n+m)^{\frac{1}{2}}\epsilon^{-1})$ under standard
smoothness assumptions in the finite-sum setting?}

\begin{table*}[t]
\renewcommand\arraystretch{1.5}
	\caption{Comparison of our methods with closely related works for {\bf nonconvex-strongly-convex} bilevel optimization under standard smoothness assumptions, without relying on high-order smoothness.
 The $\tilde{\mathcal{O}}$ notation hides a factor of $\log(\epsilon^{-1})$.
 The sample complexity corresponds to the number of calls made to stochastic gradients and Hessian (Jocobian)-vector products required to get an $\epsilon$-stationary point, i.e., $\mathbb{E}\|\nabla H(x)\|^2\leq \epsilon$.
 $^*:$ This result can be found in Appendix D of \cite{dagreou2022framework}.
	}
	\label{table1+}
	\vskip 0.15in
	\begin{center}
		\begin{small}
				\begin{tabular}{cccccccc}
					\toprule
                    Setting & Method &
                    \makecell{
                    Stochastic Estimators}  & Sample Complexity  & Batch Size\\\cline{1-8}
\multirow{5}*{ \makecell{Expection\\(Mean-squared\\ smoothness)}}
                    & \makecell{MRBO \cite{yang2021provably}} &  STORM & $\tilde{\mathcal{O}}(\epsilon^{-1.5})$ & $\tilde{\mathcal{O}}(1)$
                \\\cline{2-8}
                    & \makecell{SUSTAIN \cite{khanduri2021near}} &  STORM &$\tilde{\mathcal{O}}(\epsilon^{-1.5})$ & $\tilde{\mathcal{O}}(1)$
                \\\cline{2-8}
                   & \makecell{VRBO \cite{yang2021provably}}   & SARAH  & $\tilde{\mathcal{O}}(\epsilon^{-1.5})$ &  $\tilde{\mathcal{O}}(\epsilon^{-0.5})$
                   \\\cline{2-8}
                   & SRMBA (Ours)   & STORM  & $\tilde{\mathcal{O}}(\epsilon^{-1.5})$ &  ${\mathcal{O}}(1)$
                   \\\cline{2-8}
                   &  SPABA (Ours)  & PAGE  & $\mathcal{O}(\epsilon^{-1.5})$ &  ${\mathcal{O}}(1)$
\\\cline{1-8}
\multirow{3}*{ \makecell{Finite-Sum}}
					& \makecell{SABA \cite{dagreou2022framework}} & SAGA & $\mathcal{O}((n+m)\epsilon^{-1})^*$ & $\mathcal{O}(1)$
                \\\cline{2-8}
                & MA-SABA (Ours) & \makecell{SAGA \\
                    + $x$-Momentum} & $\mathcal{O}((n+m)^{\frac{2}{3}}\epsilon^{-1})$ &  $\mathcal{O}(1)$
                \\\cline{2-8}
                &SPABA (Ours)  & PAGE &  $\mathcal{O}((n+m)^{\frac{1}{2}}\epsilon^{-1})$ & $\mathcal{O}(1)$  \\
					\bottomrule
				\end{tabular}
		\end{small}
	\end{center}
	\vskip -0.1in
\end{table*}

\subsection{Main Contribution}

The purpose of this work is to understand these theoretical questions. Our contributions are summarized below.

\begin{itemize}
\item
\textbf{Bridging the gap between stochastic
bilevel and single-level optimization when using the SAGA.}
We first introduce a single-loop stochastic bilevel algorithm, named MA-SABA, that achieves a sample complexity of $\mathcal{O}((n+m)^{\frac{2}{3}}\epsilon^{-1})$ without the need for high-order smoothness.
It is worth noting that MA-SABA is based on SABA and inspired by MA-SOBA by integrating an additional standard momentum into the update of the UL variable.

\item \textbf{Achieving the optimal sample complexity in both the finite-sum and general expectation setting.}
We propose a fully single-loop and sample-efficient stochastic bilevel algorithm, called SPABA, that achieves an optimal sample complexity of $\mathcal{O}((n+m)^{\frac{1}{2}}\epsilon^{-1})$ under standard smoothness assumptions in the finite-sum scenario. Additionally, it attains optimal sample complexity of $\mathcal{O}(\epsilon^{-1.5})$ under the bounded variance and mean-squared smoothness in the general expectation context.
Technically, SPABA is an adaptation of the PAGE algorithm in \cite{li2021page} to the bilevel setting.

\item
\textbf{Convergence rate and complexity analysis.}
It is often difficult to analyze biased stochastic algorithms.
We provide a general and unified convergence rates and complexity analysis based on biased stochastic gradient estimator such as STORM and PAGE, which either match or improve the state-of-the-art sample complexity results.

\item  Finally, numerical experiments demonstrate the superior efficiency of our
proposed methods in bilevel optimization.

\end{itemize}

\subsection{Additional Related Work}

In the section we give a brief review of some recent works that are directly related to ours.  A summary of the comparison of the proposed
methods with closely related works is provided in Table \ref{table1+}.

\textbf{Lower Bounds for Stochastic Bilevel Optimization.}
When $H(x)$ be convex or strongly convex, the study \cite{ji2021lower} has provided lower complexity bounds for deterministic bilevel optimization, that are larger than the corresponding optimal complexities of minimax optimization. For non-convex stochastic bilevel optimization,
since nonconvex optimization can be regarded as a specific instance of a bilevel problem, it is natural to consider that lower bounds for nonconvex stochastic optimization also apply as lower bounds for bilevel counterparts.
Therefore, the $\mathcal{O}(\epsilon^{-2})$ complexity is a lower bound for non-convex stochastic bilevel optimization in general expectation setting \cite{arjevani2023lower}.
Such complexity is attained by SGD in nonconvex stochastic optimization \cite{ghadimi2013stochastic}.
If we additionally assume that the stochastic gradient satisfies a mean-squared smoothness property, the lower bound is improved to $\mathcal{O}(\epsilon^{-1.5})$ \cite{arjevani2023lower},
which is attained in nonconvex stochastic optimization by SPIDER \cite{fang2018spider}, Spiderboost \cite{wang2018spiderboost}, SARAH \cite{pham2020proxsarah}, and PAGE \cite{li2021page}.
Moreover, with the additional assumption of Lipschitz continuity, STORM \cite{cutkosky2019momentum} can also reach this complexity level.
In the nonconvex
finite-sum setting, if we assume that the objective function is averagely $L$-smooth, the lower bound becomes $\Omega(n^{1/2}\epsilon^{-1})$ \cite{fang2018spider,li2021page}. Such complexity has been achieved by SARAH \cite{nguyen2017sarah,nguyen2017stochastic,pham2020proxsarah}, SPIDER \cite{fang2018spider}, and PAGE \cite{li2021page}.

\textbf{Discussion under Stronger Smoothness Conditions.}
Some studies have been conducted based on stronger smoothness conditions, such as SOBA and SABA in \cite{dagreou2022framework}. Indeed, when the UL and LL objective functions possess high-order smoothness, their study illustrates that SABA, an adaptation of the SAGA algorithm \cite{defazio2014saga}, exhibits a sample complexity of $\mathcal{O}((n+m)^{2/3}\epsilon^{-1})$. This is consistent with the sample complexity of SAGA in the single-level counterpart.
Recently, leveraging on high-order smoothness, SRBA \cite{dagreou2023lower}, an adaptation of the SARAH algorithm to the bilevel setting, achieves the same complexity $\mathcal{O}((n+m)^{1/2}\epsilon^{-1})$ as single-level SARAH.
It is unclear whether a gap exists between stochastic bilevel and single-level optimization when utilizing the SARAH gradient estimator under standard smoothness assumptions.

\section{The Proposed Stochastic Bilevel Algorithms}

\subsection{Overview of the Framework in \cite{arbel2022amortized, dagreou2022framework}}
\label{framwork1}

In this section, we provide an overview of the algorithm design. First, we review the decoupling method employed in \cite{arbel2022amortized, dagreou2022framework}. To handle the nonlinear characteristics of $\nabla H$ within functions $f$ and $g$,
the authors in \cite{dagreou2022framework} introduce an extra variable $z\in\mathbb{R}^{d_y}$ to effectively decouple the nonlinear structure in $\nabla H$. This allows us to utilize $\nabla_1 f(x,y)
-\nabla_{12}^2 g(x,y) z$
to approximate the hypergradient $\nabla H(x)$, where $y$ represents an approximate solution to the LL problem, while $z$ serves as an inexact solution to the linear system $\left[\nabla_{22}^2g(x,y)\right] z - \nabla_2 f(x,y)=0$, which can also be seen as optimizing the following quadratic problem:
\begin{eqnarray}\label{z*}
\min_{z}\
\frac{1}{2}\langle\nabla_{22}^2g(x,y)z,z\rangle-\langle\nabla_2 f(x,y),z\rangle.
\end{eqnarray}

In summary, to solve the upper-level optimization problem $\min H(x)$, we decompose the search direction (or hypergradient estimate) of $x$ into three steps, as follows:
\begin{align}
D_x(x,y,z) =& \nabla_1 f(x,y)-\nabla_{12}^2g(x,y)z,\label{DX} \\
D_y(x,y,z) =& \nabla_2 g(x,y),           \label{DY} \\
D_z(x,y,z) =& \nabla_{22}^2g(x,y)z- \nabla_2 f(x,y).     \label{DZ}
\end{align}
Notably, all search directions are linear within functions $f$ and $g$. The latter two directions align with two strongly convex optimization problems: the lower-level optimization problems (\ref{LL}) and (\ref{z*}). In addition, as detailed in Section 2.1 of \cite{Averaged2023}, the search directions presented in (\ref{DX}-\ref{DZ}) precisely correspond to the KKT condition of the equality-constrained optimization
reformulation of (\ref{UL}):
\begin{equation*}
    \min_{x,y} \ f(x,y)\quad
    \mathrm{s.t.}
    \quad
    \nabla_2 g(x,y)=0.
\end{equation*}
Consequently, $z$ can be interpreted as the dual multiplier.

Now, we provide a comprehensive description of the framework in \cite{arbel2022amortized, dagreou2022framework}, referred to as the Decoupling stochastic Bilevel Optimizer (DecBO).
In each iteration, we sample $\mathcal{S}_k^f$ for $f(x,y)$ and $\mathcal{S}_k^g$ for $g(x,y)$. We then construct unbiased or biased stochastic estimators, denoted as $v_k^x$, $v_k^y$ and $v_k^z$, for $D_x(x_k, y_k, z_k)$, $D_y(x_k, y_k, z_k)$ and $D_z(x_k, y_k, z_k)$ in equations (\ref{DX})-(\ref{DZ}), respectively. These gradient estimators are constructed using the samples from $\mathcal{S}_k^f$ and $\mathcal{S}_k^g$, as well as past gradient estimators.
We provide a pseudo code to illustrate this (see Algorithm \ref{alg1}).

\begin{algorithm}[tb]
  \caption{Pseudocode for a generic Decoupling stochastic Bilevel Optimizer (DecBO)}
  \label{alg1}
\begin{algorithmic}[1]
   \STATE {\bfseries Input:} Initializations $(x_{-1},y_{-1},z_{-1})$ and $(x_{0},y_{0},z_{0})$, number of total iterations $K$, step size $\{\alpha_k, \beta_k,\gamma_k\}$;
   \FOR{$k=0$ {\bfseries to} $K-1$}
   \STATE Sample $\mathcal{S}_k^f$ for $f$ and $\mathcal{S}_k^g$ for $g$;

    \STATE Construct an unbiased or biased estimator $v_k^x$ of $D_x(x_k, y_k, z_k)$ in (\ref{DX}) using $\mathcal{S}_k^f, \mathcal{S}_k^g$ and past gradient estimators;
   \STATE Update
   \begin{equation}
       x_{k+1} \leftarrow x_{k}-\alpha_k v_k^x;
   \end{equation}

   \STATE Construct an unbiased or biased estimator $v_k^y$ of $D_y(x_k, y_k, z_k)$ in (\ref{DY}) using $\mathcal{S}_k^g$ and past gradient estimators;
   \STATE Update
   \begin{equation}
       y_{k+1} \leftarrow y_{k}-\beta_k v_k^y;
   \end{equation}

   \STATE Construct an unbiased or biased estimator $v_k^z$ of $D_z(x_k, y_k, z_k)$ in (\ref{DZ}) using $\mathcal{S}_k^f, \mathcal{S}_k^g$ and past gradient estimators;
   \STATE Update
   \begin{equation}
       z_{k+1} \leftarrow z_{k}-\gamma_k v_k^z.
   \end{equation}
   \ENDFOR
\end{algorithmic}
\end{algorithm}

The proposed framework opens opportunities for developing new algorithms in stochastic bilevel optimization. These algorithms can integrate diverse stochastic gradient estimation techniques from stochastic single-level optimization. For example, the aforementioned unbiased or biased gradient estimators can be efficiently constructed by combining variance-reduced gradient estimators like SAGA, SVRG, SPIDER or SARAH with momentum. Alternatively, one can utilize accelerated variance-reduced gradient estimators such as STORM or PAGE.
We focus in this work on the loopless variance-reduced estimators because they share handy theoretical properties.
As a result, the framework DecBO also benefits from a loopless structure. In the subsequent sections, we delve into the study of three such techniques.

\subsection{MA-SABA: Bridging the Gap between Stochastic
Bilevel and Single-level Optimization when Using the SAGA}

For the finite-sum setting, we present MA-SABA, which is based on SABA \cite{dagreou2022framework} and inspired by MA-SOBA \cite{chen2023optimal} by integrating an additional standard momentum into the update of the UL variable.

The SAGA method \cite{defazio2014saga} achieves variance reduction by updating historical gradients and performing gradient correction. Define two memory variables $w_{k,i}=(w_{k,i}^x,w_{k,i}^y, w_{k,i}^z)$ for $i \in [n]$ and $ w_{k,j}=(w_{k,j}^x,w_{k,j}^y, w_{k,j}^z)$  for $j \in[m]$ corresponding to calls to $f$ and $g$, respectively.
At each iteration $k$, we draw two random independent indices $i \in [n] $ and $j \in [m]$ uniformly, for $i^{\prime}\neq i$, do
\begin{eqnarray*}
\left\{\begin{array}{l}
(w_{k+1,i}^x,w_{k+1,i}^y, w_{k+1,i}^z)=(x_k,y_k,z_k),\\
(w_{k+1,i^{\prime}}^x,w_{k+1,i^{\prime}}^y, w_{k+1,i^{\prime}}^z)=(w_{k,i^{\prime}}^x,w_{k,i^{\prime}}^y, w_{k,i^{\prime}}^z),
\end{array}\right.
\end{eqnarray*}
and similarly for $ (w_{k+1,j}^x,w_{k+1,j}^y, w_{k+1,j}^z)$ .

At each iteration $k$, we randomly select $i_k\in [n]$ and $j_k\in[m]$.
In order to facilitate gradient correction, for $u_k:=(x_k,y_k,z_k)$, we define two operations
\begin{align*}
\mathcal{M}_f(\phi,k,u,w):=&\phi_{i_k}(x_k,y_k)-\phi_{i_k}(w_{k,i_k}^x,w_{k,i_k}^y)\\&+\sum_{i=1}^n\frac{\phi_{i}(w_{k,i}^x,w_{k,i}^y)}{n} ,\\
\mathcal{M}_g(\phi,k,u,w):=&\phi_{j_k}(u_k)-\phi_{j_k}(w_{k,j_k})+\sum_{j=1}^m\frac{\phi_{j}(w_{k,j})}{m}.
\end{align*}

Then we update $x_k$ using an additional
standard momentum.
The specific form of the iteration directions of MA-SABA are as follows:
\begin{align*}
    v^y_{k}=&\mathcal{M}_g(\nabla_{2}G,k,u,w),\\
    v^z_{k}=&\mathcal{M}_g(\nabla_{22}^2 Gz,k,u,w)-\mathcal{M}_f(\nabla_2 F,k,u,w),\\
    v^x_k=&(1-\rho_{k-1})v_{k-1}^x+\rho_{k-1} \mathcal{M}_f(\nabla_1 F,k-1,u,w)\\
    &-\rho_{k-1}\mathcal{M}_g(\nabla_{12}^2 Gz,k-1,u,w).
\end{align*}

\subsection{SPABA: Stochastic ProbAbilistic Bilevel Algorithm}

Now we introduce SPABA, an adaptation of the PAGE algorithm in \cite{li2021page} to the
bilevel setting. To start, we present the algorithm description within the finite sum setting.
During each iteration, we sample $I\subset [n]$ for $f$ and $J\subset [m]$ for $g$, with a minibatch size of $b$.
The PAGE method is utilized for stochastic gradient estimators in all three directions as follows:
\begin{align*}
    v^y_{k}&=v_k(\nabla_{2}G;b),\\
    v^z_{k}&=v_k(\nabla_{22}^2 Gz;b)-v_k(\nabla_2 F;b),\\
    v^x_k&=v_k(\nabla_1 F;b)-v_k(\nabla_{12}^2 Gz;b),
\end{align*}
where $\phi(u_k;b)=\frac{1}{b}\sum_{i'\in I} \phi_{i'}(u_k)$ and
\begin{align*}
v_k(\phi;b)=
\begin{cases} \phi(u_k) &  \text{ w.p.   }\, p,\\
v_{k-1}^{x}+\phi(u_k;b)-\phi(u_{k-1};b)
&  \text{w.p.} \,1-p.
\end{cases}
\end{align*}
Recall that PAGE uses the vanilla minibatch SGD update with probability (w.p.) $p$, and reuses the previous gradient with a momentum-based minibatch SGD w.p. $1-p$.

Furthermore, similar to PAGE, SPABA is adaptable to the general expectation setting by replacing the full gradient with another vanilla minibatch SGD using a minibatch size of $\tau'$. Refer to Section \ref{algos} for more details.

\subsection{SRMBA: Stochastic Recursive Momentum Bilevel Algorithm}

The STORM method \cite{cutkosky2019momentum} does not require the maintenance of anchor points or the use of large batches. Next, we propose SRMBA, which is a combination of the idea of STORM and the framework DecBO in Algorithm \ref{alg1}.

At each iteration $k$, we randomly select $\xi$ and $\zeta$ for the functions $f$ and $g$, respectively. Define $D^y_{k}=D_y\left(x_{k},y_{k},z_{k};\zeta\right)$, $D^z_{k}=D_z\left(x_{k},y_{k},z_{k};\xi,\zeta\right)$ and $D^x_{k}=D_x\left(x_{k},y_{k},z_{k};\xi,\zeta\right )$. The iteration directions of SRMBA take the specific form as follows:
\begin{align*}
    v^y_k&=\rho^y_k D_{k}^y + \left(1-\rho^y_k\right)\left(D_k^y-D_{k-1}^y+v^y_{k-1}\right),\\
    v^z_k&= \rho^z_k D_{k}^z + \left(1-\rho^z_k\right)\left(D_k^z-D_{k-1}^z+v^z_{k-1} \right),\\
    v^x_k&= \rho^x_k D_{k}^x + \left(1-\rho^x_k\right)\left(D_k^x-D_{k-1}^x+v^x_{k-1}\right).
\end{align*}

\section{Complexity Analysis}

In this section, we will present the theoretical results for MA-SABA, SPABA and SRMBA, which
either match or improve the state-of-the-art sample complexity results.

We say that $\bar{x}$ is a $\epsilon$-stationary point if  $\mathbb{E}\|\nabla H(\bar{x})\|^2\leq \epsilon$.
The sample complexity corresponds to the total number of calls made to stochastic gradients and Hessian (Jocobian)-vector products required to get an $\epsilon$-stationary point.

\subsection{Structure Assumptions}

In order to provide convergence rates and complexity analysis, one usually needs the following standard assumptions depending on the setting \cite{ghadimi2018approximation,guo2021stochastic,yang2021provably,khanduri2021near,chen2021closing,arbel2022amortized,dagreou2022framework,hong2023two,chen2023optimal,huang2023momentumbased}.

\begin{assumption}\label{assump UL}
(1) $\nabla f$ is Lipschitz continuous in $(x,y)$ with Lipschitz constant $L^f$; (2) There exists $C^f>0$, such that $\|\nabla_2 f(x,y^*(x))\|\leq C^f$ for any $x$.
\end{assumption}

\begin{assumption}\label{assump LL}
(1) $\nabla g$ and $\nabla^2 g$ are $L^g_1$ and $L^g_2$ Lipschitz continuous in $(x,y)$, respectively;
(2) $g(x,\cdot)$ is $\mu$-strongly convex for any $x$.
\end{assumption}

Such assumptions are classical and sufficient to ensure the Lipschitz continuity of
$y^*(x)$ and $z^*(x)$, the boundedness of $z^*(x)$, and the $L$-smoothness of $H(x)$. Next, we discuss assumptions made on the stochastic oracles.

\begin{assumption}\label{assumporacle}
(\textbf{Bounded Variance})
In the general expectation setting, there exist positive constants $\sigma_f$, $\sigma_{g,1}$ and $\sigma_{g,2}$ such that
\begin{eqnarray*}
\mathbb{E}[\|\nabla F(x,y;\xi)-\nabla f(x,y)\|^2]\leq (\sigma_f)^2,\\
\mathbb{E}[\|\nabla G(x,y;\zeta)-\nabla g(x,y)\|^2]\leq (\sigma_{g,1})^2,\\
\mathbb{E}[\|\nabla^2 G(x,y;\zeta)-\nabla^2 g(x,y)\|^2]\leq (\sigma_{g,2})^2.
\end{eqnarray*}
\end{assumption}

Furthermore, to achieve a better sample complexity results, we need to adopt the mean-squared smoothness assumption in \cite{arbel2022amortized,chen2023optimal}.

\begin{assumption}\label{assumptionijsmooth}
(\textbf{Mean-Squared Smoothness})
Stochastic functions $\nabla F(x,y;\xi)$, $\nabla G(x,y;\zeta)$ and $\nabla^2G(x,y;\zeta)$ are $L^f$, $L^g_1$ and $L^g_2$ Lipschitz continuous in $(x,y)$, respectively.
\end{assumption}

\subsection{Convergence Analysis}\label{analysisframework}

We provide a general and unified convergence rates and complexity analysis and then illustrate it through the proposed methods. Let us identify what the crucial steps are.
A clearer exposition of the analytical process is provided in Figures \ref{fig:proofsketchunbaised} and \ref{fig:proofsketchbaised} of Appendix \ref{proofsketch}.

\textbf{General approach.} One of the most important steps is to establish a recursive estimate often generated by two or three consecutive iterates:
\begin{align}\label{ana}
    \alpha_k\mathbb{E}\left[\left\|\nabla H(x_k)\right\|^2\right]\leq L_k-L_{k+1}+\Delta_k,
\end{align}
where $L_k$, $L_{k+1}$ and $\Delta_k$ are all nonnegative quantities,
$\alpha_k$ is the step size used for updating $x_k$.
Denote $\theta=\min_{k\in[K]}\{\alpha_k\}$.
By induction, we have
\begin{align*}
    \frac{1}{K}\sum_{k=0}^{K-1}\mathbb{E}\left[\left\|\nabla H(x_k)\right\|^2\right]\leq \frac{L_0}{K\theta}+\frac{\sum_{k=0}^{K-1}\Delta_k}{K\theta}.
\end{align*}
This allows us to estimate the convergence rates of the underlying algorithm.

Usually, the recursive estimate (\ref{ana}) is derived through a series of recursive inequalities in conditional expectation:
\begin{equation}\label{recursiveinequality_text}
    \mathbb{E}\big[\widetilde{D}_{k+1}\,|\, \mathcal{F}_k\big] + \Lambda_k
    \leq \omega_k \widetilde{D}_{k} + \Omega_k,
\end{equation}
where $\widetilde{D}_{k}$, $\Lambda_k$, $\Omega_k$ are all nonnegative quantities, and $\omega_k\in[0,1]$ is a contraction factor.
We can now divide the proof of the recursive estimate (\ref{ana}) into four main steps:

\textbf{(1)} We begin by bounding the descent of $H(x)$ as follows:
\begin{equation}
\begin{aligned}
\frac{\alpha_k}{2} \mathbb{E}\big[\left\|\nabla H\left(x_k\right)\right\|^2\big]
\leq \mathbb{E}\left[H\left(x_k\right)\right]-\mathbb{E}\left[H\left(x_{k+1}\right)\right]\\+
 \Big(\frac{L^H\alpha_k^2}{2}-\frac{\alpha_k}{2}\Big) \mathbb{E}\big[\left\|v_k^x\right\|^2\big]
\\ +\frac{\alpha_k}{2} \mathbb{E}\big[\left\|\nabla H\left(x_k\right)-v_k^x\right\|^2\big],
\end{aligned}
\label{HinHOWi}
\end{equation}
which is a recursive inequality, as demonstrated in (\ref{recursiveinequality_text}). It is established in Lemma \ref{Hlemma} by the $L^{H}$-smoothness of $H(x)$.

\textbf{(2)}
We investigate the descent property of the mean-squared error $\mathbb{E}\big[\left\|\nabla H\left(x_k\right)-v_k^x\right\|^2\big]$ on the right-hand side of equation (\ref{HinHOWi}).
When integrating a standard momentum or a variation of momentum, such as those found in PAGE and STORM, into the update of $x_k$, we can establish a recursive inequality in the form of (\ref{recursiveinequality_text}) for $\widetilde{D}_{k}:=\left\|\nabla H\left(x_k\right)-v_k^x\right\|^2$. This inequality is derived from two or three consecutive iterates.
For example, this result is proven in Lemma \ref{vh} when standard momentum is utilized in $v_k^x$. It's important to highlight that the contraction factor $\omega_k=1-\rho_k$, where $\rho_k$ is the ``momentum" parameter and will tend to approach $0$ in the subsequent setting. If there is no momentum, it is only possible to obtain an upper bound for the mean-squared error.

\textbf{(3)} To gain better control over the $\Omega_k$-type terms in the descent of the mean-squared error $\mathbb{E}\big[\left\|\nabla H\left(x_k\right)-v_k^x\right\|^2\big]$,
it is essential to investigate the descent of $\mathbb{E}\big[\left\|D_x(x_k,y_k,z_k)-\nabla H\left(x_k\right)\right\|^2\big]$. Leveraging the key point presented in Lemma \ref{DH}, that
\begin{eqnarray*}
&&\mathbb{E}\big[\left\|D_x(x_k,y_k,z_k)-\nabla H\left(x_k\right)\right\|^2\big] \\
& &\leq c_1 \mathbb{E}\big[\left\|y_k-y^*\left(x_k\right)\right\|^2\big] +c_2 \mathbb{E}\big[\left\|z_k-z^*\left(x_k\right)\right\|^2\big],
\end{eqnarray*}
our analysis extends to a thorough examination of the descent of approximation errors $\mathbb{E}\big[\left\|y_k-y^*\left(x_k\right)\right\|^2\big]$ and $\mathbb{E}\big[\left\|z_k-z^*\left(x_k\right)\right\|^2\big]$.
By leveraging the strong convexity of the LL problem and the quadratic problem (\ref{z*}), one can readily derive recursive inequalities akin to (\ref{recursiveinequality_text}) for the approximate errors.
Additionally, $(1-\omega_k)$ exhibits a similar order of magnitude as the step sizes for both $y$ and $z$; please refer to Lemmas \ref{y,z} and \ref{eg30} for illustrations.

\textbf{(4)} In all the recursive inequalities mentioned above, the remaining terms include only
the variances of the stochastic gradient estimators, such as $\mathbb{E}\big[\left\|D_x(x_k,y_k,z_k)-D_k^x\right\|^2\big]$.
If the stochastic gradient estimators used lack variance reduction properties, like SGD, it is only feasible to attain a constant upper bound, even when we consider Assumption \ref{assumporacle}.
To further reduce sampling complexity, one can integrate unbiased or biased variance reduction techniques into the algorithm. For example, MA-SABA aligns its sampling complexity with that of single-level optimization using SAGA. For an illustration, please refer to Lemma \ref{var}.

\subsection{Convergence Results}\label{Results}

In this section, we provide the convergence results for the proposed stochastic bilevel algorithms. The detailed proofs of the results are deferred to the appendix.

We first provide the theoretical analysis of MA-SABA leading to a sample complexity in $\mathcal{O}((n+m))^{{2}/{3}}\epsilon^{-1})$ under standard smoothness assumptions in the finite-sum setting. This result bridges the gap between stochastic bilevel and single-level optimization when using the SAGA.

\begin{theorem}\textbf{(Convergence Rate of MA-SABA.)}\label{thsaba}
Fix an iteration $K>1$ and assume that Assumptions \ref{assump UL} to \ref{assump LL} and \ref{assumptionijsmooth} hold.
Then there exist positive constants $c_1$, $c_2$, $c_3$ and $c_4$ such that if $\alpha_k=c_1(n+m)^{-2/3}$, $\beta_k=c_2(n+m)^{-2/3}$, $\gamma_k=c_3(n+m)^{-2/3}$, $\rho_k=c_4 (n+m)^{-2/3}$, the iterates in MA-SABA satisfy
\begin{eqnarray*}
\frac{1}{K} \sum_{k=0}^{K-1} \mathbb{E}\left[\left\|\nabla H\left(x_k\right)\right\|^2\right]=\mathcal{O}\left((n+m)^{\frac{2}{3}}K^{-1}\right).
\end{eqnarray*}
\end{theorem}
\begin{remark}\it\textbf{ (Sample Complexity of MA-SABA.)}
To achieve the $\epsilon$-stationary point, the sampling complexity of MA-SABA is $\mathcal{O}((n+m))^{{2}/{3}}\epsilon^{-1})$, which is analogous to the sample complexity of SAGA in the nonconvex finite-sum setting.
\end{remark}

Next, we present the theoretical analysis of SPABA in both the finite-sum and general expectation settings.

\begin{theorem}\textbf{(Convergence Rate of SPABA in Finite-Sum Setting.)}\label{thpagefini}
Fix an iteration $K>1$ and assume that Assumptions \ref{assump UL} to \ref{assump LL} and \ref{assumptionijsmooth} hold.
Then there exist positive constants $c$, $c_{\beta}$, and $c_{\gamma}$, such that if
\begin{gather*}
    \alpha_k\leq  \frac{c}{1+\sqrt{\frac{1-p}{pb}}},\quad\beta_k=c_{\beta}\alpha_k,\quad  \gamma_k=c_{\gamma}\alpha_k,
\end{gather*}
the iterates in SPABA satisfy
$$ \frac{1}{K}\sum_{k=0}^{K-1}  \mathbb{E}\left[\left\|\nabla H\left(x_k\right)\right\|^2\right]
=\mathcal{O}\left(\frac{1+\sqrt{\frac{1-p}{pb}}}{K}\right).$$
\end{theorem}
\begin{remark}\it\textbf{(Sample Complexity of SPABA in Finite-Sum setting.)}
If we take $p=b/(n+m+b)$ and $b=\mathcal{O}((n+m)^{1/2})$, then the sample complexity of SPABA is $\mathcal{O}((n+m)^{1/2}\epsilon^{-1})$. This implies that there is no gap between stochastic bilevel
and single-level optimization in the context of PAGE implementation. The lower bound established in \cite{dagreou2023lower} for bilevel optimization indicates that SPABA attains optimal sample complexity in the finite-sum setting when $m=\mathcal{O}(n)$ and $\epsilon=\mathcal{O}(n^{-1/2})$.
\end{remark}

\begin{theorem}\textbf{(Convergence Rate of SPABA in Expectation Setting.)}\label{thpage1.5}
Fix an iteration $K>1$ and assume that Assumptions \ref{assump UL} to \ref{assumptionijsmooth} hold.
Choose minibatch size $\tau'$ and $b<\tau'$, the probability $p\in (0,1]$.
Then there exist positive constants $c$, $c_{\beta}$, $c_{\gamma}$ and $\sigma$, such that if
\begin{gather*}
    \alpha_k\leq  \frac{c}{1+\sqrt{\frac{1-p}{pb}}},\quad\beta_k=c_{\beta}\alpha_k,\quad  \gamma_k=c_{\gamma}\alpha_k,
\end{gather*}
the iterates in SPABA satisfy
\begin{gather*}
\frac{1}{K}\sum_{k=0}^{K-1}  \mathbb{E}\left[\left\|\nabla H\left(x_k\right)\right\|^2\right]
\\=\mathcal{O}\left(\frac{1+\sqrt{\frac{1-p}{pb}}}{K}+\frac{1}{Kp\tau'}+\frac{\sigma}{\tau'}\right).\end{gather*}
\end{theorem}
\begin{remark}\it\textbf{(Sample Complexity of SPABA in Expectation Setting.)}
    If we take $p=b/(n+m+b)$,  $\tau'=\mathcal{O}(\epsilon^{-1})$ and $b\leq\sqrt{\tau'}$, then the sample complexity of SPABA is $\mathcal{O}(\epsilon^{-1.5})$.
    This means that there is no gap between stochastic bilevel and single-level optimization when implementing PAGE. And SPABA achieves optimal sample complexity in the general expectation scenario.
\end{remark}

Finally, we state the convergence rate and sample complexity of SRMBA, an adaptation of the STORM method to the bilevel setting.

\begin{theorem}\textbf{(Convergence Rate of SRMBA in Expectation Setting.)}\label{thstorm}
Fix an iteration $K>1$ and assume that Assumptions \ref{assump UL} to \ref{assumptionijsmooth} hold.
Then there exist positive constants $\eta$, $c_{\beta}$, $c_{\gamma}$, $c_x$, $c_y$ and $c_z$ such that if
\begin{gather*}
    \alpha_k=\frac{1}{(\eta+k)^{1/3}},\quad \beta_k=c_{\beta}\alpha_k,\quad \gamma_k=c_{\gamma}\alpha_k;\\
    \rho_k^x=c_x \alpha_k^2,\quad  \rho_k^y=c_y \alpha_k^2,\quad \rho_k^z=c_z \alpha_k^2,
\end{gather*}
the iterates in SRMBA satisfy
\begin{align*}
\frac{1}{K}\sum_{k=0}^{K-1} \mathbb{E}\left[\left\|\nabla H\left(x_k\right)\right\|^2\right]
=\mathcal{O}\left(\frac{\log(K-1)}{K^{2/3}}\right).
\end{align*}
\end{theorem}
\begin{remark}\it\textbf{(Sample Complexity of SRMBA in Expectation Setting.)}
Theorem \ref{thstorm} implies that the sample complexity of SRMBA is $\mathcal{O}(\epsilon^{-1.5}\log(\epsilon^{-1}))$, which is analogous to the sample complexity of STORM in the nonconvex optimization. This tells us that there is no gap between stochastic bilevel and single-level optimization when implementing STORM.
\end{remark}

\section{Numerical Experiments}

While our contribution is mostly theoretical, we conducted a series of experiments to compare our proposed algorithms (SRMBA, SPABA, and MA-SABA) with their corresponding counterparts, namely, AmIGO \cite{arbel2022amortized}, SUSTAIN \cite{khanduri2021near}, SABA \cite{dagreou2022framework}, SOBA \cite{dagreou2022framework}, SRBA  \cite{dagreou2023lower}, MRBO \cite{yang2021provably}, and VRBO \cite{yang2021provably}.
Further elaboration on these experiments is available in the Appendix.

\begin{figure}[ht]
	\centering
	\subfigure[Datacleaning]{
		\centering
		\includegraphics[width=0.45\linewidth]{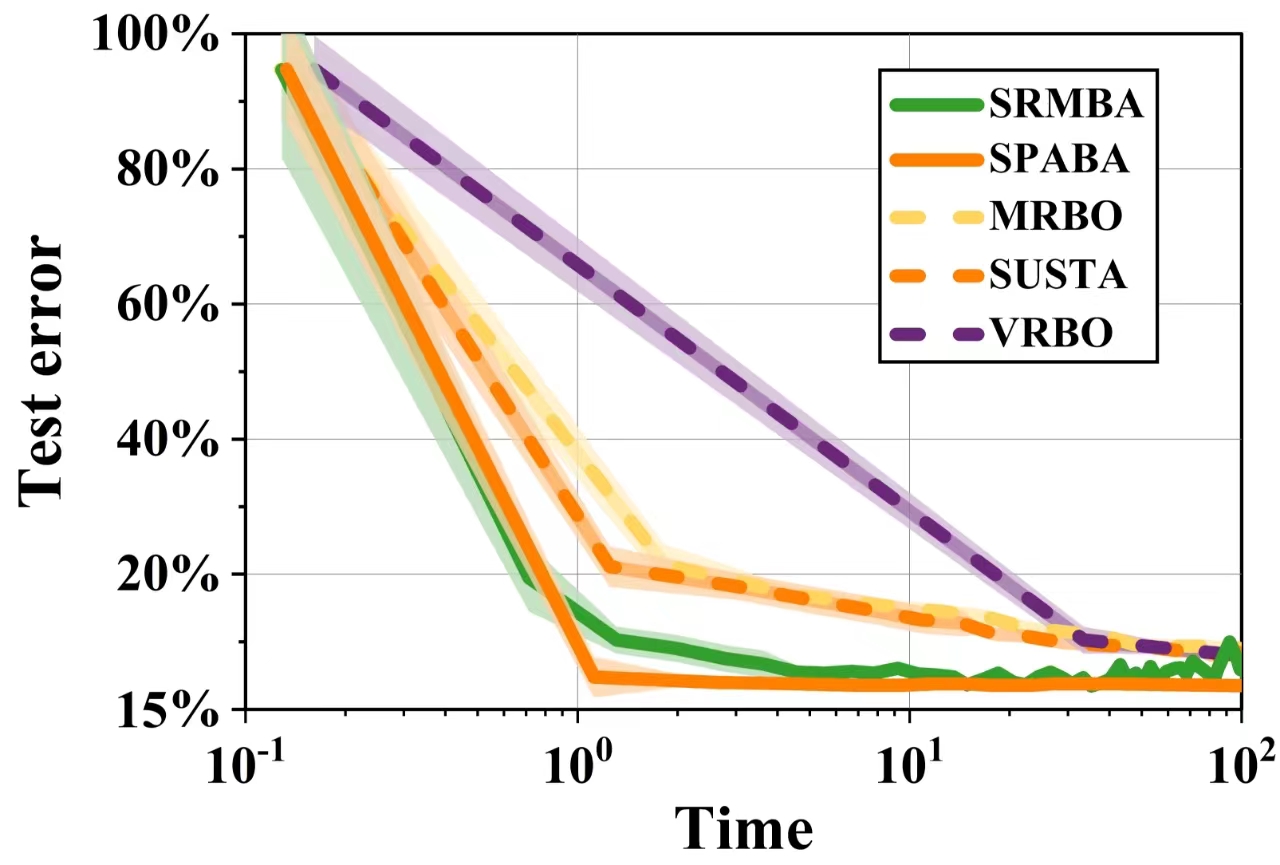}
		\label{datacleaning}
	}
	\subfigure[Logistic regression]{
		\centering
		\includegraphics[width=0.45\linewidth]{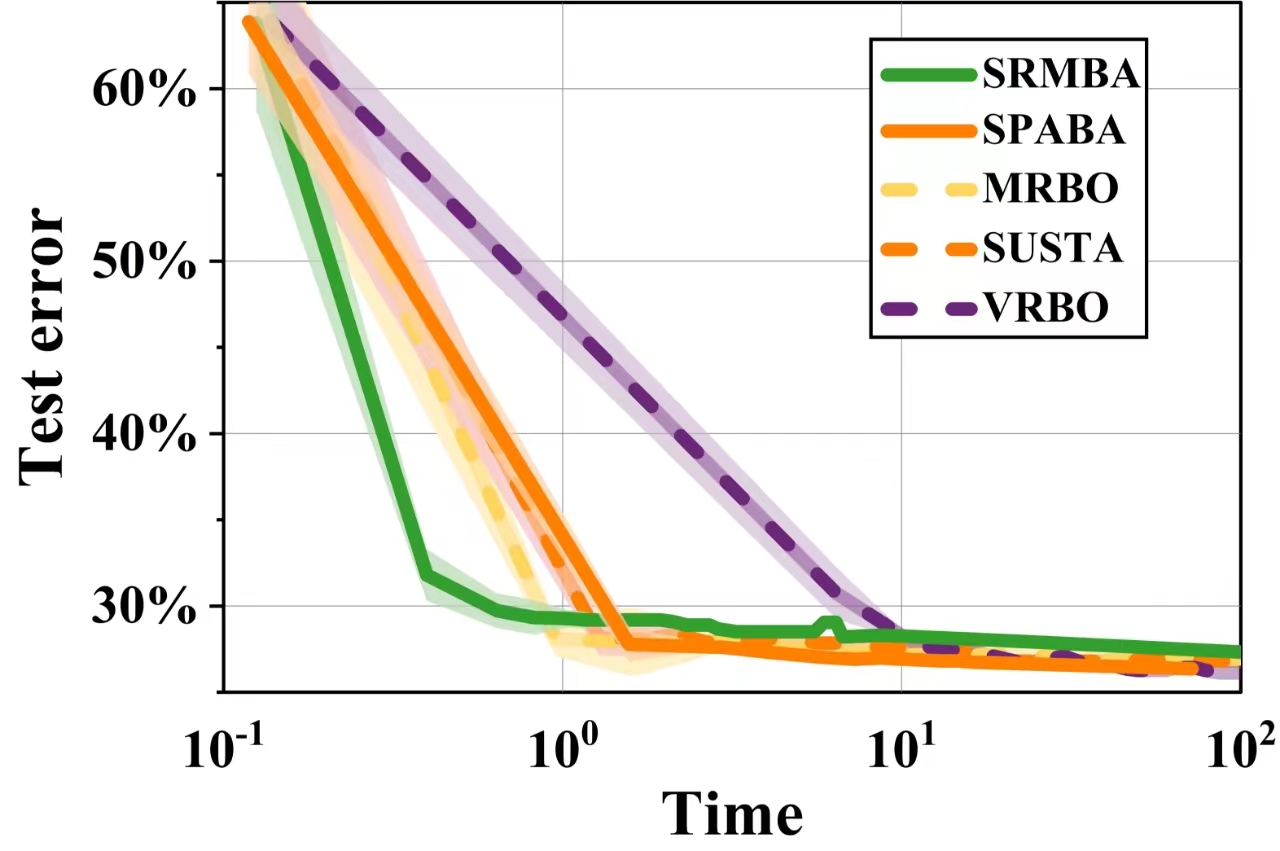}
		\label{covtype}
	}
	\caption{\textbf{Left:} Compare SRMBA and SPABA with other accelerated algorithms in a data hypercleaning experiment on the MINST dataset. \textbf{Right:} Compare SRMBA and SPABA with other accelerated algorithms in a hyperparameter selection experiment on the covtype dataset.
 }
	\label{sT}
\end{figure}

\begin{figure}
	\centering
	\includegraphics[width=5cm]{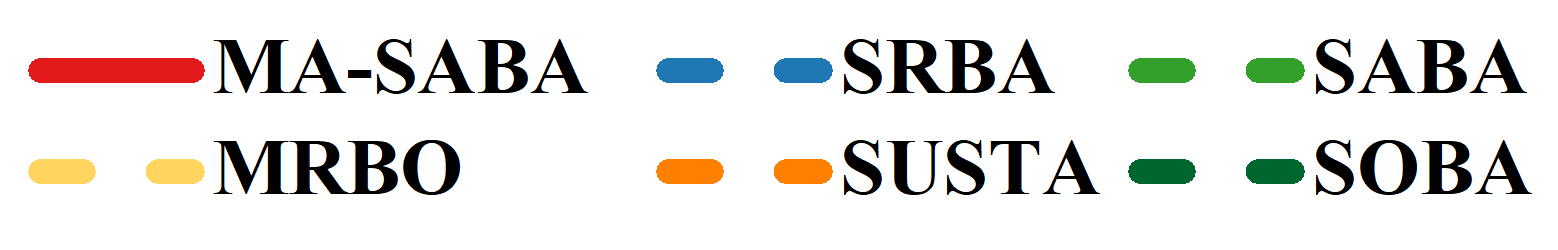}
	\subfigure{
		\includegraphics[width=0.45\linewidth]{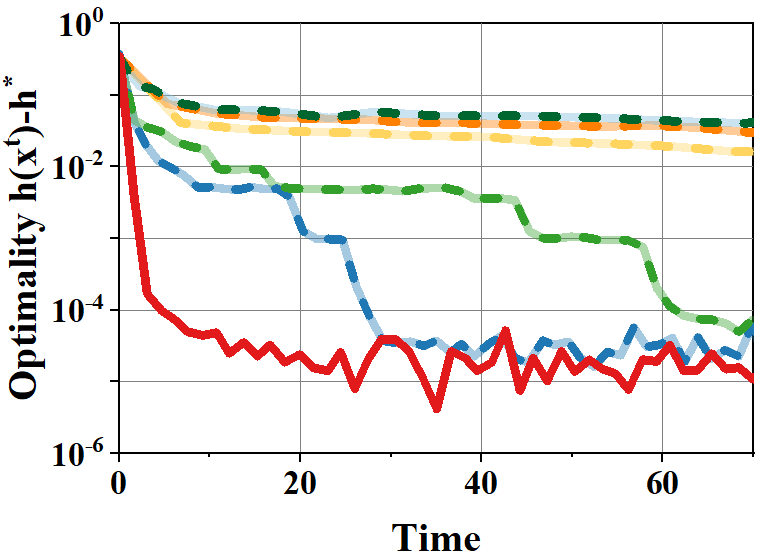}
	}
	\subfigure{
		\includegraphics[width=0.45\linewidth]{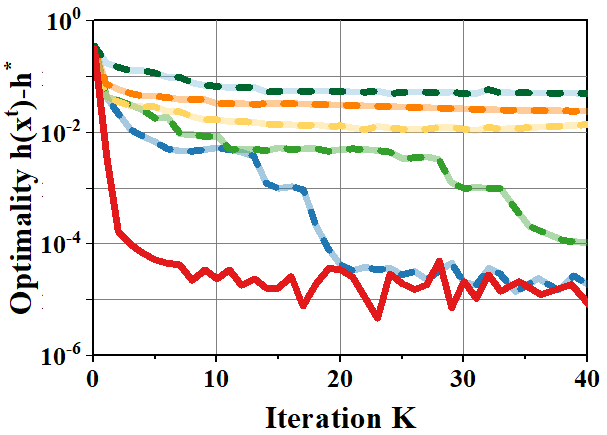}
	}
	\caption{Comparison of MA-SABA with competitors in a hyperparameter selection experiment. The results indicate that MA-SABA outperforms other methods in terms of both time and iteration.
 Solid lines depict our proposed methods, whereas dashed lines represent competitors.
 }
	\label{fig:enter-label-a}
\end{figure}

\subsection{Data Hyper-Cleaning}

The first learning task we perform is data hyper-cleaning conducted on the MNIST dataset \footnote{\url{http://yann.lecun.com/exdb/mnist/}} \cite{franceschi2017forward}. The dataset is divided into a training set $\left(d_i^{\text {train }}, y_i^{\text {train }}\right)$, a validation set $\left(d_j^{\text {val }}, y_j^{\text {val }}\right)$ and a test set.
The training set comprises 20,000 samples, the validation set contains 5,000 samples and the test set encompasses 10,000 samples. The target values $y$ range from 0 to 9, while the samples $d$ are of dimension 784. Within the training set, each sample is subject to corruption with a probability $\tilde{p}$: a sample $d_i$ is deemed corrupted when its label $y_i$ is replaced by a random label from the set $\{0, \ldots, 9\}$. Samples within the validation and test sets remain uncorrupted. The objective of data cleaning is to train a multinomial logistic regression model on the training set and ascertain a weight per training sample, ideally diminishing to 0 for corrupted samples. This is formalized by the bilevel optimization problem with $f(\lambda,\theta)=\frac{1}{m} \sum_{j=1}^m \ell\left(\theta d_j^{\text {val }}, y_j^{\text {val }}\right)$ and $g(\lambda,\theta)=\frac{1}{n} \sum_{i=1}^{n} \sigma\left(\lambda_i\right) \ell\left(\theta d_i^{\text {train }}, y_i^{\text {train }}\right)+C_r\|\theta\|^2$, where $\ell$ is the cross entropy loss and $\sigma$ is the sigmoid function. We report in Figure \ref{datacleaning} the test error, i.e., the percentage of incorrect predictions on the testing data. We utilize a corruption probability of $\tilde{p}=0.7$ (sample corruption rate) for this experiment. In this scenario, SPABA demonstrates the most favorable performance.

\subsection{Hyperparameter Selection}

We address hyperparameter selection for determining regularization parameters in $\ell^2$ logistic regression.
Let $\left(\left(d_i^{\text {train }}, y_i^{\text {train }}\right)\right){1 \leq i \leq n}$ and $\left(\left(d_j^{\text {val }}, y_j^{\text {val }}\right)\right){1 \leq j \leq m}$ denote the training and validation sets, respectively. In this context, the LL variable $\theta$ corresponds to the model parameters, while the UL variable $\lambda$ represents the regularization parameter. The functions $f$ and $g$ for bilevel optimization are defined as follows: $f(\lambda,\theta)=\frac{1}{m} \sum_{j=1}^m \varphi\left(y_j^{\mathrm{val}}\left\langle d_j^{\mathrm{val}}, \theta\right\rangle\right)$ and $g(\lambda,\theta)=\frac{1}{n} \sum_{i=1}^n \varphi\left(y_i^{\text {train }}\left\langle d_i^{\text {train }}, \theta\right\rangle\right)+\frac{1}{2} \sum_{k=1}^p e^{\lambda_k} \theta_k^2$ where $\varphi(u)=\log \left(1+e^{-u}\right)$.
In this experiment, two datasets, namely IJCNN1 and covtype, are employed, corresponding to the algorithms MA-SABA, SPABA and SRMBA, respectively. In Figure \ref{covtype}, the test error is presented alongside the corresponding running time. It is observed that SRMBA exhibits the shortest runtime, while SPABA achieves the highest accuracy promptly.
In the hyperparameter selection experiment, the suboptimality gap is depicted in Figure \ref{fig:enter-label-a} for each method. The lowest values are attained by MA-SABA, indicating its superior performance. MA-SABA reaches a considerably high final value, significantly outperforming other methods.

\section{Conclusion}

In this work we propose a loopless and sample-efficient stochastic bilevel algorithm, named SPABA, achieving optimal sample
complexity in both the finite-sum and expectation
settings.
Technically, SPABA is an adaptation of the PAGE algorithm in \cite{li2021page} within the proposed framework in \cite{arbel2022amortized, dagreou2022framework}.
More importantly, the complexity analysis of SPABA can be easily generalized to other stochastic gradient estimators. In fact, it already leads to MA-SABA and SRMBA that is an adaptation of STORM to the bilevel setting.
It's worth noting that the proposed algorithms still rely on computing Hessian and Jacobian matrices. Recent works \cite{chen2023near,yao2024constrained,kwon2023penalty} have used value function approaches to avoid querying second-order oracle information. Leveraging these approaches to develop single-loop, Hessian-free stochastic bilevel algorithms that achieve optimal or near-optimal sample complexity would be interesting and promising.

\section*{Impact Statements}
This paper presents work whose goal is to advance the field of Machine Learning. There are many potential societal consequences of our work, none which we feel must be specifically highlighted here.
\section*{Acknowledgements}
Authors listed in alphabetical order. This work is supported by National Key R \& D Program of China (2023YFA1011400), National Natural Science Foundation of China (12131003, 12222106, 12326605, 62331014, 12371305), Guangdong Basic and Applied Basic Research Foundation (No. 2022B1515020082) and Shenzhen Science and Technology Program (No. RCYX20200714114700072).
We thank Chengming Yu and all anonymous reviewers for their valuable comments and suggestions on this work.

\nocite{langley00}

\bibliography{sto_bilevel_ICML2014}
\bibliographystyle{icml2024}

\newpage
\appendix
\onecolumn

\section{Appendix}

The appendix is organized as follows:
\begin{itemize}
        \item We present a unified framework for converagence analysis and highlight the proof sketch of Theorems in Section \ref{proofsketch}.
	\item Additional experimental results are provided in Section \ref{sec:aer}.
        \item Algorithms and general lemmas are provided in Section \ref{alg&lemma}.
	\item
            Proof details are provided in Sections \ref{detial:proofmasaba} to \ref{detial:proofs trom}.

        \item The algorithm description and proof for MA-SOBA-q are provided in Section \ref{detial:proofsoba}.
\end{itemize}

\section{Convergence Analysis Framework and Proof Sketches for Theorems}\label{proofsketch}

To analyze complexity, we introduce a general and unified convergence analysis method. In this section, we provide a more detailed exposition. Furthermore, we illustrate this by proving Theorems in this paper.

\subsection{Convergence Analysis Framework}\label{anaframedetial}

Our convergence analysis relies on the following recursive inequality
\begin{align}\label{eana}
    \alpha_k\mathbb{E}\left[\left\|\nabla H(x_k)\right\|^2\right]\leq L_k-L_{k+1}+\Delta_k,
\end{align}
where $L_k$, $L_{k+1}$, and $\Delta_k$ are all positive terms,
$\alpha_k$ is the step size used for updating the UL variable $x_k$.
The term $L_k$ is referred to as the potential function or Lyapunov function.

To derive (\ref{eana}), there are two crucial considerations: first, how to construct an appropriate $L_k$; and second, analyzing the descent of each element within $L_k$.
Although the Lyapunov does not possess a uniform form, we start from the descent of the total UL objective and analyze layer-by-layer the elements it should comprise and their respective descents. This will be presented in \ref{3step}.

Assuming that, through the analysis of the decreasing properties of the elements in Lyapunov function and the selection of appropriate parameters, we have obtained (\ref{eana}). The next customary step is to define $\theta=\min_{k\in[K]}\{\alpha_k\}$, and by induction, we have
\begin{align*}
    \frac{1}{K}\sum_{k=0}^{K-1}\mathbb{E}\left[\left\|\nabla H(x_k)\right\|^2\right]\leq \frac{L_0}{K\theta}+\frac{\sum_{k=0}^{K-1}\Delta_k}{K\theta},
\end{align*}
which characterizes the convergence of the algorithm.
\subsection{Proof Sketches for Theorems }\label{3step}

In this section, we will present proof sketches for the theorems, utilizing a four-step layer-by-layer analysis to derive a recursive inequality similar to (\ref{eana}). Additionally,
\textbf{this analysis showcases how to close the gap between stochastic bilevel and single-level optimization under classical assumptions and how to effectively handle biased stochastic estimations to attain superior complexity results. }

\begin{figure}[h]
    \centering
    \includegraphics[width=1\linewidth]{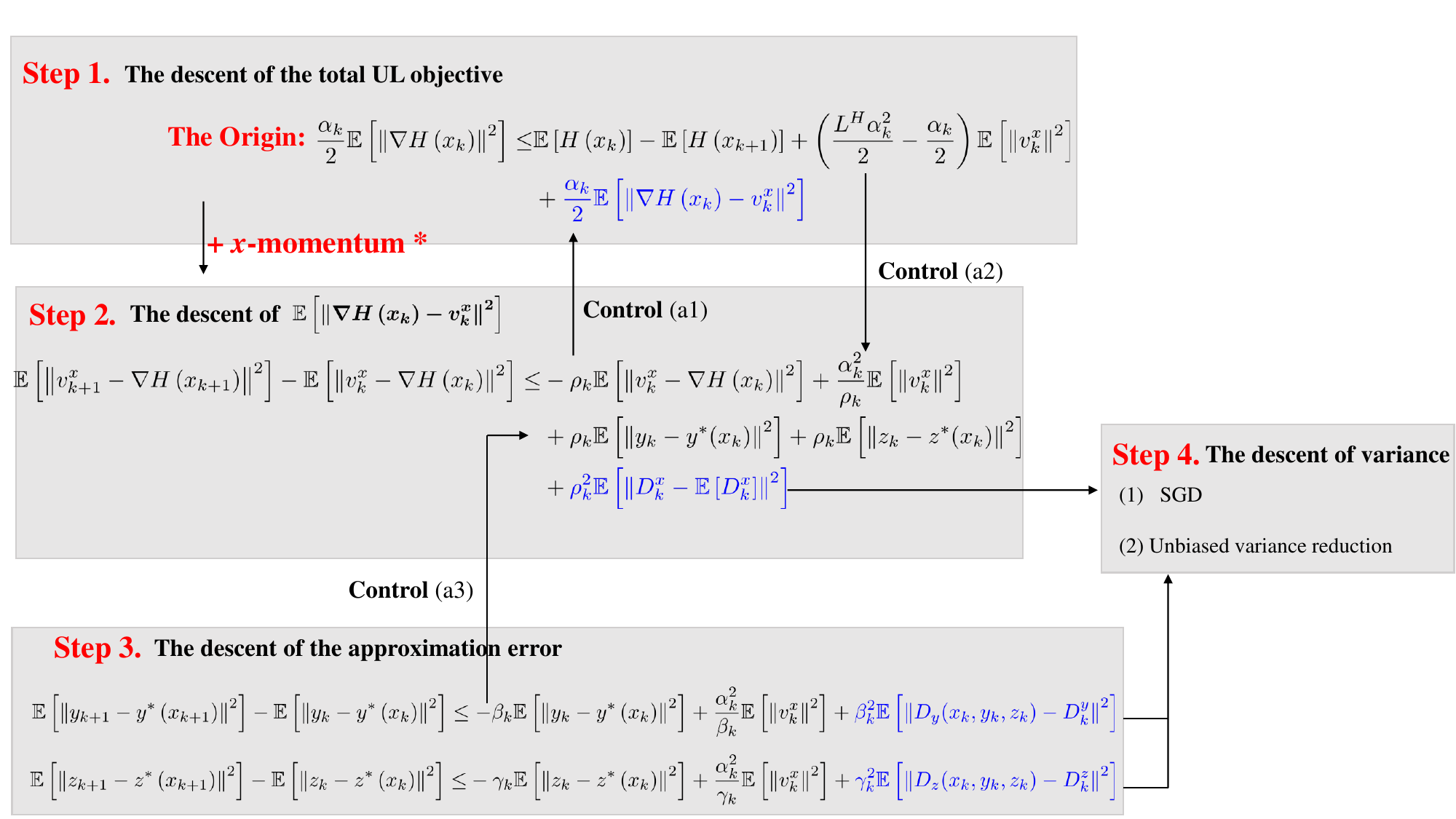}
    \caption{Proof sketch (using unbiased estimations for the iteration directions of $y$ and $z$)\\
    (*):
Fortunately, with the momentum applied to $x$, we achieve the second step, namely the descent of
$\mathbb{E}[\|\nabla H(x_k)-v_k^x\|^2]$.
Control(a1) demonstrates how Step 2 manages $\alpha_k\mathbb{E}[\|\nabla H(x_k)-v_k^x\|^2]$.
Control(a2) and Control(a3) illustrate that step 2 does not introduce new, uncontrollable terms.
Each term can be managed by inequalities found in either Step 1 or 2.
The blue section in the figure indicating Step 4 highlights the variance terms critically influencing the convergence rate and complexity. This necessitates further examination in Step 4, employing either SGD or its variants with variance reduction.
    }
    \label{fig:proofsketchunbaised}
\end{figure}

\begin{figure}[h]
    \centering
    \includegraphics[width=1\linewidth]{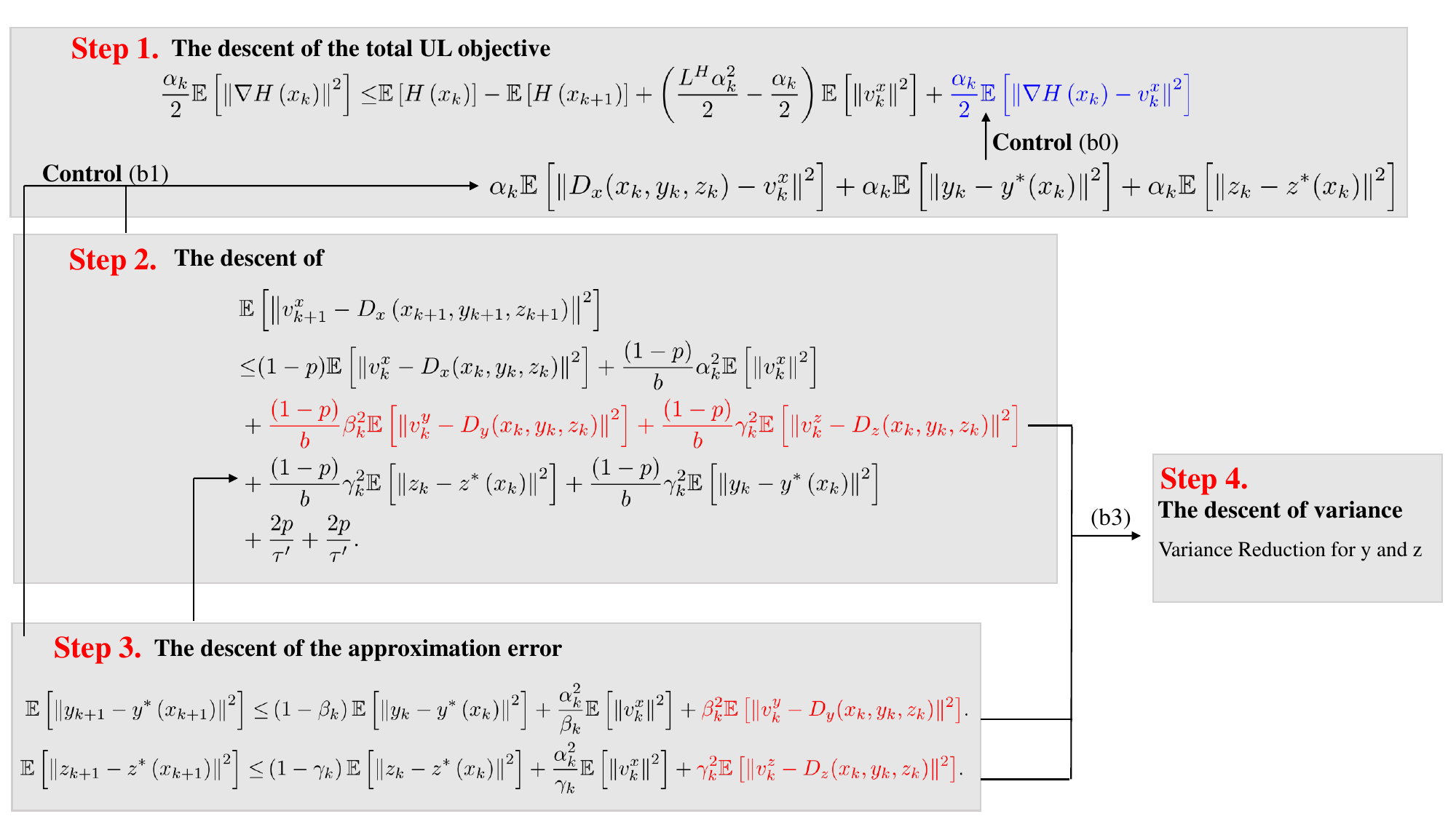}
    \caption{Proof sketch (using biased estimations for the iteration directions of $y$ and $z$).
In this analytical framework, we begin by dissecting $\mathbb{E}[\|\nabla H(x_k)-v_k^x\|^2]$ into two segments,
$\mathbb{E}[\|D_x(x_k,y_k,z_k)-v_k^x\|^2]$ and $\mathbb{E}[\|y_k-y^*(x_k)\|^2]+\mathbb{E}[\|z_k-z^*(x_k)\|^2$, utilizing
Control(b0). These segments are subsequently regulated by Step 2 and Step 3, referred to as
Control(b1). The unique aspect here is that the inequalities applied in steps 2 and 3 are specifically designed for biased estimations. This adaptation enables the integration of an expanded selection of variance reduction methods to efficiently manage the red section.}
    \label{fig:proofsketchbaised}
\end{figure}

Usually, the recursive estimate (\ref{ana}) is derived through a series of recursive inequalities in conditional expectation:
\begin{equation}\label{recursiveinequality}
    \mathbb{E}\big[\widetilde{D}_{k+1}\,|\, \mathcal{F}_k\big] + \Lambda_k
    \leq \omega_k \widetilde{D}_{k} + \Omega_k,
\end{equation}
where $\widetilde{D}_{k}$, $\Lambda_k$, $\Omega_k$ are all nonnegative quantities, and $\omega_k\in[0,1]$ is a contraction factor.
We can now divide the proof of the recursive estimate (\ref{ana}) into four main steps:

\textbf{Step 1: Originating from the descent of the total UL objective.}

We begin by bounding the descent of $H(x)$ as follows:
\begin{equation}
\begin{aligned}
\underbrace{\mathbb{E}\left[H\left(x_{k+1}\right)\right]}_{\mathbb{E}\big[\widetilde{D}_{k+1}\,|\, \mathcal{F}_k\big]}
+\underbrace{\frac{\alpha_k}{2} \mathbb{E}\big[\left\|\nabla H\left(x_k\right)\right\|^2\big]}_{\Lambda_k }
\leq
\underbrace{\mathbb{E}\left[H\left(x_k\right)\right]}_{\widetilde{D}_{k}}
+
 \underbrace{\Big(\frac{L^H\alpha_k^2}{2}-\frac{\alpha_k}{2}\Big) \mathbb{E}\big[\left\|v_k^x\right\|^2\big]
 +\frac{\alpha_k}{2} \mathbb{E}\big[\left\|\nabla H\left(x_k\right)-v_k^x\right\|^2\big]}_{\Omega_k^{(1)}},
\end{aligned}
\label{HinHOW}
\end{equation}
which is a recursive inequality with $\omega_k=1$, as demonstrated in (\ref{recursiveinequality}). It is established in Lemma \ref{Hlemma} by the $L^{H}$-smoothness of $H(x)$.
With this lemma, we can achieve the same handling of the effect of
$\alpha_k^2/\beta_k$ as in Lemma 3.9 in \cite{dagreou2022framework}, even though we introduce a new term $\mathbb{E}\left[\left\|\nabla H\left(x_k\right)-v_k^x\right\|^2\right]$. We will address this term in the next step.

\textbf{Step 2: The descent of the mentioned mean-squared error.}

Considering the presence of a mean-squared error $\mathbb{E}\big[\left\|\nabla H\left(x_k\right)-v_k^x\right\|^2\big]$ in $\Omega_k^{(1)}$ in (\ref{HinHOW}), we study the descent of the mentioned mean-squared error.
When integrating a standard momentum or a variation of momentum, such as those found in PAGE and STORM, into the update of $x_k$, we can establish a recursive inequality in the form of (\ref{recursiveinequality}) for $\widetilde{D}_{k}:=\left\|\nabla H\left(x_k\right)-v_k^x\right\|^2$. This inequality is derived from two or three consecutive iterates.

For example, when standard momentum is utilized in $v_k^x$, we derive the following recursive inequality:
\begin{equation}
\begin{aligned}\label{step2}
&\underbrace{\mathbb{E}\left[\left\|v_{k+1}^x-\nabla H\left(x_{k+1}\right)\right\|^2\right]}_{\mathbb{E}\big[\widetilde{D}_{k+1}\,|\, \mathcal{F}_k\big]}
\\
\leq&
\underbrace{\left(1-\rho_k\right)}_{\omega_k} \underbrace{\mathbb{E}\left[\left\|v_k^x-\nabla H\left(x_k\right)\right\|^2\right]}_{\widetilde{D}_{k}}
\\&+\underbrace{\frac{2\left(L^H\right)^2 \alpha_k^2}{\rho_k} \mathbb{E}\left[\left\|v_k^x\right\|^2\right]
+2\rho_k\mathbb{E}\left[\left\|D_x(x_k,y_k,z_k)-\nabla H(x_k)\right\|^2\right]
+\rho_k^2 \mathbb{E}\left[\left\|D_k^x-\mathbb{E}\left[D_k^x\right]\right\|^2\right]}_{\Omega_k^{(2)}},
\end{aligned}
\end{equation}
a result that is rigorously proven in Lemma \ref{vh}.

For Theorems \ref{thpagefini}, \ref{thpage1.5} and \ref{thstorm}, they utilize the momentum variants PAGE and STORM, which are based on the fact that
$$\mathbb{E}\left[\left\|v_k^x-\nabla H\left(x_k\right)\right\|^2\right]\leq
2\mathbb{E}\left[\left\|v_k^x-D_x(x_k,y_k,z_k)\right\|^2\right]
+2\mathbb{E}\left[\left\|D_x(x_k,y_k,z_k)-\nabla H\left(x_k\right)\right\|^2\right],$$
and we can obtain recursive inequalities about $\mathbb{E}\left[\left\|v_k^x-D_x(x_k,y_k,z_k)\right\|^2\right]$ similar to (\ref{step2}), respectively seen in Lemmas \ref{pagefinitemoment}(2), \ref{pagegeneralmonment}(2), \ref{eg33}(2).

It's important to highlight that the contraction factor $\omega_k=1-\rho_k$, where $\rho_k$ is the ``momentum" parameter and will tend to approach $0$ in the subsequent setting. If there is no momentum, it is only possible to obtain an upper bound for the mean-squared error.

\textbf{Step 3: The descent of the approximation error.}

To gain better control over the terms in $\Omega_k^{(2)}$ in the descent of the mean-squared error $\mathbb{E}\big[\left\|\nabla H\left(x_k\right)-v_k^x\right\|^2\big]$,it is essential to investigate the descent of $\mathbb{E}\left[\left\|D_x(x_k,y_k,z_k)-\nabla H\left(x_k\right)\right\|^2\right]$. Leveraging the foundational fact presented in Lemma \ref{DH}, that
\begin{eqnarray*}
\mathbb{E}\left[\left\|D_x(x_k,y_k,z_k)-\nabla H\left(x_k\right)\right\|^2\right]
& \leq c_1 \mathbb{E}\left[\left\|y_k-y^*\left(x_k\right)\right\|^2\right] +c_2 \mathbb{E}\left[\left\|z_k-z^*\left(x_k\right)\right\|^2\right],
\end{eqnarray*}
our analysis extends to a thorough examination of the descent of the approximation errors $\mathbb{E}\big[\left\|y_k-y^*\left(x_k\right)\right\|^2\big]$ and $\mathbb{E}\big[\left\|z_k-z^*\left(x_k\right)\right\|^2\big]$.
By leveraging the strongly convexity of the LL problem and the quadratic problem (\ref{z*}), one can readily derive recursive inequalities akin to (\ref{recursiveinequality}) for the approximate errors.
Additionally, $(1-\omega_k)$ exhibits a similar order of magnitude as the step sizes for both $y$ and $z$; please refer to Lemma \ref{y,z} and Lemma \ref{eg30} for illustrations.

\textbf{Step 4: The descent of variance.}

In all the recursive inequalities mentioned above, the remaining terms include only
the variances of the stochastic gradient estimators, such as $\mathbb{E}\big[\left\|D_x(x_k,y_k,z_k)-D_k^x\right\|^2\big]$.
If the stochastic gradient estimators used lack variance reduction properties, like SGD, it is only feasible to attain a constant upper bound, even when we consider Assumption \ref{assumporacle}.
To further reduce sampling complexity, one can integrate unbiased or biased variance reduction techniques into the algorithm.
MA-SABA aligns its sampling complexity with that of single-level optimization using SAGA. For an illustration, please refer to Lemma \ref{var}.
SPABA implements the variance reduction stochastic estimation technique PAGE. The recursive inequalities for variance reduction needed for Theorems \ref{thpagefini} and \ref{thpage1.5} are detailed in Lemmas \ref{pagefinitemoment}(1)(3) and \ref{pagegeneralmonment}(1)(3).
SRMBA incorporates the STORM technique, with the recursive inequalities for variance reduction elucidated in Lemma \ref{eg33}(1)(3).

Table \ref{tab:lemmath} outlines the specific lemmas utilized in the four steps integral to the proof of each theorem.
Utilizing the delineated four-step framework, we craft a Lyapunov function embedded with essential variables, judiciously select coefficients for this function, and calibrate the algorithm's step size parameters. Through this strategic approach, we establish inequalities parallel to equation (\ref{ana}), thus paving the way for substantiating convergence outcomes.

\begin{table}[H]
    \centering
    \begin{tabular}{|c|c|c|c|c|}
    \hline
        Methods and Conclusions  & Step 1 & Step 2 & Step 3 & Step 4 \\
        \hline
        MA-SABA (Th\ref{thsaba}) & Lemma \ref{Hlemma} &  Lemma \ref{vh} & Lemma \ref{y,z}&Lemma \ref{S} \\
        \hline
        SPABA (Th\ref{thpagefini}) &  Lemma \ref{Hlemma} &  Lemma \ref{pagefinitemoment}(2) & Lemma \ref{eg30} & Lemma \ref{pagefinitemoment}(1)(3) \\
         \hline
         SPABA (Th\ref{thpage1.5})&  Lemma \ref{Hlemma} &  Lemma \ref{pagegeneralmonment}(2) & Lemma \ref{eg30} & Lemma \ref{pagegeneralmonment}(1)(3) \\
        \hline
        SRMBA (Th\ref{thstorm})&  Lemma \ref{Hlemma} &  Lemma \ref{eg33}(2) &Lemma \ref{eg30}& Lemma \ref{eg33}(1)(3)\\
        \hline
    \end{tabular}
    \caption{Lemmas Aligned with Each Step in Theorem Proofs}
    \label{tab:lemmath}
\end{table}
\textbf{Comparison with \cite{dagreou2022framework}}

In \cite{dagreou2022framework}, without stronger smoothness conditions, the derivation was limited to  Lemma D.3, presenting a challenge as the coefficient of $\mathbb{E}[\|v_k^x\|]^2$ changed from $\alpha_k^2$ to
$\alpha_k^2/\beta_k$.
They pointed out that to achieve convergence, it is required that the ratio
$\alpha_k/\beta_k$ goes to zero, stating, ``This prevents us from getting rates that match rates of single-level algorithms."

Our approach uniquely addresses the challenging term $\mathbb{E}[\|v_k^x\|]^2\alpha_k^2/\beta_k$ from a new perspective.
To manage this term, we describe the descent of $H(x_k)$ through (\ref{HinHOW}), providing a characterization different from Lemma 3.10 in \cite{dagreou2022framework}, albeit introducing an additional new term $\mathbb{E}[|\nabla H\left(x_k\right)-v_k^x|^2]$.
Excitingly, by incorporating momentum into the iterative direction of $x_k$, we ensure a decrease in this term while preventing the emergence of new terms with coefficients analogous to $\mathbb{E}[|v_k^x|]^2$, thus circumventing the limitations mentioned in \cite{dagreou2022framework}.

\textbf{Extra: Utilizing Momentum-Based Biased Variance Reduction.}

Our framework is meticulously designed to adeptly address biased estimations, recognizing the efficacy of targeted variance reduction strategies in yielding superior results. A pivotal distinction of our approach is the specialized adaptation of Step 2 and Step 3, meticulously crafted to accommodate biased estimations. By harnessing the capabilities of this framework, we unlock the potential to develop stochastic algorithms specifically engineered for bilevel optimization, thereby achieving markedly lower sampling complexities.

Therefore, by selecting appropriate step sizes and coefficients for the Lyapunov functions to scale the inequality, we can derive a recursive inequality of the form similar to (\ref{ana}).

\newpage
\section{Additional experimental results}\label{sec:aer}
All experiments were conducted in Python, utilizing the Benchopt package \cite{moreau2022benchopt}, JAX \cite{jax2018github}, and Numba \cite{lam2015numba} for efficient implementation of stochastic methods. For each problem, oracles for a given function $f$ were employed, providing the quantities $\left(f(x,y), \nabla_1 f(x,y), \nabla_{22}^2 f(x,y) z, \nabla_{12}^2 f(x,y) z\right)$ to avoid redundant computation of intermediate results.

The experiments were executed using Python 3.8 on a system equipped with an Intel(R) Xeon(R) Gold 5218R CPU @ 2.10GHz and an NVIDIA A100 GPU with 40GB of memory.

\subsection{Hyperparameter selection on covtype dataset}

Similar to \cite{dagreou2022framework}, we conducted an additional experiment involving the selection of the best regularization parameter for an $\ell^2$-regularized multinomial logistic regression problem on the covtype dataset\footnote{\url{https://scikit-learn.org/stable/modules/generated/sklearn.datasets.fetch_covtype.html}}. This dataset comprises 581,012 samples with $p=54$ features and encompasses $C=7$ classes. Specifically, we utilized $n=371,847$ training samples, $m=92,962$ validation samples, and $n_{\text {test }}=116,203$ test samples.

In this experiment, we fitted a multiclass logistic regression model on this dataset, with one hyperparameter per class. Thus, if $\left(d_i^{\text {train }}, y_i^{\text {train }}\right)_{i \in[n]}$ and $\left(d^{\mathrm{val}}_j, y^{\mathrm{val}}_j\right)_{j \in[m]}$ represent the training and validation samples, respectively, we solve the following bilevel optimization:
$$
\begin{aligned}
	& f(\lambda,\theta)=\frac{1}{m} \sum_{j=1}^m \ell\left(\theta d_j^{\mathrm{val}}, y_j^{\mathrm{val}}\right) \text { and } \\
	& g(\lambda,\theta)=\frac{1}{n} \sum_{i=1}^n \ell\left(\theta d_i^{\mathrm{train}}, y_i^{\mathrm{train}}\right)+\sum_{c=1}^C e^{\lambda_c} \sum_{i=1}^p \theta_{i, c}^2,
\end{aligned}
$$
where $\lambda \in \mathbb{R}^C$ is the UL variable and $\theta \in \mathbb{R}^{p \times C}$ is the LL variable.

\textbf{Hyper-parameter setting for algorithm.} For SPABA, the probability $p=0.5$, the step-sizes are chosen as $\alpha_k = 0.2/0.01, \gamma_k = \beta_k = 0.2$. For MA-SABA, the step-sizes are chosen as $\alpha_k = 0.2, \beta_k = 0.2/0.0001, \gamma_k = \beta_k$ and $ \rho_k = 0.2$. For SRMBA, the step-sizes are chosen as $\alpha_k = \frac{5}{k^{1/3}}, \beta_k  = \frac{0.2}{k^{1/3}}, \gamma_k = \frac{0.002}{k^{1/3}}$ and $ \rho_k^{x} = \rho_k^{y} = \rho_k^{z} = \frac{0.5}{k^{2/3}} $. Other algorithms choose their step sizes according to the optimal strategy in \cite{dagreou2022framework}.

\begin{figure}[H]
	\centering
	\subfigure{
		\includegraphics[width=0.3\linewidth]{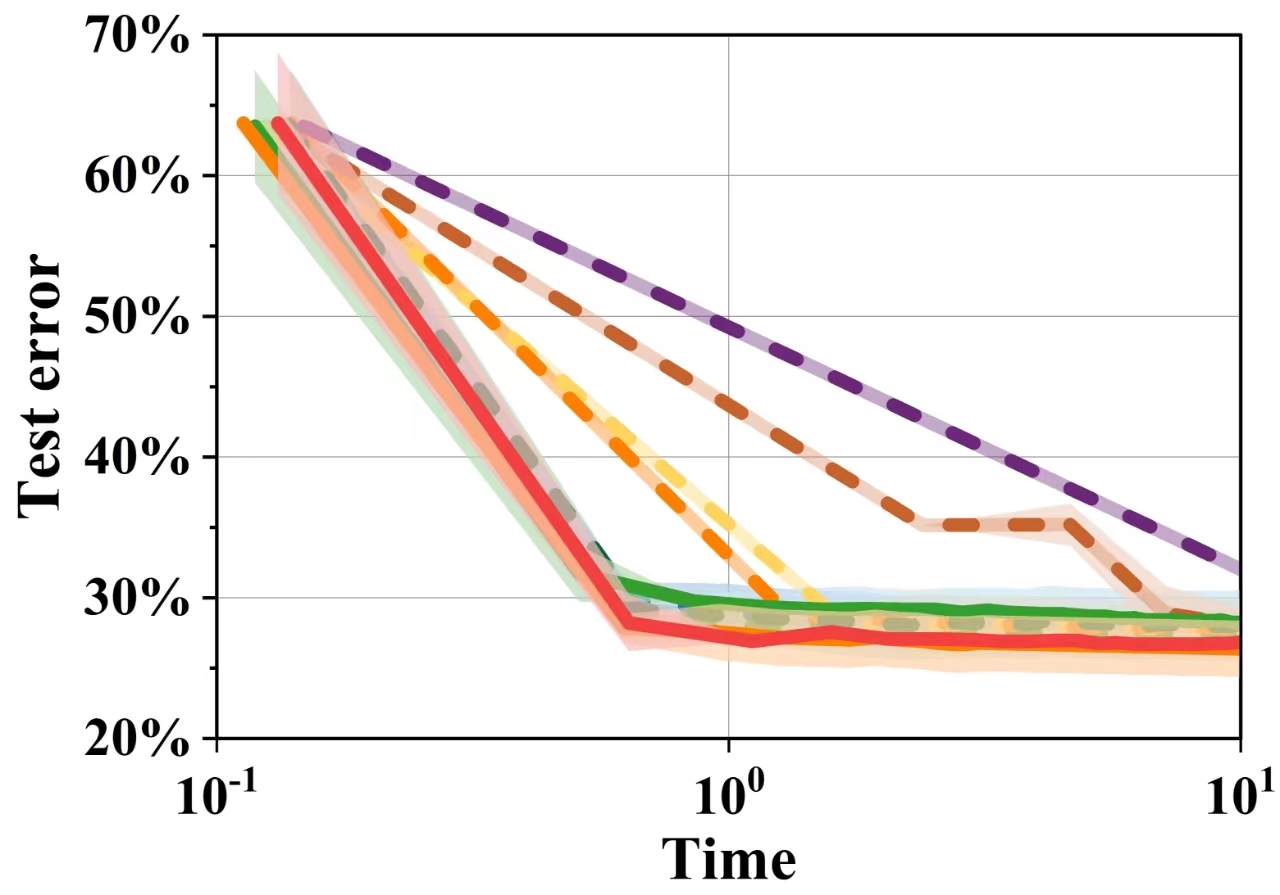}
	}
	\subfigure{
		\includegraphics[width=0.3\linewidth]{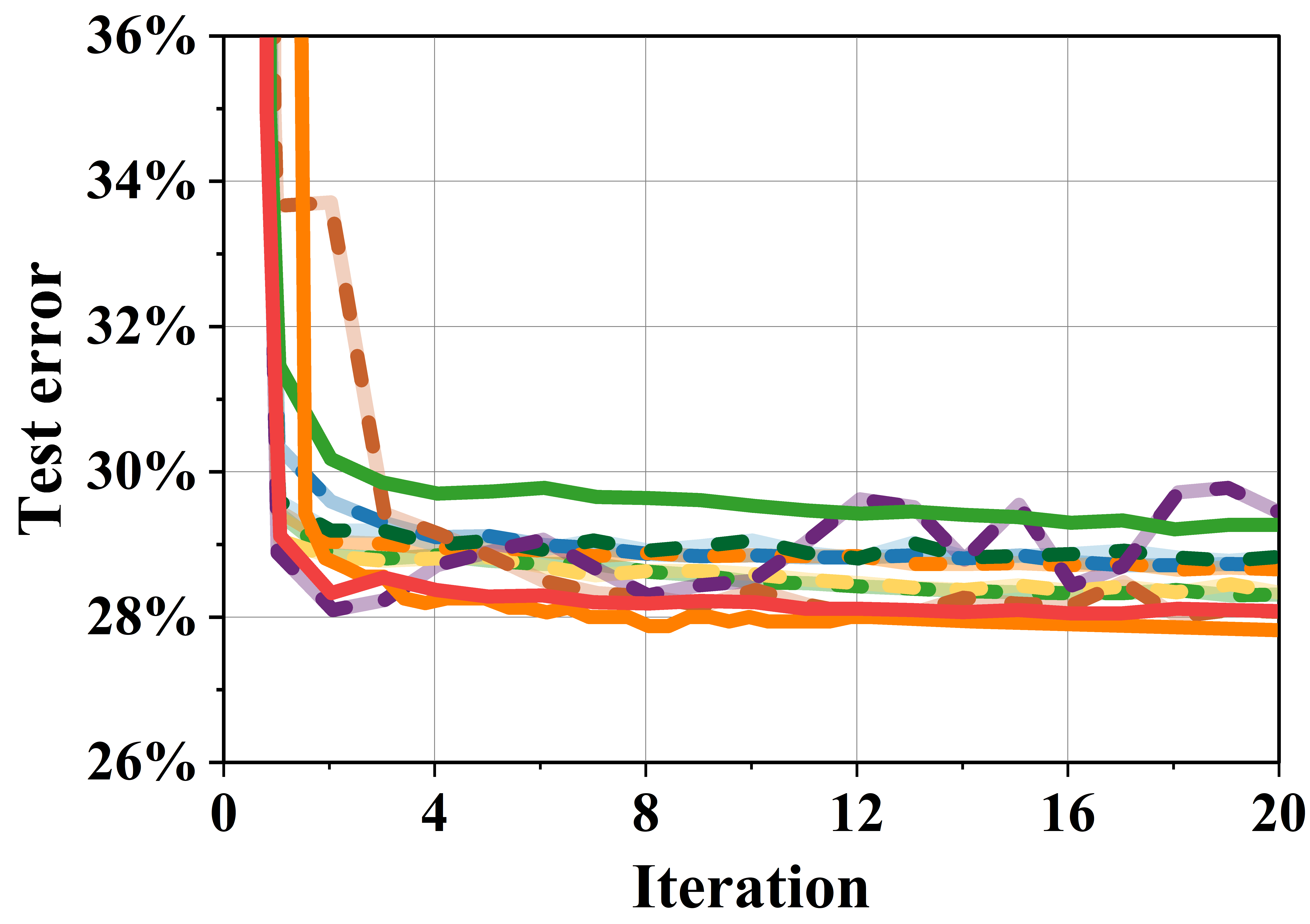}
	}
	\subfigure{
	\includegraphics[width=0.2\linewidth]{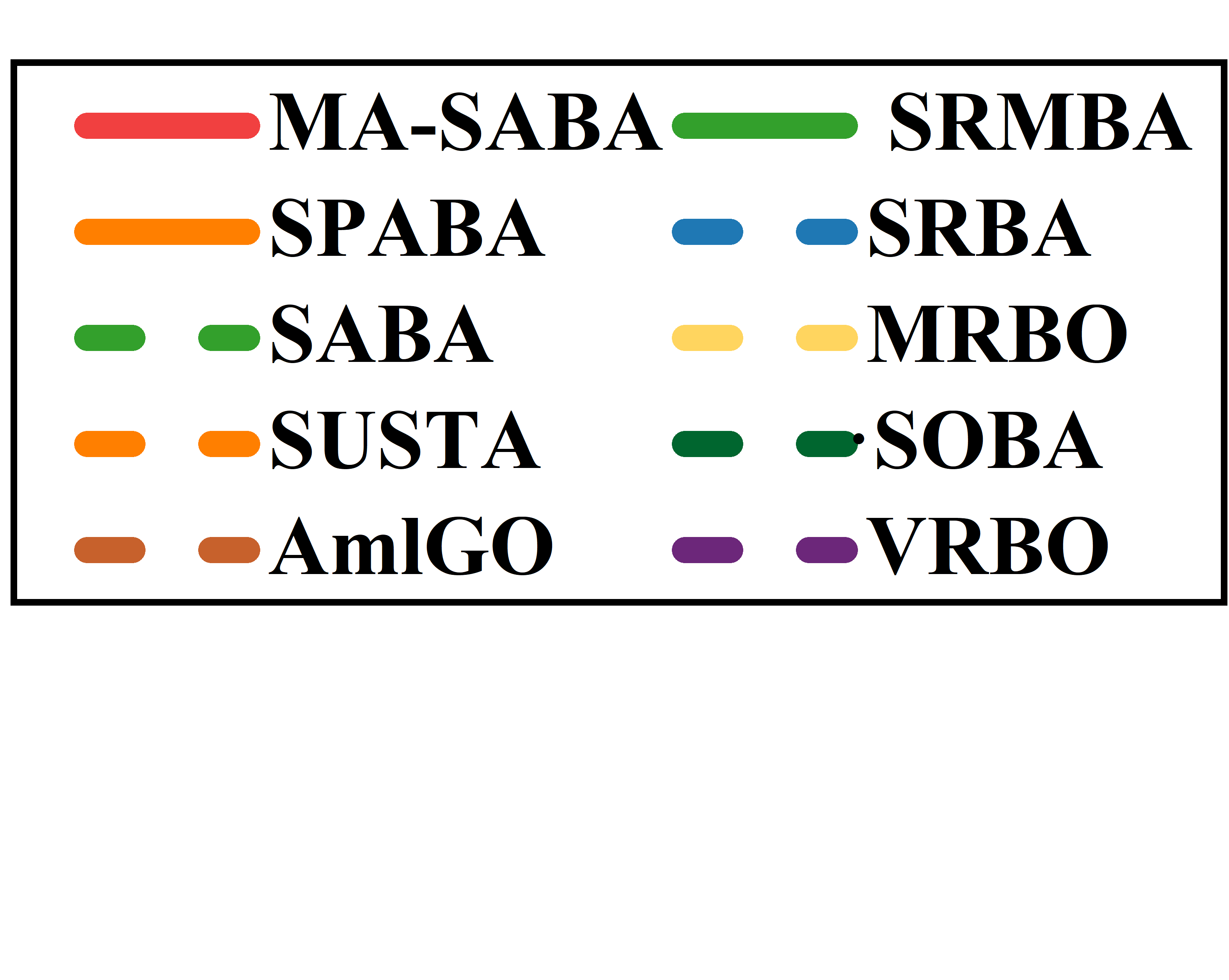}
    }
	\caption{Comparison of MA-SABA, SPABA, and SRMBA with other stochastic bilevel optimization methods in a hyperparameter selection experiment. The result reveals that MA-SABA achieves the best performance in terms of both time and iteration. The dashed lines represent other stochastic bilevel optimization methods, while the solid lines denote the proposed methods.}
	\label{fig:enter-label}
\end{figure}

\subsection{Hyperparameters selection on IJCNN1}

In this experiment, we select the regularization parameters for a multiregularized logistic regression model, where we have one hyperparameter per feature:
$$
\begin{aligned}
	& f(\lambda, \theta)=\frac{1}{m} \sum_{i=1}^m \varphi\left(y_i^{\mathrm{val}}\left\langle d_i^{\mathrm{val}}, \theta\right\rangle\right) \text { and } \\
	& g(\lambda, \theta)=\frac{1}{n} \sum_{i=1}^n \varphi\left(y_i^{\text {train }}\left\langle d_i^{\text {train }}, \theta\right\rangle\right)+\frac{1}{2} \theta^{\top} \operatorname{diag}\left(e^{\lambda_1}, \ldots, e^{\lambda_p}\right) \theta,
\end{aligned}
$$
where $\lambda,\theta$ are the UL and LL variables, respectively. The parametrization choice, using $e^\lambda$ rather than $\lambda$, ensures that there are no constraints placed on the variable $\lambda$.  It is a classical approach in the bilevel optimization literature \cite{pedregosa2016hyperparameter, ji2021bilevel, grazzi2021convergence}.

In these experiments, as in \cite{dagreou2022framework}, we employ Just-In-Time (JIT) compilation using the Numba package \cite{lam2015numba} to reduce Python overhead in the iteration loop. Additionally, to evaluate $H(\lambda)$, we utilize L-BFGS \cite{liu1989limited} to compute $y^*(x_k)$ and subsequently evaluate the function $H(x_k)=f(x_k, y^*(x_k))$.

\textbf{Hyper-parameter setting for algorithm.} For SPABA, the probability $p=0.5$, the step-sizes are chosen as $\alpha_k = 0.2/0.01, \gamma_k = \beta_k = 0.2$. For MA-SABA, the step-sizes are chosen as $\alpha_k = 0.5, \beta_k = 0.5, \gamma_k = 0.4$ and $ \rho_k = 0.2$. Other algorithms choose their step sizes according to the optimal strategy in \cite{dagreou2022framework}.

\begin{figure}[H]
	\centering
	\subfigure{
		\includegraphics[width=0.3\linewidth]{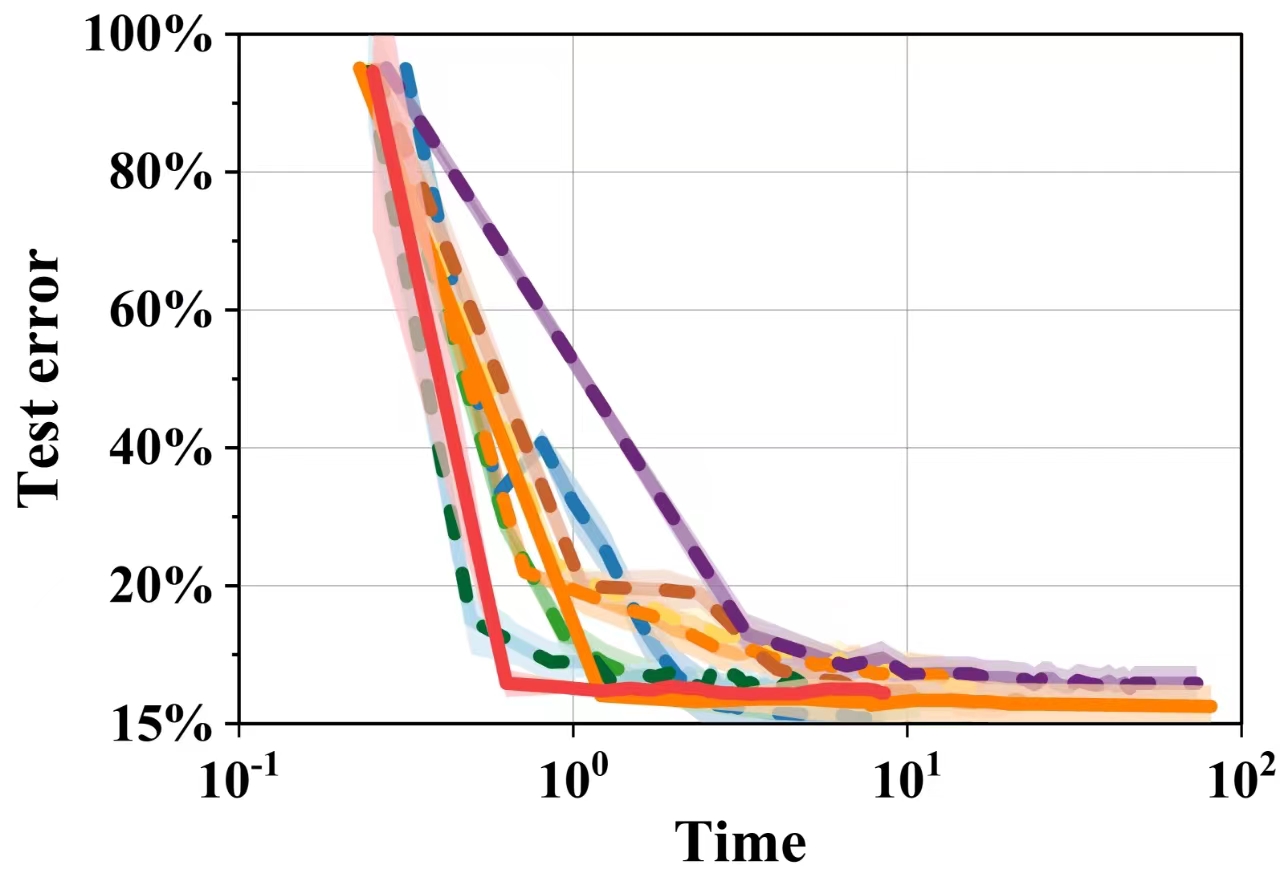}
	}
	\subfigure{
		\includegraphics[width=0.29\linewidth]{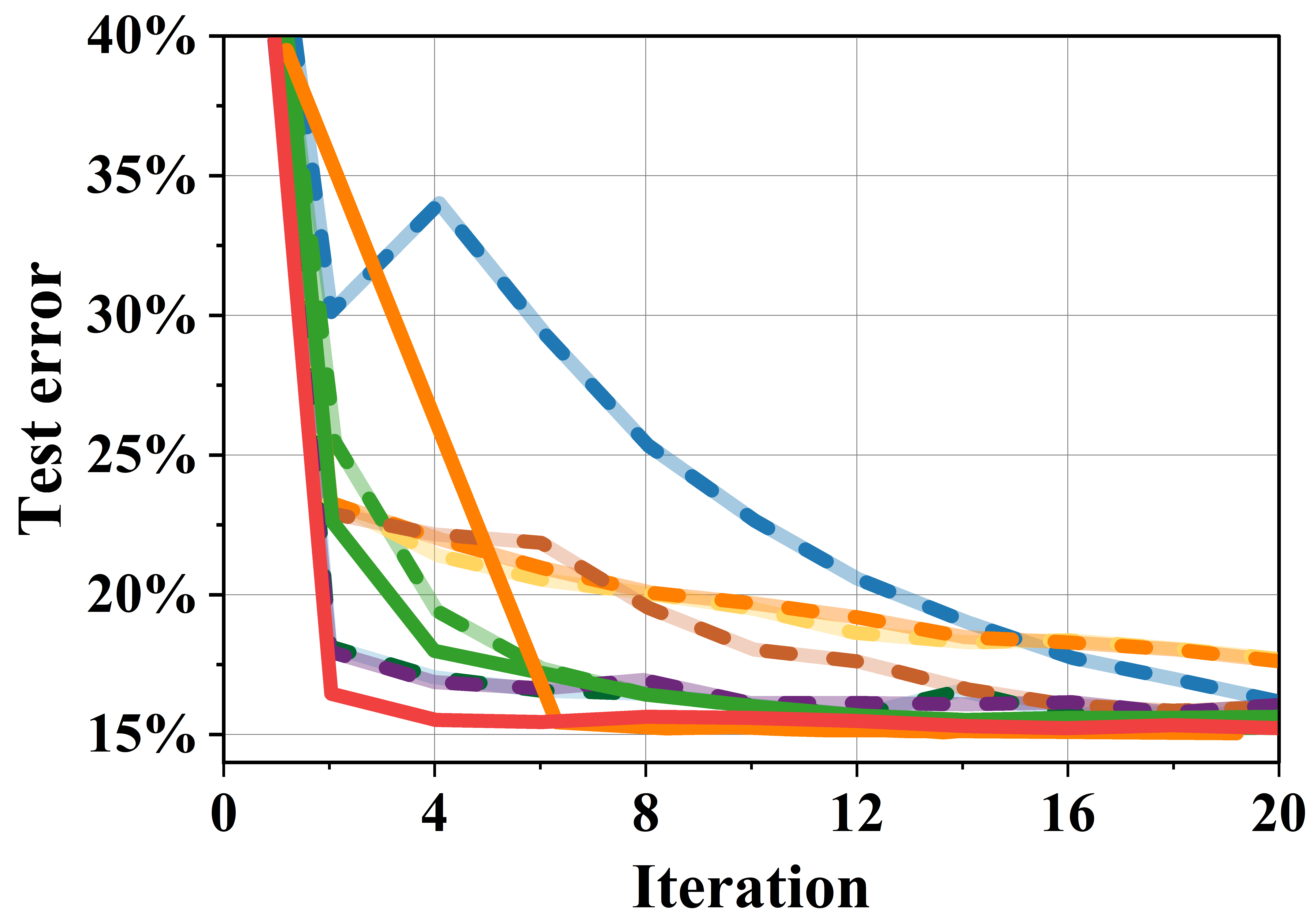}
	}
    \subfigure{
	\includegraphics[width=0.2\linewidth]{1.png}
    }
	\caption{
		Comparison of MA-SABA, SPABA, and SRMBA with other stochastic bilevel optimization methods in a data hyper-cleaning experiment. It demonstrates that MA-SABA achieves superior performance in both time and iteration. The dashed lines represent other stochastic bilevel optimization methods, while the solid lines depict the proposed methods.}
	\label{datacleaning1}
\end{figure}

\begin{figure}[H]
	\centering
	\subfigure{
		\centering
		\includegraphics[width=0.3\linewidth]{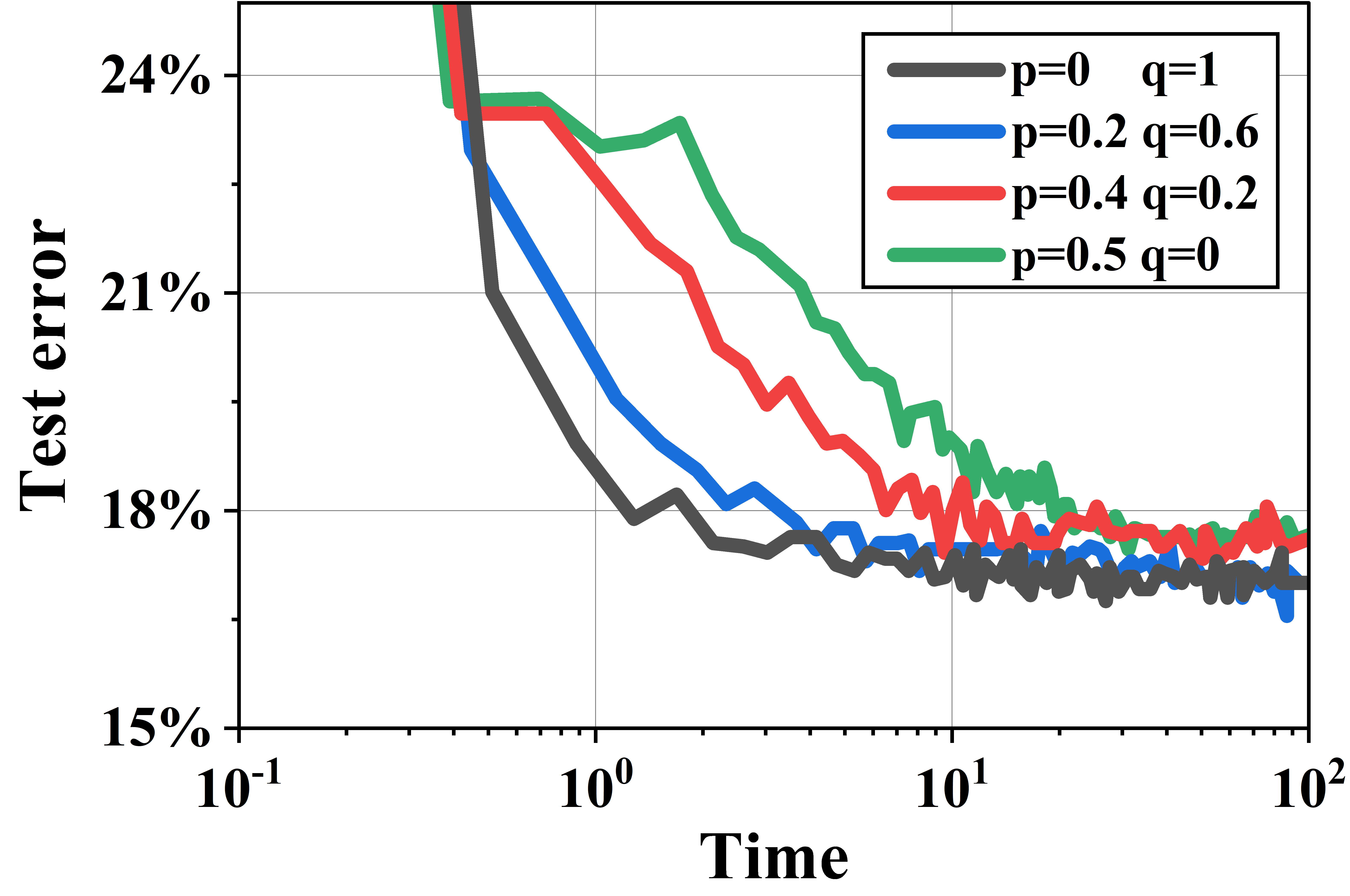}
		\label{qp}
	}
	\subfigure{
		\centering
		\includegraphics[width=0.3\linewidth]{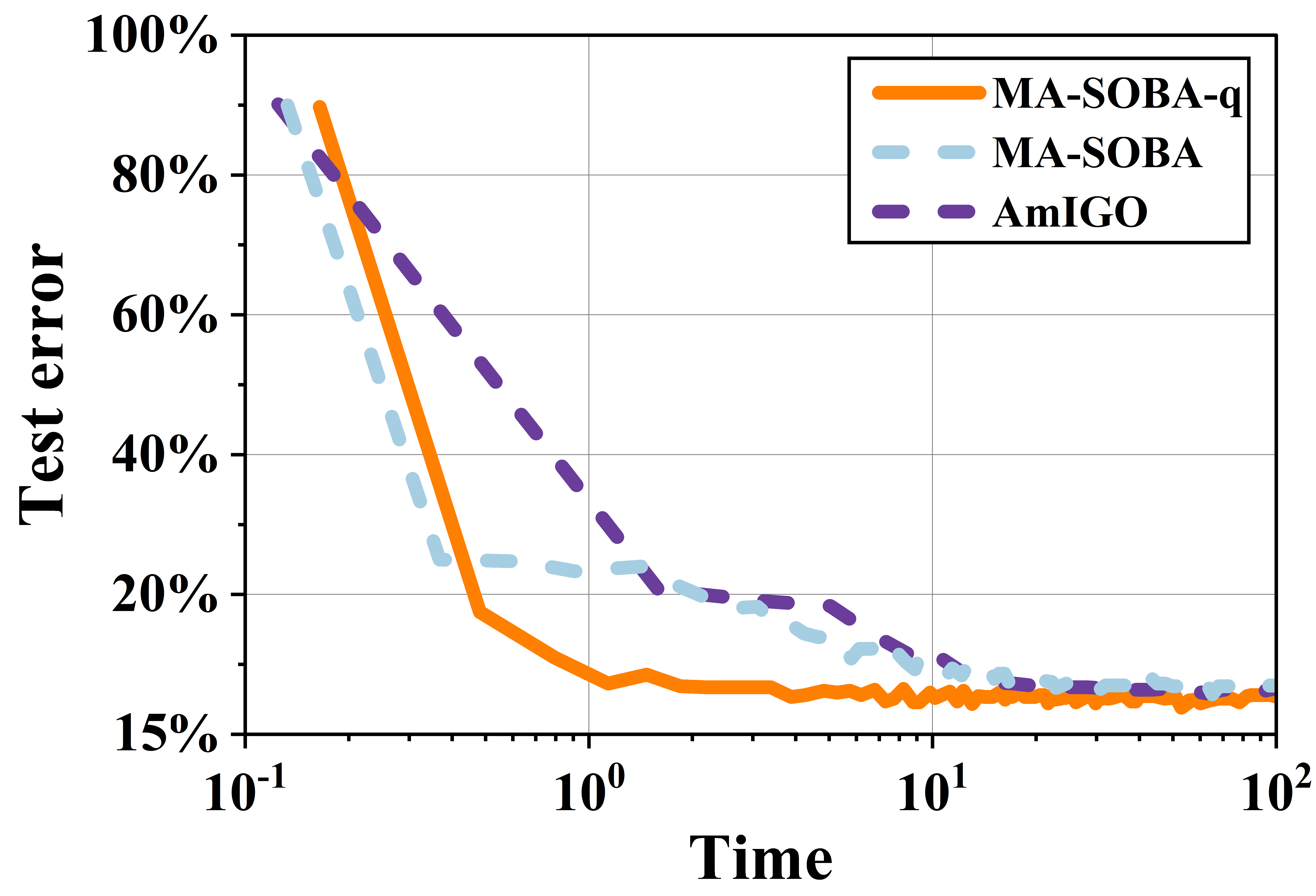}
		\label{SRMBA}
	}
	\caption{ \textbf{Left:} Compare by selecting different $q$ and $p$ in MA-SOBA-q with data hyper-cleaning on MINST. \textbf{Right:} Compare of MA-SOBA-q with other acceleration algorithms on hyper-cleaning on MINST. }
	\label{mMq}
\end{figure}

\subsection{Data hyper-cleaning}

Following the experimental setup in \cite{dagreou2022framework}, we identified the optimal value for the regularization parameter $C_r$ as 0.2 through a manual search, aiming to achieve the highest final test accuracy. It's worth noting that in this case, we were unable to utilize Just-In-Time (JIT) compilation from Numba due to the incompatibility of the softmax function from Scipy with Numba at the time of the experiment.

Figure \ref{datacleaning1} presents additional convergence curves with different methods. MA-SABA consistently emerges as the fastest algorithm to reach its final accuracy. Generally, the error decreases rapidly until it reaches a final value. Moreover, in Figure \ref{mMq}, we test the impact of the parameter $q$ on the algorithm MA-SOBA-q¡¯s performance.
We observe that as $q$ increases starting from $q = 0$, the convergence speed of the algorithm also accelerates, aligning with our theoretical expectations.

\textbf{Hyper-parameter setting for algorithm.} For MA-SABA, the step-sizes are chosen as $\alpha_k = 0.005/0.0002, \beta_k = 0.005, \gamma_k = 0.01$ and $ \rho_k = 0.2$. For SRMBA, the step-sizes are chosen as $\alpha_k = \frac{500}{k^{1/3}}, \beta_k  = \frac{0.2}{k^{1/3}}, \gamma_k = \frac{0.02}{k^{1/3}}$ and $ \rho_k^{x} = \rho_k^{y} = \rho_k^{z} = \frac{5}{k^{2/3}} $. Other algorithms choose their step sizes according to the optimal strategy in \cite{dagreou2022framework}.
In Figure \ref{mMq}, for MA-SOBA-q in Section \ref{detial:proofsoba}, the step-sizes and the batch-sizes are chosen as:
$ \alpha_k = 0.1/0.001, \beta_k = \gamma_k = 0.1, S = 1000 (p=0,q=1)$;
$ \alpha_k = 0.1/0.001, \beta_k = \gamma_k = 0.1, S = 1000 (p=0.2,q=0.6)$;
$ \alpha_k = 0.1/0.001, \beta_k = \gamma_k = 0.1, S = 1000 (p=0.4,q=0.2)$;
$ \alpha_k = 0.1/0.001, \beta_k = \gamma_k = 0.1, S = 1000 (p=0.5,q=0)$.

\section{Algorithms and General lemmas}\label{alg&lemma}
In this section, we present the specific forms of the algorithms MA-SABA, SPABA, and SRMBA and provide some general conclusions that are useful for the proof.

\subsection{Algorithms}\label{algos}

\begin{algorithm}[h]
  \caption{MA-SABA}
\begin{algorithmic}[1]
   \STATE {\bfseries Input:} Initializations $(x_{-1},y_{-1},z_{-1})$, $(x_{0},y_{0},z_{0})$, and $v_{-1}^x$, number of total iterations $K$, step size $\{\alpha_k, \beta_k,\gamma_k\}$, momentum parameter $\rho_k$;

   \FOR{$k=0$ {\bfseries to} $K-1$}
   \STATE Sample $i\in[n]$ for $f$ and $j\in[m]$ for $g$;
    \STATE
        $v^x_k=(1-\rho_{k-1})v_{k-1}^x+\rho_{k-1} D^x_{k-1}$;
   \STATE
       $x_{k+1} = x_{k}-\alpha_k v_k^x$;
    \STATE
        $D^x_{k}= \nabla_1 F_i(x_{k}, y_{k})-\nabla_1 F_i(w_{k,i}^x,w_{k,i}^y)+\frac{1}{n} \sum_{i=1}^n \nabla_1 F_i(w_{k,i}^x,w_{k,i}^y)-\nabla_{12}^2 G_j(x_{k}, y_{k}) z_{k}+\nabla_{12}^2 G_j(w_{k,j}^x,w_{k,j}^y) w_{k,j}^z +\frac{1}{m}\sum_{j=1}^{m}\nabla_{12}^2 G_j(w_{k,j}^x,w_{k,j}^y)w_{k,j}^z$
   \STATE
        $v^y_k=\nabla_2G_j(x_{k},y_{k})-\nabla_2G_j(w_{k,j}^x,w_{k,j}^y)+\frac{1}{m}\sum_{j=1}^{m}\nabla_2G_j(w_{k,j}^x,w_{k,j}^y);$
   \STATE
        $y_{k+1} =y_{k}-\beta_k v_k^y$;
   \STATE
      $v^z_k=
      \nabla_{22}^2 G_j(x_k, y_k) z_k-\nabla_{22}^2 G_j(w_{k,j}^x,w_{k,j}^y) w_{k,j}^z +\frac{1}{m}\sum_{j=1}^{m}\nabla_{22}^2 G_j(w_{k,j}^x,w_{k,j}^y)w_{k,j}^z -\nabla_2 F_i(x_k, y_k)+\nabla_2 F_i(w_{k,i}^x,w_{k,i}^y)-\frac{1}{n} \sum_{i=1}^n \nabla_2 F_i(w_{k,i}^x,w_{k,i}^y);$
   \STATE
       $z_{k+1}= z_{k}-\gamma_k v_k^z$.
   \ENDFOR
\end{algorithmic}
\end{algorithm}

\begin{algorithm}[h]
  \caption{SPABA}
\begin{algorithmic}[1]
   \STATE {\bfseries Input:} Initializations $(v^x_{-1},v^y_{-1},v^z_{-1})$, $(x_{-1},y_{-1},z_{-1})$ and $(x_{0},y_{0},z_{0})$, number of total iterations $K$, step size $\{\alpha_k, \beta_k,\gamma_k\}$, minibatch size $b$, constant $R$;
   \FOR{$k=0$ {\bfseries to} $K-1$}
   \STATE Sample $I$ for $f$ and $J$ for $g$, with minibatch size $|I|=|J|=b$;
    \STATE
        $
v_k^{x}= \begin{cases}\frac{1}{n} \sum_{i \in [n]} \nabla_{1} F_i\left(x_{k},y_{k}\right)-\frac{1}{m} \sum_{j \in [m]} \nabla^2_{12} G_j\left(x_{k},y_{k}\right)z_k, & \text { with probability } p,\\
v_{k-1}^{x}+\frac{1}{b} \sum_{i \in I}\left( \nabla_{1} F_i\left(x_{k},y_{k}\right)-\nabla_{1} F_i\left(x_{k-1},y_{k-1}\right)\right)\\ \quad- \frac{1}{b}  \sum_{j \in J}\left( \nabla^2_{12} G_j\left(x_{k},y_{k}\right)z_k- \nabla^2_{12} G_j\left(x_{k-1},y_{k-1}\right)z_{k-1}\right);
& \text { with probability } 1-p;
\end{cases}
$

   \STATE
       $x_{k+1} = x_{k}-\alpha_k v_k^x$;
    \STATE
        $
v_k^{y}= \begin{cases}\frac{1}{m} \sum_{j \in [m]} \nabla_2 G_j\left(x_{k},y_{k}\right), & \text { with probability } p,\\ v_{k-1}^{y}+\frac{1}{b} \sum_{j \in J}\left(\nabla_2 G_j\left(x_{k},y_{k}\right)-\nabla_2 G_j\left(x_{k-1},y_{k-1}\right)\right), & \text { with probability } 1-p;\end{cases}
$
   \STATE
        $y_{k+1} = y_{k}-\beta_k v_k^y$;
    \STATE
     $
v_k^{z}= \begin{cases}\frac{1}{m} \sum_{j \in [m]} \nabla^2_{22} G_j\left(x_{k},y_{k}\right)z_k- \frac{1}{n} \sum_{i \in [n]} \nabla_{2} F_i\left(x_{k},y_{k}\right), & \text { with probability } p,\\
v_{k-1}^{z}+\frac{1}{b} \sum_{j \in J}\left( \nabla^2_{22} G_j\left(x_{k},y_{k}\right)z_k- \nabla^2_{22} G_j\left(x_{k-1},y_{k-1}\right)z_{k-1}\right)\\\quad- \frac{1}{b} \sum_{i \in I}\left( \nabla_{2} F_i\left(x_{k},y_{k}\right)-\nabla_{2} F_i\left(x_{k-1},y_{k-1}\right)\right),
 & \text { with probability } 1-p;
\end{cases}
$
   \STATE
       $z_{k+1} = \text{Proj}_{\mathbb{B}(R)}(z_{k}-\gamma_k v_k^z)$.
   \ENDFOR
\end{algorithmic}
\end{algorithm}

\begin{algorithm}[h]
  \caption{SRMBA}
\begin{algorithmic}[1]
   \STATE {\bfseries Input:} Initializations $(x_{-1},y_{-1},z_{-1})$, $(x_{0},y_{0},z_{0})$, and $v_{-1}^x$, number of total iterations $K$, step size $\{\alpha_k, \beta_k,\gamma_k\}$, momentum parameter $\{\rho_k^x, \rho_k^y,\rho_k^z\}$,constant $R$ ;

   \FOR{$k=0$ {\bfseries to} $K-1$}
   \STATE Sample $\xi$ for $f$ and $\zeta$ for $g$;
    \STATE
        $D^x_{k}=\nabla_1 F(x_{k},y_{k};\xi )-\nabla_{12}^2G(x_{k},y_{k};\zeta )z_{k}$;
    \STATE
        $v^x_k=\rho^x_k D_{k}^x + (1-\rho^x_k)(v^x_{k-1}+D_k^x-D_{k-1}^x)$;
   \STATE
       $x_{k+1} = x_{k}-\alpha_k v_k^x$;
   \STATE
        $D^y_{k}=\nabla_2G(x_{k},y_{k};\zeta)$;
    \STATE
        $v^y_k=  \rho^y_k D_{k}^y + (1-\rho^y_k)(v^y_{k-1}+D_k^y-D_{k-1}^y)$;
   \STATE
        $y_{k+1} = y_{k}-\beta_k v_k^y$;
   \STATE
        $D^z_{k}=\nabla_{22}^2G(x_{k},y_{k};\zeta )z_{k}- \nabla_2 F(x_{k},y_{k};\xi)$;
    \STATE
      $ v^z_k=\rho^z_k D_{k}^z + (1-\rho^z_k)(v^z_{k-1}+D_k^z-D_{k-1}^z)$;
   \STATE
       $z_{k+1} =\text{Proj}_{\mathbb{B}(R)}( z_{k}-\gamma_k v_k^z)$.

   \ENDFOR
\end{algorithmic}
\end{algorithm}

\subsection{General lemmas}
In this section, we present general conclusions that will be used, including an important lemma on the descent of $H(x_k)$
\begin{lemma}\textbf{(Lipschitz continuity of $y^*(x)$)}\label{Ly*}

Under the Assumptions \ref{assump LL}, $y^*(x)$ is $L_{y^*}$-Lipschitz continuous,
where $L_{y^*}=\frac{L_1^g}{\mu}$.
\end{lemma}
\begin{proof}
See Lemma A.1 in \cite{dagreou2023lower}.
\end{proof}
\begin{lemma}\textbf{(Lipschitz continuity of $z^*(x)$)}

Under the Assumptions \ref{assump UL} and \ref{assump LL},
$z^*(x)$ is $L_{z^*}$ Lipschitz continuous, where $L_{z^*}=\left(\frac{L^f}{\mu}+\frac{C^f L_2^g}{\mu^2}\right)\left(1+\frac{L_1^g}{\mu}\right)$.
\end{lemma}
\begin{proof}
See Lemma A.1 in \cite{dagreou2023lower}.
\end{proof}
\begin{lemma}\label{R}\textbf{(boundness of $z^*(x)$)}

Under the Assumptions \ref{assump UL} and \ref{assump LL},
$z^*(x)$ is bounded by $R$, i.e., for each $x$, we have
$$ \|z^*(x)\|\leq \frac{C^f}{\mu}\triangleq R.$$
\end{lemma}
\begin{proof}
    See Lemma B.2 in \cite{chen2023optimal}.
\end{proof}

\begin{lemma}\label{Hsmooth}\textbf{(smoothness of function $H$)}

Suppose Assumptions \ref{assump UL} and \ref{assump LL} hold, the function $H(x)$ is $L^H$-smooth, where $$L^H=L^f+\frac{2L^f L^g_2+\left(C^f\right)^2L^g_2}{\mu}+\frac{L^f \left(L_1^g\right)^2+2C^fL_1^gL_2^g}{\mu^2}+\frac{C^f\left(L_1^g\right)^2L_2^g}{\mu^3}.$$
\end{lemma}
\begin{proof}
    See Lemma 2.2 in \cite{ghadimi2018approximation}.
\end{proof}

\begin{lemma}\label{dL}
Suppose Assumptions \ref{assump UL} and \ref{assump LL} hold.
Then the following inequalities hold:
\begin{align}
 \label{dl1}\mathbb{E} \left[\left\|D_y(x_k,y_k,z_k)\right\|^2\right] &\leq\left(L_1^g\right)^2 \mathbb{E}\left[\| y_k-y^*\left(x_k\right) \|^2\right],\\
 \mathbb{E}\left[\left\|D_z(x_k,y_k,z_k)\right\|^2\right]
&\leq L_z^2\mathbb{E}\left[\left\|z_k-z^*\left(x_k\right)\right\|^2\right]+ L_z^2\mathbb{E}\left[\left\|y_k-y^*\left(x_k\right)\right\|^2\right],\label{dl2}
\end{align}
where
$L_z^2=\max
\{3\left(L_1^g\right)^2,
3R^2\left(L_2^g\right)^2+3\left(L^f\right)^2\}$.
\end{lemma}

\begin{proof}
Proof of \eqref{dl1}:
Based on the fact that $D_y(x,y,z)=\nabla_2 g(x,y)$ and $D_y(x_k,y^*(x_k),z^*(x_k))=\nabla_2 g(x_k,y^*(x_k))=0$,
we have
\begin{align*}
\mathbb{E} \left[\left\|D_y(x_k,y_k,z_k)\right\|^2\right]
&= \mathbb{E} \left[\left\|D_y(x_k,y_k,z_k)-D_y(x_k,y^*(x_k),z^*(x_k))\right\|^2\right] \\
&= \mathbb{E} \left[\left\|\nabla_2 g(x_k,y_k)-\nabla_2 g(x_k,y^*(x_k))\right\|^2\right]\\
&\leq \left(L_1^g\right)^2\mathbb{E} \left[\left\|y_k-y^*(x_k)\right\|^2\right],
\end{align*}
where the last inequality utilizes the fact that $\nabla g$ is $L^g_1$-Lipschitz continuous, as stated in Assumption \ref{assump LL}.

Proof of \eqref{dl2}:
Based on the fact that $D_z(x,y,z)=\nabla_{22}^2 g(x,y)z-\nabla_2 f(x,y)$ and $D_z(x_k,y^*(x_k),z^*(x_k))=\nabla_{22}^2 g(x_k,y^*(x_k))z^*(x_k)-\nabla_2 f(x_k,y^*(x_k))=0,$
we have
\begin{align*}
&\mathbb{E}\left[\left\|D_z(x_k,y_k,z_k)\right\|^2\right]\\
&= \mathbb{E} \left[\left\|D_z(x_k,y_k,z_k)-D_z(x_k,y^*(x_k),z^*(x_k))\right\|^2\right]
\\
&= \mathbb{E} \left[\left\|\nabla_{22}^2 g(x_k,y_k)z_k-\nabla_2 f(x_k,y_k)-\nabla_{22}^2 g(x_k,y^*(x_k))z^*(x_k)+\nabla_2 f(x_k,y^*(x_k))\right\|^2\right]
\\
&\leq 3\mathbb{E} \left[\left\|\nabla_{22}^2 g(x_k,y_k)z_k-\nabla_{22}^2 g(x_k,y_k)z^*(x_k)\right\|^2\right]
\\
&\quad+3\mathbb{E} \left[\left\|\nabla_{22}^2 g(x_k,y_k)z^*(x_k)-\nabla_{22}^2 g(x_k,y^*(x_k))z^*(x_k)\right\|^2\right]
\\
&\quad+3\mathbb{E} \left[\left\|\nabla_2 f(x_k,y^*(x_k))-\nabla_2 f(x_k,y_k)\right\|^2\right]
\\
&\leq 3\left(L_1^g\right)^2\mathbb{E} \left[\left\|z_k-z^*(x_k)\right\|^2\right]
+3R^2\left(L_2^g\right)^2\mathbb{E} \left[\left\|y_k-y^*(x_k)\right\|^2\right]
+3\left(L^f\right)^2\mathbb{E} \left[\left\|y_k-y^*(x_k)\right\|^2\right]
\\
&=3\left(L_1^g\right)^2\mathbb{E} \left[\left\|z_k-z^*(x_k)\right\|^2\right]
+3\left(R^2\left(L_2^g\right)^2+\left(L^f\right)^2\right) \mathbb{E}\left[\left\|y_k-y^*(x_k)\right\|^2\right]\\
&\leq
L_z^2\mathbb{E}\left[\left\|z_k-z^*\left(x_k\right)\right\|^2\right]
+ L_z^2\mathbb{E}\left[\left\|y_k-y^*\left(x_k\right)\right\|^2\right],
\end{align*}
where the validity of the second inequality is based on the application of Lemma \ref{R},
along with the assumptions that $\nabla f$ is $L^f$-Lipschitz continuous as stated in Assumption \ref{assump UL},
and $\nabla^2 g$ is $L^g_2$-Lipschitz continuous as mentioned in Assumption \ref{assump LL}.
The last inequality is due to the fact that
$L_z^2=\max
\{3\left(L_1^g\right)^2,
3R^2\left(L_2^g\right)^2+3\left(L^f\right)^2\}.$
\end{proof}

\begin{lemma}\label{DH}
Suppose Assumptions \ref{assump UL} and \ref{assump LL} hold. Then we have
\begin{eqnarray*}
\mathbb{E}\left[\left\|D_x(x_k,y_k,z_k)-\nabla H\left(x_k\right)\right\|^2\right]
 \leq 3\left(\left(L^f\right)^2+\left(L^g_2 R\right)^2\right) \mathbb{E}\left[\left\|y_k-y^*\left(x_k\right)\right\|^2\right] +3\left(L^g_1\right)^2 \mathbb{E}\left[\left\|z_k-z^*\left(x_k\right)\right\|^2\right].
\end{eqnarray*}
\end{lemma}
\begin{proof}
Using the  unbiasedness of $D_k^x$ and the Cauchy-Schwarz inequality, we have
\begin{eqnarray*}
&&\mathbb{E}\left[\left\|D_x(x_k,y_k,z_k)-\nabla H\left(x_k\right)\right\|^2\right]\\
&=&\mathbb{E}\left[\left\| \nabla_1 f\left(x_k, y_k\right)
+
\nabla_{12}^2 g\left(x_k, y_k\right) z_k
-
\nabla_1 f\left(x_k, y^*\left(x_k\right)\right)\right.\right.
\\
&&\quad\left.\left.
-
\nabla_{12}^2 g\left(x_k, y_k\right) z^*\left(x_k\right)
+
\nabla_{12}^2 g\left(x_k, y_k\right) z^*\left(x_k\right)
-
\nabla_{12}^2 g\left(x_k, y^*\left(x_k\right)\right) z^*\left(x_k\right) \right\|^2\right]\\
& \leq&
3\left(\mathbb{E} \left[\| \nabla_1 f\left(x_k, y_k\right)-\nabla_1 f\left(x_k, y^*\left(x_k\right) )\|^2\right]
+
\mathbb{E}\left[\left\|\nabla_{12}^2 g\left(x_k,y_k\right)\left(z_k-z^*\left(x_k\right)\right)\right\|^2\right]\right.\right. \\
&&\quad\left.+\mathbb{E}\left[\left\|\nabla_{12}^2 g\left(x_k, y_k\right)-\nabla_{12}^2 g\left(x_k, y^*\left(x_k\right)\right) \right\|^2\left\|z^*\left(x_k\right) \right\|^2\right]\right).
\end{eqnarray*}
The three terms on the right-hand side of the above inequality can be bounded by utilizing Assumption \ref{assumporacle} and Lemma \ref{R}.
Thus, the lemma is proven.
\end{proof}
\begin{lemma}\label{Hlemma}
Suppose Assumptions \ref{assump UL} and \ref{assump LL} hold. Then we have
\begin{eqnarray*}
\mathbb{E}\left[H\left(x_{k+1}\right)\right]\leq \mathbb{E}\left[H\left(x_k\right)\right]
-\frac{\alpha_k}{2} \mathbb{E}\left[\left\|\nabla H\left(x_k\right)\right\|^2\right]
+
 \left(\frac{L^H\alpha_k^2}{2}-\frac{\alpha_k}{2}\right) \mathbb{E}\left[\left\|v_k^x\right\|^2\right]
 +\frac{\alpha_k}{2} \mathbb{E}\left[\left\|\nabla H\left(x_k\right)-v_k^x\right\|^2\right].
\end{eqnarray*}
\end{lemma}
\begin{proof}
    The $L^H$-smoothness of $H(x)$ in Lemma \ref{Hsmooth} implies
\begin{eqnarray*}
H\left(x_{k+1}\right)-H\left(x_k\right)
&\leq &
\left\langle\nabla H\left(x_k\right), x_{k+1}-x_k\right\rangle
+
\frac{L^H}{2}\left\|x_{k+1}-x_k\right\|^2 \\
&= & -\alpha_k\left\langle\nabla H\left(x_k\right), v_k\right\rangle
+\frac{L^H}{2} \alpha_k^2\left\|v^x_k\right\|^2 \\
& =&
-\frac{\alpha_k}{2}\left[\left\|\nabla H\left(x_k\right)\right\|^2
+\left\|v_k\right\|^2-\left\|\nabla H\left(x_k\right)-v^x_k\right\|^2\right]+\frac{L^H}{2} \alpha_k^2\left\|v^x_k\right\|^2 \\
& =&-\frac{\alpha_k}{2}\left\|\nabla H\left(x_k\right)\right\|^2+\left(\frac{L^H}{2} \alpha_k^2-\frac{\alpha_k}{2}\right)\left\|v^x_k\right\|^2 +\frac{\alpha_k}{2}\left\|\nabla H\left(x_k\right)-v^x_k\right\|^2,
\end{eqnarray*}
where the second equation uses the fact that $\langle a, b\rangle=\frac{1}{2}\left(\|a\|^2+\|b\|^2-\|a-b\|^2\right)$. Taking expectation on both sides, we know (\ref{Hlemma}) holds.
\end{proof}

\section{Proof of Theorem \ref{thsaba}}\label{detial:proofmasaba}

\begin{lemma}\label{vh}
Suppose Assumptions \ref{assump UL} and \ref{assump LL} hold. Then we have
\begin{align*}
\mathbb{E}\left[\left\|v_{k+1}^x-\nabla H\left(x_{k+1}\right)\right\|^2\right] \leq& \left(1-\rho_k\right) \mathbb{E}\left[\left\|v_k^x-\nabla H\left(x_k\right)\right\|^2\right]+\frac{2\left(L^H\right)^2 \alpha_k^2}{\rho_k} \mathbb{E}\left[\left\|v_k^x\right\|^2\right]
\\&+2\rho_k\mathbb{E}\left[\left\|E[D_k^x]-\nabla H(x_k)\right\|^2\right]
+\rho_k^2 \mathbb{E}\left[\left\|D_k^x-\mathbb{E}\left[D_k^x\right]\right\|^2\right],
\end{align*}
where we require that $0\leq\rho_k\leq 1$.
\end{lemma}
\begin{proof}
Due to iteratively updating $v_k^x$, we have
\begin{eqnarray*}
&&\mathbb{E}\left[\left\|v^x_{k+1}-\nabla H\left(x_{k+1}\right)\right\|^2\right] \\
&= & \mathbb{E}\left[\left\|\left(1-\rho_k\right) v^x_k+\rho_k D_k^x-\nabla H\left(x_{k+1}\right)\right\|^2\right] \\
&=&\mathbb{E}  \left[\|\left(1-\rho_k\right)\left(v^x_k-\nabla H\left(x_k\right)\right)-\rho_k \nabla H\left(x_k\right)+\rho_k \mathbb{E}\left[D_k^x\right]-\rho_k \mathbb{E}\left[D_k^x\right]+\rho_k D_k^x+\nabla H\left(x_k\right)-\nabla H\left(x_{k+1}\right) \|^2\right] \\
&=& \mathbb{E}\left[\left\|\left(1-\rho_k\right)\left(v^x_k-\nabla H\left(x_k\right)\right)-\rho_k \nabla H\left(x_k\right)+\rho_k \mathbb{E}\left[D_k^x\right]+\nabla H\left(x_k\right)-\nabla H\left(x_{k+1}\right)\right\|^2\right]  \\
&&+  \rho_k^2 \mathbb{E}\left[\left\|D_k^x-\mathbb{E}\left[D_k^x\right]\right\|^2\right] \\
&\leq & \left(1-\rho_k\right) \mathbb{E}\left[\left\|v^x_k-\nabla H\left(x_k\right)\right\|^2\right]+\rho_k \mathbb{E}\left[\left\|\mathbb{E}\left[D_k^x\right]-\nabla H\left(x_k\right)+\frac{\nabla H\left(x_k\right)-\nabla H\left(x_{k+1}\right)}{\rho_k}\right\|^2\right] \\
& &+\rho_k^2 \mathbb{E}\left[\left\|D_k^x-\mathbb{E}\left[D_k^x\right]\right\|^2\right] \\
&\leq & \left(1-\rho_k\right) \mathbb{E}\left[\left\|v^x_k-\nabla H\left(x_k\right)\right\|^2\right]+2 \rho_k \mathbb{E}\left[\left\|\mathbb{E}\left[D_k^x\right]-\nabla H\left(x_k\right)\right\|^2\right] \\
& &+\frac{2 \alpha_k^2\left(L^H\right)^2}{\rho_k} \mathbb{E}\left[\left\|v^x_k\right\|^2\right]+\rho_k^2 \mathbb{E}\left[\left\|D_k^x-\mathbb{E}\left[D_k^x\right]\right\|^2\right],
\end{eqnarray*}
where the third equation uses the unbiasedness of $D_k^x$, the first inequality is due to the convexity of $\|\cdot\|^2$, and the second inequality uses the $L^H-$smoothness of $H$.
\end{proof}

\begin{lemma}\label{y,z}
Suppose Assumption \ref{assump UL} and \ref{assump LL} hold and  the step size satisfy $$\beta_k\leq1/(\mu+L^g_1),\quad  \gamma_k\leq1/(10\mu).$$
Then we have
\begin{align*}
\mathbb{E}\left[\left\| y_{k+1}-y^*\left(x_{k+1}\right)\right\|^2\right]-\mathbb{E}\left[\left\|y_k-y^*\left(x_k\right)\right\|^2\right]
 \leq&-\beta_k \mu \mathbb{E}\left[\left\|y_k-y^*\left(x_k\right)\right\|^2\right] +\frac{2 L_{y^*}^2 \alpha_k^2}{\beta_k \mu} \mathbb{E}\left[\left\|v_k^x\right\|^2\right]
 \\&+2 \beta_k^2\mathbb{E}\left[\left\|D_y(x_k,y_k,z_k)-D_k^y\right\|^2\right].
\end{align*}
\begin{align*}
\mathbb{E}\left[\left\|z_{k+1}-z^*\left(x_{k+1}\right)\right\|^2\right]-\mathbb{E}\left[\left\|z_k-z^*\left(x_k\right)\right\|^2\right]
\leq&-\gamma_k \mu \mathbb{E}\left[\left\|z_k-z^*\left(x_k\right)\right\|^2\right]
+8\Delta \gamma_k \mathbb{E}\left[\left\| y_k-y^* (x_k)\right\|^2\right]
\\&+2 \gamma_k^2\mathbb{E}\left[\left\|D_z(x_k,y_k,z_k)-D_k^z\right\|^2\right]
 +\frac{3 L_{z^*}^2 \alpha_k^2}{\gamma_k \mu} \mathbb{E}\left[\left\|v_k^x\right\|^2\right],
\end{align*}
where $\Delta=\left(\left(L_2^g R\right)^2+\left(L^f\right)^2\right)/\mu$.
\end{lemma}
\begin{proof}
    \textbf{Inequality for $y$}

We use Young inequality to start by expanding the square
\begin{eqnarray*}
\left\|y_{k+1}-y^*\left(x_{k+1}\right)\right\|^2
&= & \left\|y_{k+1}-y^*\left(x_k\right)+y^*\left(x_k\right)-y^*\left(x_{k+1}\right)\right\|^2 \\
&= & \left\|y_{k+1}-y^*\left(x_k\right)\right\|^2+\left\|y^*\left(x_k\right)-y^*\left(x_{k+1}\right)\right\|^2+2\left\langle y_{k+1}-y^*\left(x_k\right), y^*\left(x_k\right)-y^*\left(x_{k+1}\right)\right\rangle \\
&\leq & \left(1+\beta_k \mu\right)\left\|y_{k+1}-y^*\left(x_k\right)\right\|^2+\left(1+\frac{1}{\beta_k \mu}\right)\left\|y^*\left(x_k\right)-y^*\left(x_{k+1}\right)\right\|^2
\\ & \leq & \left(1+\beta_k \mu\right)\left\|y_{k+1}-y^*\left(x_k\right)\right\|^2
+ \left(1+\frac{1}{\beta_k \mu}\right)L_{y^*}^2\alpha_k^2\left\|v_k^x\right\|^2,
\end{eqnarray*}
where the last inequality is due to Lemma \ref{Ly*}.

Taking the expectation conditionally on $x_k$, $y_k$, $z_k$ yields
\begin{eqnarray}\label{y1}
E_k\left[\left\|y_{k+1}-y^*\left(x_{k+1}\right)\right\|^2\right]
\leq
\left(1+\beta_k\mu\right)E_k\left[\left\|y_{k+1}-y^*\left(x_k\right)\right\|^2\right]
+\left(1+\frac{1}{\beta_k \mu}\right)L_{y^*}^2\alpha_k^2E_k[\left\|v_k^x\right\|^2].
\end{eqnarray}
For the first member, we have
\begin{eqnarray*}
&&E_k\left[\left\|y_{k+1}-y^*\left(x_k\right)\right\|^2\right]\\
& =&E_k\left[\left\|y_k-y^*\left(x_k\right)-\beta_k v_k^y\right\|^2\right]
 \\
& =&E_k\left[\left\|y_k-\beta_k D_y(x_k,y_k,z_k)-y^*\left(x_k\right)-\beta_k\left(v_k^y-D_y(x_k,y_k,z_k)\right)\right\|^2\right] \\
& =&E_k\left[\left\|y_k-\beta_k D_y(x_k,y_k,z_k)-y^*\left(x_k\right)\right\|^2\right]+E_k\left[\| \beta_k\left(v_k^y-D_y(x_k,y_k,z_k) \right)\|^2\right] \\
& &+2 \left.E_k\left[\langle y_k-\beta_k D_y(x_k,y_k,z_k)-y^*\left(x_k\right) , \beta_k\left(v_k^y-D_y(x_k,y_k,z_k)_k^*\right)\right\rangle\right] \\
& =&E_k\left[\left\|y_k-\beta_k D_y(x_k,y_k,z_k)-y^*\left(x_k\right)\right\|^2\right]+E_k\left[\| \beta_k\left(v_k^y-D_y(x_k,y_k,z_k)\right) \|^2\right] \\
& \leq&\left(1-\beta_k \mu\right)^2\left\|y_k-y^*\left(x_k\right)\right\|^2+\beta_k^2 E_k\left[\|v_k^y-D_y(x_k,y_k,z_k) \|^2\right],
\end{eqnarray*}
where the first inequality holds because $D_k^y$ is an unbiased estimate of $D_y(x_k,y_k,z_k)$.
The first inequality utilizes Lemma 10 in \cite{qu2017harnessing} which requires that
$g$ is strongly convex and Lipschitz smooth.
Plugging it into (\ref{y1}) and taking the total expectation, we have
\begin{eqnarray*}
&&\mathbb{E}\left[\left\|y_{k+1}-y^*\left(x_{k+1}\right)\right\|^2\right]\\
&\leq&
\left(1+\beta_k\mu\right)\left(1-\beta_k \mu\right)^2\mathbb{E}[\left\|y_k-y^*\left(x_k\right)\right\|^2]\\&&+\left(1+\beta_k\mu\right)\beta_k^2 \mathbb{E}\left[\|v_k^y-D_y(x_k,y_k,z_k) \|^2\right]
+\left(1+\frac{1}{\beta_k \mu}\right)L_{y^*}^2\alpha_k^2\mathbb{E}[\left\|v_k^x\right\|^2]\\
&\leq&
\left(1-\beta_k \mu\right)\mathbb{E}[\left\|y_k-y^*\left(x_k\right)\right\|^2+2\beta_k^2 \mathbb{E}\left[\|v_k^y-D_y(x_k,y_k,z_k) \|^2\right]
+\frac{2L_{y^*}^2\alpha_k^2}{\beta_k \mu}\mathbb{E}[\left\|v_k^x\right\|^2],
\end{eqnarray*}
where the last inequality is due to $\beta_k\leq1/(\mu+L^g_1)$.

\textbf{Inequality for $z$}

Similar to the analysis of $\mathbb{E}\left[\|y(x_k)-y^*(x_k)\|^2\right]$, we analyze the auxiliary variable $z$
\begin{eqnarray} \left\|z_{k+1}-z^*\left(x_{k+1}\right)\right\|^2&=&\left\|z_{k+1}-z^*\left(x_k\right)+z^*\left(x_k\right)-z^*\left(x_{k+1}\right)\right\|^2 \notag\\
& \leq&\left(1+\frac{\gamma_k \mu}{2}\right)\left\|z_{k+1}-z^*\left(x_k\right)\right\|^2 +\left(1+\frac{2}{\gamma_k \mu}\right)\left\|z^*\left(x_k\right)-z^*\left(x_{k+1}\right)\right\|^2.\label{z0}
\end{eqnarray}
For the second term, taking total expectation and utilizing the Lipschitz continuity of $z*(x)$, we have
\begin{eqnarray}
\mathbb{E}\left[\left\|z^*\left(x_k\right)-z^*\left(x_{k+1}\right)\right\|^2\right] \leq L_{z^*}^2 \alpha_k^2 \mathbb{E}\left[\left\|v^x_k\right\|^2\right] . \label{z1}\end{eqnarray}
The analysis of the first term is more complex.
Based on the definition of
$z_{k+1}$ and the fact that
$D_k^z$ is an unbiased estimate of
$D_z(x_k,y_k,z_k)$, we have
\begin{eqnarray}
&&E_k\left[\left\|z_{k+1}-z^*\left(x_k\right)\right\|^2\right]\\
&=&E_k\left[\left\|z_k-\gamma_k v_z^k-z^*\left(x_k\right)\right\|^2\right] \notag \\
& =&E_k\left[\left\|z_k-\gamma_k D_z(x_k,y_k,z_k)-z^*\left(x_k\right)-\gamma_k\left(v_k^z-D_z(x_k,y_k,z_k)\right)\right\|^2\right]\notag  \\
& =&E_k\left[\left\|z_k-\gamma_k D_z(x_k,y_k,z_k)-z^*\left(x_k\right)\right\|^2\right]+\gamma_k^2E_k\left[\| v_k^z-D_z(x_k,y_k,z_k) \|^2\right].\label{z2}
\end{eqnarray}
According to the definition of $D_z$ and $D_z(x_k,y^*(x_k),z^*(x_k))=0$, we have
\begin{eqnarray}
&&E_k\left[\left\|z_k-\gamma_k D_z(x_k,y_k,z_k)-z^*\left(x_k\right)\right\|^2\right]\notag\\
&& = E_k\left[\left\|z_k-\gamma_k \left[\nabla_{22}^2g(x_k,y_k)z_k-\nabla_2f(x_k,y_k)\right]-z^*\left(x_k\right)\right\|^2\right]
\notag\\
&& =E_k\left[\left\|z_k-z^*(x_k)-\gamma_k \nabla_{22}^2g(x_k,y_k)z_k
+\gamma_k \nabla_{22}^2g(x_k,y_k)z^*_k
-\gamma_k \nabla_{22}^2g(x_k,y_k)z^*_k
+\gamma_k \nabla_{22}^2g(x_k,y^*(x_k))z^*_k\right.
\right.
\notag\\&&\left.
\left.\quad+\gamma_k\nabla_2f(x_k,y_k)
-\gamma_k\nabla_2f(x_k,y^*(x_k))
\right\|^2\right]\notag\\
&& = E_k\left[\left\|\right.
\right.(I-\gamma_k\nabla_{22}^2g(x_k,y_{k}))(z_k-z^*(x_k))
\notag\\
&&\left.\left.\quad+\gamma_k\left[(\nabla_{22}^2g(x_k,y^*(x_k))-\nabla_{22}^2g(x_k,y_k))z^*(x_k)+\nabla_2f(x_k,y_k)-\nabla_2f(x_k,y^*(x_k))\right] \right\|^2\right]\notag\\
&& \leq \left(1+\frac{\gamma_k\mu}{3}\right)E_k\left[\left\|(I-\gamma_k\nabla_{22}^2g(x_k,y_{k}))(z_k-z^*(x_k))\right\|^2\right]\notag\\
&&\quad+\left(2+\frac{6}{\gamma_k\mu}\right)\gamma_k^2\left[E_k\|\nabla_{22}^2g(x_k,y^*(x_k))-\nabla_{22}^2g(x_k,y_k)\|^2\|z^*(x_k)\|^2+E_k\|\nabla_2f(x_k,y_k)-\nabla_2f(x_k,y^*(x_k))\|^2\right]\notag\\
&&\leq \left(1+\frac{\gamma_k\mu}{3}\right)(1-\gamma_k\mu)^2E_k\left[\left\|z_k-z^*(x_k))\right\|^2\right]+\left(2+\frac{6}{\gamma_k\mu}\right)\gamma_k^2\left(\left(L_2^gR\right)^2+(L^f)^2\right)E_k[\|y_k-y^*(x_k)\|^2].\label{z3}
\end{eqnarray}
Combining (\ref{z0}), (\ref{z1}), (\ref{z2}), and (\ref{z3}) and taking the total expectation, we have
\begin{eqnarray*}
\mathbb{E}\left[\left\|z_{k+1}-z^*(x_{k+1})\right\|^2\right]
&\leq&
\left(1+\frac{\gamma_k \mu}{2}\right)\left(1+\frac{\gamma_k\mu}{3}\right)(1-\gamma_k\mu)^2\mathbb{E}\left[\left\|z_k-z^*(x_k))\right\|^2\right]\\
&&+\left(1+\frac{\gamma_k \mu}{2}\right)\left(2+\frac{6}{\gamma_k\mu}\right)\gamma_k^2\left(\left(L_2^gR\right)^2+(L^f)^2\right)\mathbb{E}[\|y_k-y^*(x_k)\|^2]
\\&&+\left(1+\frac{\gamma_k \mu}{2}\right)\gamma_k^2\mathbb{E}\left[\| v_k^z-D_z(x_k,y_k,z_k) \|^2\right]
\\&&+\left(1+\frac{2}{\gamma_k \mu}\right)L_{z^*}^2 \alpha_k^2 \mathbb{E}\left[\left\|v^x_k\right\|^2\right]
\\
&\leq&
(1-\gamma_k\mu)\mathbb{E}\left[\left\|z_k-z^*(x_k))\right\|^2\right]+ 8\gamma_k \Delta\mathbb{E}[\|y_k-y^*(x_k)\|^2]
\\&&+2\gamma_k^2\mathbb{E}\left[\| v_k^z-D_z(x_k,y_k,z_k) \|^2\right]
+\frac{3L_{z^*}^2\alpha_k^2}{\gamma_k \mu}  \mathbb{E}\left[\left\|v^x_k\right\|^2\right].
\end{eqnarray*}
For convenience of expression, let's denote $\Delta=\left(\left(L_2^g R\right)^2+\left(L^f\right)^2\right)/\mu$,
and the last inequality is based on the choice of $\gamma_k\leq 1/(10\mu)$.
\end{proof}

To facilitate the discussion, we define $S_k=E_{k,f}^z+E_{k,f}^y+E_{k,f}^x+E_{k,g}^z+E_{k,g}^y+E_{k,g}^x$,
where
\begin{align*}
    E_{k,f}^z=\frac{1}{n}\sum_{i=1}^n\mathbb{E}[\|z_k-w_{k,i}^z\|^2],
\end{align*}
and similarly $E_{k,f}^z$, $E_{k,f}^y$,$E_{k,f}^x$, $E_{k,g}^z$, $E_{k,g}^y$ and $E_{k,g}^x$.
Additionally, let
$$\tau=\min\left\{\frac{1}{2n},\frac{1}{2m}\right\}.$$
\begin{lemma}\label{var}
Suppose Assumption \ref{assump UL}, \ref{assump LL} and \ref{assumptionijsmooth} hold,
there exist positive constants  $L_x'$, $L_x''$, $L_z'$ and $L_z''$ such that
\begin{align*}
    \mathbb{E}\left[\left\|D_y(x_k,y_k,z_k)-D_k^y\right\|^2\right]
    &\leq \left(L_{1}^g\right)^2S_k, \\
    \mathbb{E}\left[\left\|D_z(x_k,y_k,z_k)-D_k^z\right\|^2\right]
    &\leq L_z'S_k+L_z''\left(\mathbb{E}\left[\left\|y_k-y^*(x_k)\right\|^ 2\right]+\mathbb{E}\left[\left\|z_k-z^*(x_k)\right\|^ 2\right]\right),\\
    \mathbb{E}\left[\left\|D_x(x_k,y_k,z_k)-D_k^x\right\|^2\right]
    &\leq L_x'S_k+L_x''\left(\mathbb{E}\left[\left\|y_k-y^*(x_k)\right\|^ 2\right]+\mathbb{E}\left[\left\|z_k-z^*(x_k)\right\|^ 2\right]\right).
\end{align*}
where
\begin{eqnarray*}
L_x^{\prime}=L_z^{\prime}&=&\max\{16\left(L_2^g\right)^2 R^2,16\left(L_1^g\right)^2,2(L^f)^2\}, \\
L_x^{\prime \prime}=L_z^{\prime \prime}&=&\max\{24\left(L_1^g\right)^2,24 R^2\left(L_2^g\right)^2 \}.
\end{eqnarray*}
\end{lemma}

\begin{proof}
Assuming we sample $i$ and $j$ from $[n]$ and $[m]$ at iteration $k$ respectively, then we have
\begin{align*}
& {E}_k\left[\left\|D_y\left(x_k, y_k, z_k\right)-D_k^y\right\|^2\right] \\
& =E_k[\| \nabla_2 g\left(x_k, y_k\right)-\nabla_2 G_j\left(x_k, y_k\right)+\nabla_2 G_j\left(w_{k, j}^x-w_{k, j}^y\right)-\frac{1}{m} \sum_{j=1}^m \nabla_2 G_j\left(w_{k, j}^x, w_{k, j}^y\right) \|^2] \\
& \leq E_k\left[\left\|\nabla_2 G_j\left(x_k, y_k\right)-\nabla_2 G_j\left(w_{k, j}^x, w_{k, j}^y\right)\right\|^2\right] \\
& =\frac{1}{m} \sum_{j=1}^m E_k\left[\left\|\nabla_2 G_j\left(x_k, y_k\right)-\nabla_2 G_j\left(w_{k, j}^x, w_{k, j}^y\right)\right\|^2\right] \\
& \leq \left(L_1^g\right)^2\left(\frac{1}{m} \sum_{j=1}^m \mathbb{E}\left[\left\|x_k-w_{k, j}^x\right\|^2\right]+\frac{1}{m} \mathbb{E}\left[\left\|y_k- w_{k, j}^y\right\|^2\right]\right) ,
\end{align*}
where the first inequality uses the fact that $E[(X-E[X])^2]\leq E[X^2]$, the second inequality uses the Lipschitz continuity of $\nabla G_j$.
Taking the total expectation and by the definition of $S_k$, we can obtain
\begin{align*}
    \mathbb{E}\left[\left\|D_y(x_k,y_k,z_k)-D_k^y\right\|^2\right]
& =\left(L_1^g\right)^2\left(E_{k,g}^x+E_{k,g}^y\right)  \leq  \left(L_1^g\right)^2 S_k.
\end{align*}
For $x$, we have
\begin{align*}
& E_k\left[\left\|\nabla_x\left(x_k, y_k, z_k\right)-\nabla_k^x\right\|^2\right] \\
=&E_k\left[\| \nabla_{22}^2 g\left(x_k, y_k\right) z_k-\nabla_2 f\left(x_k, y_k\right)-\nabla_{22}^2 G_j\left(x_k, y_k\right) z_k+\nabla_{22}^2 G_j\left(w_{k, j}^x, w_{k, j}^y\right) w_{k,j}^z\right.
\\&\left.\quad -\frac{1}{m} \sum_{j=1}^m \nabla_{22}^2 G_j\left(w_{k, j}^x, w_{k, j}^y\right) w_{k, j}^z
+\nabla_2 F_i\left(x_k, y_k\right)-\nabla_2 F_i\left(w_{k, i}^x, w_{k, i}^y\right) +\frac{1}{n} \sum_{i=1}^n \nabla_2 F_i\left(w_{k, i}^x, w_{k, i}^y\right) \|^2\right] \\
& \leq E_k\left[\| \nabla_{22}^2 G_j\left(w_{k, j}^x, w_{k, j}^y\right) w_{k, j}^z-\nabla_{22}^2 G_j\left(x_k, y_k\right) z_k+\nabla_2 F_i\left(x_k, y_k\right)-\nabla_2 F_i\left(w_{k, i}^x, w_{k, i}^y\right) w_{k, i}^z\|^2\right] \\
& \leq 2 E_k\left[\left\|\nabla_{22}^2 G_j\left(w_{k, j}^x, w_{k, j}^y\right) w_{k, j}^z-\nabla_{22}^2 G_j\left(x_k, y_k\right) z_k\right\|^2\right] +2 E_k\left[\left\|\nabla F_i\left(x_k, y_k\right)-\nabla F_i\left(w_{k, i}^x, w_{k, i}^y\right)\right\|^2\right] .
\end{align*}
Taking the total expectation of the first term and we have
\begin{align*} &\mathbb{E}\left[E_k\left[\left\|\nabla_{22}^2 G_j\left(w_{k, j}^x, w_{k, j}^y\right) w_{k, j}^z-\nabla_{22}^2 G_j\left(x_k, y_k\right) z_k\right\|^2\right]\right]\\
=&
\mathbb{E}[E_k[\| \nabla_{22}^2 G_j\left(x_k, y_k\right) z_k-\nabla_{22}^2 G_j\left(x_k, y_k\right) z^*\left(x_k\right)+\nabla_{22}^2 G_j\left(x_k, y_k\right) z^*\left(x_k\right)  -\nabla_{22}^2 G_j\left(x_k, y^*\left(x_k\right)\right) z^*\left(x_k\right)
\\&+\nabla_{22}^2 G_j\left(x_k, y^*\left(x_k\right)\right) z^*\left(x_k\right) -\nabla_{22}^2 G_j\left(w_{k, j}^x, w_{k, j}^y\right) z^*\left(x_k\right)
\\
& \left.\left.+\nabla_{22}^2 G_j\left(w_{k, j}^x, w_{k, j}^y\right) z^*\left(x_k\right) -\nabla_{22}^2 G_j\left(w_{k, j}^x, w_{k, j}^y\right) w_{k, j}^z \|^2\right]\right] \\
\leq & 4\left(L_1^g\right)^2 \mathbb{E}\left[\left\|z_k-z^*\left(x_k\right)\right\|^2\right]+4 R^2\left(L_2^g\right)^2 \mathbb{E}\left[\left\|y_k-y^*\left(x_k\right)\right\|^2\right] \\
&+  4\left(L_2^g\right)^2 R^2\left(\frac{1}{m} \sum_{j=1}^m\mathbb{E}\left[\left\|x_k-w_{k, j}^x\right\|^2\right]+2 E\left[\left\|y^*\left(x_k\right)-y_k\right\|^2\right]+ \frac{2}{m} \sum_{j=1}^mE\left[\left\|y_k-w_{k, j}^k\right\|^2\right]\right) \\
&+  4\left(L_1^g\right)^2\left(2 \mathbb{E}\left[\left\|z^*\left(x_k\right)-z_k\right\|^2\right]+ \frac{2}{m} \sum_{j=1}^m\mathbb{E}\left[\left\|z_k-w_{k, j}^z\right\|^2\right]\right)\\
=&12\left(L_1^g\right)^2
\mathbb{E}\left[\left\|z_k-z^*\left(x_k\right)\right\|^2\right]
+12 R^2\left(L_2^g\right)^2 \mathbb{E}\left[\left\|y_k-y^*\left(x_k\right)\right\|^2\right]
\\&+4\left(L_2^g\right)^2 R^2E_{k,g}^x
+8\left(L_2^g\right)^2 R^2E_{k,g}^y
+8\left(L_1^g\right)^2E_{k,g}^z
\end{align*}
where the inequality is due to Assumption \ref{assumptionijsmooth} and Lemma \ref{R}.

The second term can be bounded as
\begin{align*}
\mathbb{E}\left[\left\|\nabla F_i\left(x_k, y_k\right)-\nabla F_i\left(w_{k, i}^x, w_{k, i}^y\right)\right\|^2\right]\leq (L^f)^2(E_{k,f}^x+E_{k,f}^y).
\end{align*}
Combining the above inequalities, we have
\begin{align*}
&E_k\left[\left\|D_x\left(x_k, y_k, z_k\right)-D_k^x\right\|^2\right] \\
\leq&
24\left(L_1^g\right)^2
\mathbb{E}\left[\left\|z_k-z^*\left(x_k\right)\right\|^2\right]
+24 R^2\left(L_2^g\right)^2 \mathbb{E}\left[\left\|y_k-y^*\left(x_k\right)\right\|^2\right]
\\&+8\left(L_2^g\right)^2 R^2E_{k,g}^x
+16\left(L_2^g\right)^2 R^2E_{k,g}^y
+16\left(L_1^g\right)^2E_{k,g}^z
+2(L^f)^2(E_{k,f}^x+E_{k,f}^y)\\
\leq& L_x'S_k+L_x''\left(\mathbb{E}\left[\left\|y_k-y^*(x_k)\right\|^ 2\right]+\mathbb{E}\left[\left\|z_k-z^*(x_k)\right\|^ 2\right]\right),
\end{align*}
where $$L_x'=\max\{16\left(L_2^g\right)^2 R^2,16\left(L_1^g\right)^2,2(L^f)^2\},\quad L_x''=\max\{24\left(L_1^g\right)^2,24 R^2\left(L_2^g\right)^2 \}.$$
Similarly, we can obtain the inequality for $\mathbb{E}\left[\left\|D_z(x_k,y_k,z_k)-D_k^z\right\|^2\right]$.
\end{proof}
\begin{lemma}\label{E}
 For the error between the iterates and the memories, we have the following inequalities:
\begin{eqnarray*}
E_{k+1,f}^x &\leq&\left(1-\frac{1}{2 n}\right) E_{k,f}^x +(2 n+1) \alpha_k^2 \mathbb{E}\left[\left\|v_k^x\right\|^2\right],\\
E_{k+1,g}^x  &\leq&\left(1-\frac{1}{2 m}\right)E_{k,g}^x +(2 m+1) \alpha_k^2 \mathbb{E}\left[\left\|v_k^x\right\|^2\right], \\
E_{k+1,f}^y  & \leq&\left(1-\frac{1}{2 n}\right)E_{k,f}^y +\beta^2_k \mathbb{E}\left\|D^y_k\right\|^2+2 n \beta^2_k \mathbb{E}\left[\left\|D_y\left(x_k, y_k, z_k\right)\right\|^2\right], \\
E_{k+1,g}^y  & \leq&\left(1-\frac{1}{2 m}\right)E_{k,g}^y+\beta^2_k \mathbb{E}\left\|D^y_k\right\|^2+2 m \beta^2_k \mathbb{E}\left[\left\|D_y\left(x_k, y_k, z_k\right)\right\|^2\right], \\
E_{k+1,f}^z & \leq&\left(1-\frac{1}{2 n}\right) E_{k,f}^z+\gamma^2_k \mathbb{E}\left\|D^z_k\right\|^2+2 n \gamma^2_k \mathbb{E}\left[\left\|D_z\left(x_k, y_k, z_k\right)\right\|^2\right], \\
E_{k+1,g}^z & \leq&\left(1-\frac{1}{2 m}\right) E_{k,g}^z+\gamma^2_k \mathbb{E}\left\|D^z_k\right\|^2+2 m \gamma^2_k\mathbb{E}\left[\left\|D_z\left(x_k, y_k, z_k\right)\right\|^2\right].
\end{eqnarray*}
\end{lemma}
\begin{proof}
According to the definition of $x_k^i$, we have
\begin{eqnarray*}
E_k\left[\left\|x_{k+1}-w_{k+1,i}^x\right\|^2\right]&=&\frac{1}{n} E_k\left[\left\|x_{k+1}-x_k\right\|^2\right]+\frac{n-1}{n} E_k\left[\left\|x_{k+1}-w_{k,i}^x\right\|^2\right] \\
& =&\frac{\alpha_k^2}{n} E\left[\left\|v_k^x\right\|^2\right]+\frac{n-1}{n} E_k\left[\left\|x_{k+1}-w_{k,i}^x\right\|^2\right].
\end{eqnarray*}
For the second term, we use the Young's inequality, then
\begin{eqnarray*}
E_k\left[\left\|x_{k+1}-w_{k,i}^x\right\|^2\right]&=&E_k\left[\left\|x_k-\alpha_k v_k^x- w_{k,i}^x\right\|^2\right] \\
& =&E_k\left[\left\|x_k-w_{k,i}^x\right\|^2+\alpha_k^2\left\|v_k^x\right\|^2-2 \alpha_k\left\langle v_k^x, x_k-w_{k,i}^x\right\rangle\right] \\
& \leq& E_k\left[\left\|x_k-w_{k,i}^x\right\|^2+\alpha_k^2\left\|v_k^x\right\|^2+\frac{\alpha_k}{2 n \alpha_k}\left\|x_k-w_{k,i}^x\right\|^2+2 n \alpha_k^2\left\|v_k^x\right\|^2\right].
\end{eqnarray*}
Thus, we obtain
\begin{eqnarray*}
&&E_k\left[\left\|x_{k+1}-w_{k+1,i}^x\right\|^2\right]\\&=&\frac{\alpha_k^2}{n} E_k\left[\left\|v_k^x\right\|^2\right]+\frac{n-1}{n}\left[\left( 1 + \frac { 1 } { 2 n } \right) E\left[\left\|x_k-w_{k,i}^x\right\|^2\right]+(2 n+1) \alpha_k^2 E_k\left[\left\|v_k^x\right\|^2\right]\right] \\
& =&\left(\frac{\alpha_k^2}{n}+\frac{(n-1)(2 n+1) \alpha_k^2}{n}\right) E_k\left[\left\|v_k^x\right\|^2\right] +\frac{n-1}{n}\left(1+\frac{1}{2 n}\right) E_k\left[\left\|x_k-w_{k,i}^x\right\|^2\right] \\
& \leqslant&(2 n+1) \alpha_k^2 E\left[\left\|v_k^x\right\|^2\right]+\left(1-\frac{1}{2 n}\right) E_k\left[\left\|x_k-w_{k,i}^x\right\|^2\right]. \\
\end{eqnarray*}
Taking the full expectation yields the desired inequality in the lemma. Similarly, we can obtain the result regarding $E_{k+1,g}^x$.
The proof for $E_{k+1,f}^y$, $E_{k+1,g}^y$, $E_{k+1,f}^z$ and $E_{k+1,g}^z$ can be found in Lemma C.5 of \cite{dagreou2022framework}.
\end{proof}

\begin{lemma}\label{S}
Suppose Assumption \ref{assump UL}, \ref{assump LL} and \ref{assumptionijsmooth} hold,
if $4\beta_k^2(L_1^g)^2+4\gamma_k^2L_z'\leq \tau/2$,
then
\begin{align*}
S_{k+1} \leq&\left(1-\frac{\tau}{2}\right) S_k+(P_1 \gamma_k^2+P_2\beta_k^2) \mathbb{E}\left[\left\|y_k-y^*\left( x_k\right)\right\|^{2}\right] +P_3\gamma_k^2 \mathbb{E}\left[\left\| z_k-z^{*} ( x_k) \right\|^2\right] +P_4 \alpha_k^2 \mathbb{E}\left[\left\|v_k^x\right\|^2\right],
\end{align*}
where
\begin{align*}
    P_1=(2(m+n)+4)L_z^2+4  L^{\prime \prime}_z,\quad P_2=(2(m+m)+4)\left(L_1^g\right)^2,\\
    P_3=(2(m+n)+4) L_z^2+4 L^{\prime \prime}_z,\quad P_4=2(m+n)+2.
\end{align*}
\end{lemma}
\begin{proof}
By adding the inequalities in Lemma \ref{E}, we obtain
\begin{eqnarray*}
S_{k+1} &\leq&(1-\tau) S_{k}+\mathbb{E}\left[2 \beta_k^2\left[\left\|D_k^y\right\|^2\right]+2 \gamma_k^2\left[\left\|D_k^z\right\|^2\right]\right) \\&& +2(m+n)\left(\beta_k^2 \mathbb{E}\left[\left\|D_y(x_k,y_k,z_k)\right\|^2\right]+\gamma_k^2 \mathbb{E}\left[\left\|D_z(x_k,y_k,z_k)\right\|^2\right]\right) \\
&& +2(m+n+1) \alpha_k^2 \mathbb{E}\left[\left\|v_k^x\right\|^2\right] \\
& \leq&(1-\tau) S_{k}+2 \beta_k^2\left(2 \mathbb{E}\left[\left\|D_y(x_k,y_k,z_k)\right\|^{2}\right]+2 \left(L_1^g\right)^2 S_{k}\right) \\
&&+2 \gamma_k^2\left(2 \mathbb{E}\left[\left\|D_z(x_k,y_k,z_k)\right\|^2\right]+2 L^{\prime}_z S_{k}+2 L^{\prime \prime}_z \mathbb{E}[\|y_k-y^*(x_k)\|^2]+2 L^{\prime \prime}_z \mathbb{E}[\|z_k-z^*(x_k)\|^2]\right) \\
&& +2(m+n)\left(\beta_k^2 \mathbb{E}\left[\| D_y(x_k,y_k,z_k) \|^2\right]+\gamma_k^2 \mathbb{E}\left[\left\|D_z(x_k,y_k,z_k)\right\|^2\right]\right)
\\&&+2(m+n+1) \alpha_k^2 \mathbb{E}\left[\left\|v_k^x\right\|^2\right] \\
& =&\left(1-\tau+4 \beta_k^2\left( L_1^g\right)^2+4 \gamma_k^2 L^{\prime}_z\right)S_{k}  +(2(m+n)+4) \beta_ k^2 \mathbb{E}\left[\left\|D_y(x_k,y_k,z_k)\right\|^2\right] \\
&& +(2(m+n)+4) \gamma_k^2 \mathbb{E}\left[\left\|D_z(x_k,y_k,z_k)\right\|^2\right]  +4 \gamma_k^2 L^{\prime \prime}_z\left(\mathbb{E}[\|y_k-y^*(x_k)\|^2]+\mathbb{E}[\|z_k-z^*(x_k)\|^2]\right) \\
&& +(2(m+n)+2) \alpha_k^2 \mathbb{E}\left[\left\|v_k^x\right\|^2\right]    \\
&\leq & \left(1-\tau+4 \beta_k^2\left(L_1^{g}\right)^2+4 \gamma_k^2 L^{\prime}_z\right) S_{k} +(2(m+n)+2)\alpha_k^2 \mathbb{E}\left[\left\| v_k^x\right\|^{2}\right]\\
&&+\left[(2(m+m)+4) \beta_k^2\left(L_1^g\right)^2+(2(m+n)+4) \gamma_k^2 L_z^2+4 \gamma_k^2 L^{\prime \prime}_z\right] \mathbb{E}\left[\| y_k-y^* ( x_k) \|^2\right] \\
& &+\left[(2(m+n)+4) \gamma_k^2 L_z^2+4 \gamma_k^2 L^{\prime \prime}_z\right] \mathbb{E}\left[\left\|z_k-z^*\left(x_k\right)\right\|^2\right],
\end{eqnarray*}
where the second and third inequalities use Lemma \ref{var} and Lemma \ref{dL}, respectively.

Suppose $4 \beta_k^2\left( L_1^g\right)^2+4 \gamma_k^2 L^{\prime}\leq \tau/2$, we have
\begin{eqnarray*}
S_{k+1}& \leq&\left(1-\frac{\tau}{2}\right) S_{k} +(2(m+n)+2)\alpha_k^2 \mathbb{E}\left[\left\| v_k^x\right\|^{2}\right]\\
&&+\left[(2(m+m)+4) \beta_k^2\left(L_1^g\right)^2+(2(m+n)+4) \gamma_k^2 L_z^2+4 \gamma_k^2 L^{\prime \prime}_z\right] \mathbb{E}\left[\| y_k-y^* ( x_k) \|^2\right] \\
& &+\left[(2(m+n)+4) \gamma_k^2 L_z^2+4 \gamma_k^2 L^{\prime \prime}_z\right] \mathbb{E}\left[\left\|z_k-z^*\left(x_k\right)\right\|^2\right].
\end{eqnarray*}
\end{proof}

\begin{theorem}\textbf{(Restatement of Theorem \ref{thsaba})}

Fix an iteration $K>1$ and assume that Assumptions \ref{assump UL} to \ref{assump LL} and \ref{assumptionijsmooth} hold.
Let the step sizes be $\alpha_k=c_1N^{-2/3}$, $\beta_k=c_2N^{-2/3}$, $\gamma_k=c_3N^{-2/3}$, $\rho_k=c_4 N^{-2/3}$.
Take $c_1$, $c_2$, $c_3$ and $c_4$ satisfy
\begin{align*}
    c_2&\leq \min\left\{\frac{\mu}{16 c''},\sqrt{\frac{c' }{16(L_1^g)^2}}\right\},\\
    c_3&\leq\min\left\{\sqrt{\frac{c'}{16L_z'}},\sqrt{\frac{\mu c_2}{16c''}},\frac{\mu}{16\Delta}c_2\right\},\\
    c_4&\leq \min\left\{\sqrt{\frac{\mu c_3}{8L_x''}},\sqrt{2}c_3,\frac{2}{3\mu}c_3,\frac{\mu }{12(L_1^g)^2}c_3\right\},\\
    c_1&\leq\min\left\{\frac{1}{32c''},\frac{1}{2L^H},\frac{\mu}{16L_{y^*}^2}c_2,\frac{\mu}{48L_{z^*}^2}c_3,\frac{1}{64(L^H)^2}c_4,2c_4\right\},
\end{align*}
where $c'=2$ and $c''=\max\left\{ 6L_z^2+4  L^{\prime \prime}_z,6\left(L_1^g\right)^2,4\right\}$
are constants that make
$\tau\leq c' N^{-1}$ and $P_1,P_2, P_3, P_4\leq c''N$
hold true, respectively.
Then the iterates in MA-SABA satisfy
\begin{eqnarray*}
\frac{1}{K} \sum_{k=0}^{K-1} \mathbb{E}\left[\left\|\nabla H\left(x_k\right)\right\|^2\right]=\mathcal{O}\left(N^{\frac{2}{3}}K^{-1}\right).
\end{eqnarray*}
\end{theorem}

\begin{proof}
First, we introduce the notation $N=n+m$ and set $c'=2$, $c''=\max\left\{ 6L_z^2+4  L^{\prime \prime}_z,6\left(L_1^g\right)^2,4\right\}$.
From Lemma E.5, it is known that
$\tau\leq c' N^{-1}$ and $P_1,P_2, P_3, P_4\leq c''N$
hold true (see the original text lines 1145, 1317-1318).

Then, we consider the Lyapunov function
\begin{eqnarray}
L_k=\mathbb{E}\left[H\left(x_k\right)\right]
+A\mathbb{E}\left[\left\|y_k-y^*\left(x_k\right)\right\|^2\right]
+B\mathbb{E}\left[\left\|z_k-z^*\left(x_k\right)\right\|^2\right]
+C\mathbb{E}[\|\nabla H(x_k)-v^x_k\|^2]+D S_k.
\end{eqnarray}
Using Lemma \ref{DH} Lemma \ref{Hlemma}, Lemma \ref{y,z}, Lemma \ref{vh} and Lemma \ref{S}, we get
\begin{align*}
 L_{k+1}-L_k=&\mathbb{E}\left[H\left(x_{k+1}\right)\right]-E\left[H\left(x_k\right)\right]
+A\left(\mathbb{E}\left[\left\|y_{k+1}-y^*\left(x_{k+1}\right)\right\|^2\right]-\mathbb{E}\left[\left\|y_k-y^*\left(x_k\right)\right\|^2\right]\right) \\
&+B\left(\mathbb{E}\left[\left\|z_{k+1}-z^*\left(x_{k+1}\right)\right\|^2\right]-\mathbb{E}\left[\left\|z_k-z^*\left(x_k\right)\right\|^2\right]\right)
\\&+C\left(\mathbb{E}[\|\nabla H(x_{k+1})-v^x_{k+1}\|^2]-\mathbb{E}[\|\nabla H(x_k)-v^x_k\|^2]\right)
\\&+D \left(S_{k+1}-S_k\right)
\\
\leq&
-\frac{\alpha_k}{2} \mathbb{E}\left[\left\|\nabla H\left(x_k\right)\right\|^2\right]-\frac{\tau}{2}D S_k
\\&
+\left(\frac{L^H\alpha_k^2}{2}-\frac{\alpha_k}{2}+A\frac{2 L_{y^*}^2 \alpha_k^2}{\beta_k \mu}+B\frac{3 L_{z^*}^2 \alpha_k^2}{\gamma_k \mu}
+ C\frac{2\left(L^H\right)^2 \alpha_k^2}{\rho_k}+P_4 \alpha_k^2D\right)
\mathbb{E}\left[\left\|v_k^x\right\|^2\right]
\\&+\left(-A\beta_k \mu +8\Delta B\gamma_k+6C\Delta\mu\rho_k+(P_1 \gamma_k^2+P_2\beta_k^2) D \right)
\mathbb{E}\left[\left\| y_k-y^* (x_k)\right\|^2\right]
\\&
+\left(-B\gamma_k\mu +6C\left(L_1^g\right)^2\rho_k+P_3\gamma_k^2 D\right)
\mathbb{E}\left[\left\|z_k-z^*\left(x_k\right)\right\|^2\right]
\\&
+\left(\frac{\alpha_k}{2}-C\rho_k\right)
\mathbb{E}\left[\left\|\nabla H\left(x_k\right)-v_k^x\right\|^2\right]\\
&+2A \beta_k^2\mathbb{E}\left[\left\|D_y(x_k,y_k,z_k)-D_k^y\right\|^2\right]
\\&+2B \gamma_k^2\mathbb{E}\left[\left\|D_z(x_k,y_k,z_k)-D_k^z\right\|^2\right]
\\&+C\rho_k^2 \mathbb{E}\left[\left\|D_k^x-\mathbb{E}\left[D_k^x\right]\right\|^2\right],
\end{align*}
For the variance terms in the above inequality, using Lemma \ref{var}, we have
\begin{align*}
 L_{k+1}-L_k
 \leq&
-\frac{\alpha_k}{2} \mathbb{E}\left[\left\|\nabla H\left(x_k\right)\right\|^2\right]\\
&+\left(-\frac{\tau}{2}D+2A\left(L_{1}^g\right)^2\beta_k^2+2B\gamma_k^2L_z^{\prime}+C\rho_k^2L_x^{\prime}\right) S_k
\\&
+\left(\frac{L^H\alpha_k^2}{2}-\frac{\alpha_k}{2}+A\frac{2 L_{y^*}^2 \alpha_k^2}{\beta_k \mu}+B\frac{3 L_{z^*}^2 \alpha_k^2}{\gamma_k \mu}
+ C\frac{2\left(L^H\right)^2 \alpha_k^2}{\rho_k}+P_4 \alpha_k^2D\right)
\mathbb{E}\left[\left\|v_k^x\right\|^2\right]
\\&+\left(-A\beta_k \mu +8\Delta B\gamma_k+6C\Delta\mu\rho_k+(P_1 \gamma_k^2+P_2\beta_k^2) D+2B\gamma_k^2L_z^{\prime\prime}+C\rho_k^2L_x^{\prime\prime} \right)
\mathbb{E}\left[\left\| y_k-y^* (x_k)\right\|^2\right]
\\&
+\left(-B\gamma_k\mu +6C\left(L_1^g\right)^2\rho_k+P_3\gamma_k^2 D+2B\gamma_k^2L_z^{\prime\prime}+C\rho_k^2L_x^{\prime\prime} \right)
\mathbb{E}\left[\left\|z_k-z^*\left(x_k\right)\right\|^2\right]
\\&
+\left(\frac{\alpha_k}{2}-C\rho_k\right)
\mathbb{E}\left[\left\|\nabla H\left(x_k\right)-v_k^x\right\|^2\right],
\end{align*}

We choose the coefficients of the Lyapunov function as
$A=1$, $B=1$, $C=1$, $D={N^{-1/3}},$
and the step sizes $\alpha_k=c_1N^{-2/3}$, $\beta_k=c_2N^{-2/3}$, $\gamma_k=c_3N^{-2/3}$, $\rho_k=c_4 N^{-2/3}$.

Based on our choice of
\begin{align*}
    c_2&\leq \min\left\{\frac{\mu}{16 c''},\sqrt{\frac{c' }{16(L_1^g)^2}}\right\},\\
    c_3&\leq\min\left\{\sqrt{\frac{c'}{16L_z'}},\sqrt{\frac{\mu c_2}{16c''}},\frac{\mu}{16\Delta}c_2\right\},\\
    c_4&\leq \min\left\{\sqrt{\frac{\mu c_3}{8L_x''}},\sqrt{2}c_3,\frac{2}{3\mu}c_3,\frac{\mu }{12(L_1^g)^2}c_3\right\},\\
    c_1&\leq\min\left\{\frac{1}{32c''},\frac{1}{2L^H},\frac{\mu}{16L_{y^*}^2}c_2,\frac{\mu}{48L_{z^*}^2}c_3,\frac{1}{64(L^H)^2}c_4,2c_4\right\},
\end{align*}
we proceed with the following derivation:

Since $c_2 \leq \frac{\mu}{16 c''}$ and $c_3\leq\sqrt{\frac{\mu c_2}{16c''}}$, it follows that $c_3\leq\sqrt{\frac{\mu c_2}{16c''}}\leq \frac{\mu }{16c''} \leq \frac{\mu }{4c''}$;

Since $c_3\leq\sqrt{\frac{\mu c_2}{16c''}}$ and $c''=\max\left\{ 6L_z^2+4  L^{\prime \prime}_z,6\left(L_1^g\right)^2,4\right\}\geq 4L^{\prime \prime}_z$, it follows that $c_3\leq\sqrt{\frac{\mu c_2}{16c''}}\leq\sqrt{\frac{\mu c_2}{16\cdot 4L^{\prime \prime}_z }}\leq\sqrt{\frac{\mu c_2}{32L_z''}}$;

Since $c_2 \leq \frac{\mu}{16 c''}$, $c_3\leq\sqrt{\frac{\mu c_2}{32L_z''}}$ and $c''\geq 4L^{\prime \prime}_z$, it follows that $c_3\leq\sqrt{\frac{\mu c_2}{32c''}}\leq \sqrt{\frac{\mu ^2}{32\cdot16\cdot L_z'' c''}}\leq \frac{\mu }{16L_z''}$;

Since $c_4\leq\min\{\sqrt{2}c_3,\frac{2}{3\mu}c_3\}$, $c_3 \leq \min\left\{\sqrt{\frac{c'}{16L_z'}}, \frac{\mu}{16\Delta}c_2, \sqrt{\frac{\mu c_2}{32L_z''}}\right\}$, $L_x'=L_z'$ and $L_x''=L_z''$, it follows that $c_4\leq \min\left\{\sqrt{\frac{c'}{8L_x'}},\frac{ c_2}{24\Delta},\sqrt{\frac{\mu c_2}{16 L_x''}}\right\}$.

Therefore, we have
\begin{align*}
    c_2&\leq \min\left\{\frac{\mu}{16 c''},\sqrt{\frac{c' }{16(L_1^g)^2}}\right\},\\
    c_3&\leq\min\left\{\frac{\mu }{16L_z''},\sqrt{\frac{c'}{16L_z'}},\frac{\mu}{4c''},\sqrt{\frac{\mu c_2}{32L_z''}},\sqrt{\frac{\mu c_2}{16c''}},\frac{\mu}{16\Delta}c_2\right\},\\
    c_4&\leq \min\left\{\sqrt{\frac{c'}{8L_x'}},\frac{ c_2}{24\Delta},\sqrt{\frac{\mu c_2}{16 L_x''}},\frac{\mu c_3}{12(L_1^g)^2},\sqrt{\frac{\mu c_3}{8L_x''}} \right\},\\
    c_1&\leq\min\left\{\frac{1}{32c''},\frac{1}{2L^H},\frac{\mu}{16L_{y^*}^2}c_2,\frac{\mu}{48L_{z^*}^2}c_3,\frac{1}{64(L^H)^2}c_4,2c_4\right\},
\end{align*}
it can be deduced that
\begin{align*}
    \alpha_k&\leq\min\left\{\frac{1}{2L^H},\frac{\mu}{16L_{y^*}^2}\beta_k,\frac{\mu}{48L_{z^*}^2}\gamma_k,\frac{1}{64(L^H)^2}\rho_k,2\rho_k\right\},\\
    \gamma_k&\leq\min\left\{\frac{\mu}{16\Delta}\beta_k,\frac{\mu}{16L_z''}\right\},
    \quad \gamma_k^2\leq \frac{\mu}{32L_z''}\beta_k,\\
    \rho_k&\leq \min\left\{\frac{1}{24\Delta}\beta_k,\frac{\mu}{12(L_1^g)^2}\gamma_k\right\},
    \quad \rho_k^2\leq \min\left\{\frac{\mu}{8L_x''}\gamma_k,\frac{\mu}{16 L_x''}\beta_k\right\}\\
    c_1 &\leq \frac{1}{32c''},
    \quad c_2\leq \sqrt{\frac{c'}{8\left(L_1^g\right)^2}},
    \quad c_3\leq \min\left\{\sqrt{\frac{c'}{16L_z'}},\frac{\mu}{4c''}\right\},
    \quad c_4\leq \sqrt{\frac{c'}{8L_x'}},
    \quad c''(c_2^2+c_3^2)\leq \frac{c_2\mu}{8},\\
    &4\beta_k^2(L_1^g)^2+4\gamma_k^2L_z'\leq \tau/2.
\end{align*}
 Then we have the following set of inequalities established:
\begin{align*}
\left\{\begin{aligned}
-\frac{\tau}{2}D+2A\left(L_{1}^g\right)^2\beta_k^2+2B\gamma_k^2L_z^{\prime}+C\rho_k^2L_x^{\prime}\leq 0, \\
\frac{L^H\alpha_k^2}{2}-\frac{\alpha_k}{2}+A\frac{2 L_{y^*}^2 \alpha_k^2}{\beta_k \mu}+B\frac{3 L_{z^*}^2 \alpha_k^2}{\gamma_k \mu}
+ C\frac{2\left(L^H\right)^2 \alpha_k^2}{\rho_k}+P_4 \alpha_k^2D \leq 0, \\
-A\beta_k \mu +8\Delta B\gamma_k+6C\Delta\mu\rho_k+(P_1 \gamma_k^2+P_2\beta_k^2) D+2B\gamma_k^2L_z^{\prime\prime}+C\rho_k^2L_x^{\prime\prime}  \leq 0, \\
-B\gamma_k\mu +6C\left(L_1^g\right)^2\rho_k+P_3\gamma_k^2 D+2B\gamma_k^2L_z^{\prime\prime}+C\rho_k^2L_x^{\prime\prime} \leq 0, \\
\frac{\alpha_k}{2}-C\rho_k \leq 0.
\end{aligned}\right.
\end{align*}
To make the proof more comprehensive, we will verify the validity of each inequality one by one.
\begin{itemize}
\item
\textbf{Inequality 1:}
\begin{align*}
&-\frac{\tau}{2}D+2A\left(L_{1}^g\right)^2\beta_k^2+2B\gamma_k^2L_z^{\prime}+C\rho_k^2L_x^{\prime}
\\&\leq -\frac{c^{\prime}}{2}N^{-1-\frac{1}{3}}+2\left(L_{1}^g\right)^2c_2^2N^{-\frac{4}{3}}+2c_3^2L_z^{\prime}N^{-\frac{4}{3}}+c_4^2L_x^{\prime}N^{-\frac{4}{3}}
\\&\leq
-\frac{c^{\prime}}{4}N^{-\frac{4}{3}}+2c_3^2L_z^{\prime}N^{-\frac{4}{3}}+c_4^2L_x^{\prime}N^{-\frac{4}{3}}
\\&\leq
-\frac{c^{\prime}}{8}N^{-\frac{4}{3}}+c_4^2L_x^{\prime}N^{-\frac{4}{3}}
\\&
\leq 0,
\end{align*}
where the justification for the four inequalities holding true are, respectively,
$\tau\leq c' N^{-1}$,
$c_2\leq \sqrt{\frac{c'}{8\left(L_1^g\right)^2}}$,
$c_3\leq \sqrt{\frac{c'}{16L_z'}}$,
and $c_4\leq \sqrt{\frac{c'}{8L_x'}}$.

\item
\textbf{Inequality 2:}
\begin{align*}
&\frac{L^H\alpha_k^2}{2}-\frac{\alpha_k}{2}+A\frac{2 L_{y^*}^2 \alpha_k^2}{\beta_k \mu}+B\frac{3 L_{z^*}^2 \alpha_k^2}{\gamma_k \mu}
+ C\frac{2\left(L^H\right)^2 \alpha_k^2}{\rho_k}+P_4 \alpha_k^2D\\
&=\frac{L^H\alpha_k^2}{2}-\frac{\alpha_k}{2}+\frac{2 L_{y^*}^2 \alpha_k^2}{\beta_k \mu}+\frac{3 L_{z^*}^2 \alpha_k^2}{\gamma_k \mu}
+ \frac{2\left(L^H\right)^2 \alpha_k^2}{\rho_k}+P_4 \alpha_k^2 N^{-\frac{1}{3}}\\
&\leq
-\frac{\alpha_k}{4}+\frac{2 L_{y^*}^2 \alpha_k^2}{\beta_k \mu}+\frac{3 L_{z^*}^2 \alpha_k^2}{\gamma_k \mu}
+ \frac{2\left(L^H\right)^2 \alpha_k^2}{\rho_k}+P_4 \alpha_k^2 N^{-\frac{1}{3}}\\
&\leq
-\frac{\alpha_k}{8}+\frac{3 L_{z^*}^2 \alpha_k^2}{\gamma_k \mu}
+ \frac{2\left(L^H\right)^2 \alpha_k^2}{\rho_k}+P_4 \alpha_k^2 N^{-\frac{1}{3}}\\
&\leq
-\frac{\alpha_k}{16}
+ \frac{2\left(L^H\right)^2 \alpha_k^2}{\rho_k}+P_4 \alpha_k^2 N^{-\frac{1}{3}}\\
&\leq
-\frac{\alpha_k}{32}+P_4 \alpha_k^2 N^{-\frac{1}{3}}\\
&\leq -\frac{c_1}{32}N^{-\frac{2}{3}}+c'' c_1^2 N^{1-\frac{1}{3}-\frac{4}{3}}\\
&\leq 0,
\end{align*}
where the justification for the six inequalities holding true are, respectively,
$\alpha_k \leq \frac{1}{2L^H}$,
$\alpha_k \leq  \frac{\mu}{16L_{y^*}^2}\beta_k$,
$\alpha_k \leq \frac{\mu}{48L_{z^*}^2}\gamma_k$,
$\alpha_k \leq \frac{1}{64(L^H)^2}\rho_k$,
$P_4\leq c''N$,
and
$c_1\leq \frac{1}{32c''}$.

\item
\textbf{Inequality 3:}
\begin{align*}
&-A\beta_k \mu +8\Delta B\gamma_k+6C\Delta\mu\rho_k+(P_1 \gamma_k^2+P_2\beta_k^2) D+2B\gamma_k^2L_z^{\prime\prime}+C\rho_k^2L_x^{\prime\prime}
\\
&=-\beta_k \mu +8\Delta \gamma_k+6\Delta\mu\rho_k+(P_1 \gamma_k^2+P_2\beta_k^2) D+2\gamma_k^2L_z^{\prime\prime}+\rho_k^2L_x^{\prime\prime}
\\&\leq
-\frac{\beta_k \mu}{2} +6\Delta\mu\rho_k+(P_1 \gamma_k^2+P_2\beta_k^2) D+2\gamma_k^2L_z^{\prime\prime}+\rho_k^2L_x^{\prime\prime}
\\&\leq
-\frac{\beta_k \mu}{4} +(P_1 \gamma_k^2+P_2\beta_k^2) D+2\gamma_k^2L_z^{\prime\prime}+\rho_k^2L_x^{\prime\prime}
\\&\leq
-\frac{\beta_k \mu}{8} +2\gamma_k^2L_z^{\prime\prime}+\rho_k^2L_x^{\prime\prime}
\\&\leq
-\frac{\beta_k \mu}{16} +\rho_k^2L_x^{\prime\prime}
\\&\leq 0
\end{align*}
where the justification for the five inequalities holding true are, respectively,
$\gamma_k\leq\frac{\mu}{16\Delta}\beta_k$,
$\rho_k\leq\frac{1}{24\Delta}\beta_k$
$c''(c_2^2+c_3^2)\leq \frac{c_2\mu}{8}$,
$\gamma_k^2\leq \frac{\mu}{32L_z''}\beta_k$,
and $\rho_k^2\leq \frac{\mu}{16 L_x''}\beta_k$.

To prevent any confusion, we additionally note that
the second inequality arises because $P_1,P_2\leq c''N$ and $c''(c_2^2+c_3^2)\leq \frac{c_2\mu}{4}$
ensure that
\begin{align*}
-\frac{\beta_k \mu}{4} +(P_1 \gamma_k^2+P_2\beta_k^2) D&\leq -\frac{c_2 \mu}{4}N^{-\frac{2}{3}} +(c'' c_3^2+c''c_2^2) N^{1-\frac{4}{3}-\frac{1}{3}}\\
&= \left(-\frac{c_2 \mu}{4} +(c'' c_3^2+c''c_2^2)\right) N^{-\frac{2}{3}}\\
&\leq -\frac{c_2 \mu}{8}N^{-\frac{2}{3}}\\
&= -\frac{\beta_k \mu}{8}.
\end{align*}
The condition
$c''(c_2^2+c_3^2)\leq \frac{c_2\mu}{8}$ is also reasonable. This can be achieved, for instance, by requiring that the coefficients of the step sizes adhere to
$c_2\leq \frac{\mu}{16 c''}$ and $c_3^2\leq \frac{\mu c_2}{16c''}$.
\item
\textbf{Inequality 4:}
\begin{align*}
&-B\gamma_k\mu +6C\left(L_1^g\right)^2\rho_k+P_3\gamma_k^2 D+2B\gamma_k^2L_z^{\prime\prime}+C\rho_k^2L_x^{\prime\prime}\\
&=-\gamma_k\mu +6\left(L_1^g\right)^2\rho_k+P_3\gamma_k^2 D+2\gamma_k^2L_z^{\prime\prime}+\rho_k^2L_x^{\prime\prime}
\\
&\leq-\frac{\gamma_k\mu}{2} +P_3\gamma_k^2 D+2\gamma_k^2L_z^{\prime\prime}+\rho_k^2L_x^{\prime\prime}
\\
&\leq-\frac{\gamma_k\mu}{4} +2\gamma_k^2L_z^{\prime\prime}+\rho_k^2L_x^{\prime\prime}
\\
&\leq-\frac{\gamma_k\mu}{8} +\rho_k^2L_x^{\prime\prime}
\\&\leq 0,
\end{align*}
where the justification for the four inequalities holding true are, respectively,
$\rho_k\leq\frac{\mu}{12(L_1^g)^2}\gamma_k$,
$c_3\leq\frac{\mu}{4c''}$,
$\gamma_k\leq\frac{\mu}{16L_z''}$,
and $\rho_k^2\leq \frac{\mu}{8L_x''}\gamma_k$.
For a complete proof, the detailed process by which the second inequality holds is as follows:
\begin{align*}
-\frac{\gamma_k\mu}{2} +P_3\gamma_k^2 D
&\leq-\frac{c_3\mu}{2}N^{-\frac{2}{3}} +c''c_3^2 N^{1-\frac{4}{3}-\frac{1}{3}}\\
&=\left(-\frac{c_3\mu}{2} +c''c_3^2 \right)N^{-\frac{2}{3}}\\
&\leq -\frac{c_3\mu}{4}N^{-\frac{2}{3}}
\\
&=-\frac{\gamma_k\mu}{4}.
\end{align*}
\item
\textbf{Inequality 5:}
Given that $C=1$ and $\alpha_k\leq2\rho_k$, we can affirm that the last inequality holds true, which is:
\begin{align*}
\frac{\alpha_k}{2}-C\rho_k=\frac{\alpha_k}{2}-\rho_k\leq 0.
\end{align*}
\end{itemize}
Up to this point, we've confirmed that each inequality in the system holds.

Consequently, the inequality of the difference in the Lyapunov function can be simplified to
\begin{align*}
 L_{k+1}-L_k
&\leq
-\frac{\alpha_k}{2} \mathbb{E}\left[\left\|\nabla H\left(x_k\right)\right\|^2\right]
\end{align*}
Summing and rearranging the above expressions yields
\begin{align*}
\frac{1}{K}\sum_{k=0}^{K-1}\mathbb{E}\left[\left\|\nabla H\left(x_k\right)\right\|^2\right]
&\leq \frac{L_0}{\alpha_kK}=\mathcal{O}\left(\frac{N^{2/3}}{K}\right).
\end{align*}
\end{proof}

\section{Proof of Theorems \ref{thpagefini}}\label{detial:proofpagefini}

\begin{corollary}\label{Hlemmabaised}
Suppose Assumptions \ref{assump UL} and \ref{assump LL} hold. Then we have
\begin{eqnarray*}
\mathbb{E}\left[H\left(x_{k+1}\right)\right]&\leq &\mathbb{E}\left[H\left(x_k\right)\right]
-\frac{\alpha_k}{2} \mathbb{E}\left[\left\|\nabla H\left(x_k\right)\right\|^2\right]
+
 \left(\frac{L^H\alpha_k^2}{2}-\frac{\alpha_k}{2}\right) \mathbb{E}\left[\left\|v_k^x\right\|^2\right]
 +{\alpha_k} \mathbb{E}\left[\left\|D_x(x_k,y_k,z_k)-v_k^x\right\|^2\right]
\\&&+3{\alpha_k} \left(\left(L^f\right)^2+\left(L^g_2 R\right)^2\right) \mathbb{E}\left[\left\|y_k-y^*\left(x_k\right)\right\|^2\right]
+3{\alpha_k} \left(L^g_1\right)^2 \mathbb{E}\left[\left\|z_k-z^*\left(x_k\right)\right\|^2\right]
\end{eqnarray*}
\end{corollary}
\begin{proof}
    By combining Lemmas \ref{Hlemma} and \ref{DH}, the proof can be established.
\end{proof}
\begin{lemma}\label{stronglysmooth}
    If $\phi$ is $\alpha$-strongly convex and $\beta$-smooth, then $$\langle\nabla\phi(x)-\nabla\phi(y),x-y\rangle\geq \frac{\alpha\beta}{\alpha+\beta}\|x-y\|^2+\frac{1}{\alpha+\beta}\|\nabla\phi(x)-\nabla\phi(y)\|^2.$$
\end{lemma}
\begin{proof}
    See Lemma C.2. in \cite{khanduri2021near}.
\end{proof}
\begin{lemma}\label{eg30}
Suppose Assumption \ref{assump UL} and \ref{assump LL} hold and the step sizes  satisfy
\begin{align*}
\beta_k,\gamma_k\leq\min\left\{\frac{\mu+L_1^g}{\mu L_1^g},\frac{1}{\mu+L_1^g}\right\}.
\end{align*}
Then we have
\begin{align*}
\mathbb{E}\left[\left\|y_{k+1}-y^*\left(x_{k+1}\right)\right\|^2\right]
\leq&
\left(1-\frac{\mu L_1^g\beta_k}{2(\mu+L_1^g)}\right)\mathbb{E}\left[\left\|y_k-y^*\left(x_k\right)\right\|^2\right]-\frac{1}{\mu+L_1^g}\beta_k\mathbb{E}\left[\left\| D_y(x_k,y_k,z_k)\right\|^2\right]\\
&+6\frac{\mu+L_1^g}{\mu L_1^g}\beta_k\mathbb{E}\left[\| v_k^y-D_y(x_k,y_k,z_k)\|^2\right]\\
&+ \frac{2(\mu+L_1^g)L_{y^*}^2\alpha_k^2}{\mu L_1^g\beta_k}\mathbb{E}\left[\left\|v_k^x\right\|^2\right].
\end{align*}
\begin{align*}
	\mathbb{E}\left[\left\|z_{k+1}-z^*\left(x_{k+1}\right)\right\|^2\right]
	\leq&
	\left(1-\frac{\gamma_k\mu}{4}\right)\mathbb{E}\left[\left\|z_k-z^*\left(x_k\right)\right\|^2\right]
	+
	\frac{18((L_2^gR)^2+(L^f)^2)\gamma_k}{\mu}\mathbb{E}\left[\left \|y_k-y^*(x_k)\right\|^2\right]
	\\
	&+\frac{12\gamma_k}{\mu}\mathbb{E}\left[\| v_k^z-D_z(x_k,y_k,z_k)\|^2\right]
	+ \frac{4L_{z^*}^2\alpha_k^2}{\mu\gamma_k }\mathbb{E}\left[\left\|v_k^x\right\|^2\right].
\end{align*}

\end{lemma}

The proof of this lemma is similar to that of Lemma \ref{y,z}. The main difference is that $v^{\cdot}_k$ is no longer an unbiased estimate of $D_{\cdot}(x_k,y_k,z_k)$. Below, we present the specific proof process.
\begin{proof}

\textbf{Inequality for $y$.}

By utilizing the Young's inequality and the $L_{y^*}$-Lipschitz continuity of $y^*(x)$, we have
\begin{eqnarray*}
\left\|y_{k+1}-y^*\left(x_{k+1}\right)\right\|^2
&= & \left\|y_{k+1}-y^*\left(x_k\right)+y^*\left(x_k\right)-y^*\left(x_{k+1}\right)\right\|^2 \\
&\leq & \left(1+\delta_k\right)\left\|y_{k+1}-y^*\left(x_k\right)\right\|^2+\left(1+\frac{1}{\delta_k}\right)\left\|y^*\left(x_k\right)-y^*\left(x_{k+1}\right)\right\|^2
\\ & \leq & \left(1+\delta_k\right)\left\|y_{k+1}-y^*\left(x_k\right)\right\|^2
+ \left(1+\frac{1}{\delta_k}\right)L_{y^*}^2\alpha_k^2\left\|v_k^x\right\|^2
\end{eqnarray*}
Taking the expectation conditionally on $x_k$, $y_k$, $z_k$ yields
\begin{eqnarray}\label{y11}
E_k\left[\left\|y_{k+1}-y^*\left(x_{k+1}\right)\right\|^2\right]
&\leq&
\left(1+\delta_k\right)E_k\left[\left\|y_{k+1}-y^*\left(x_k\right)\right\|^2\right]
+\left(1+\frac{1}{\delta_k}\right)L_{y^*}^2\alpha_k^2E_k[\left\|v_k^x\right\|^2].
\end{eqnarray}
For the first term, once again employing Young's inequality, we have
\begin{eqnarray*}
E_k\left[\left\|y_{k+1}-y^*\left(x_k\right)\right\|^2\right]
& =&E_k\left[\left\|y_k-y^*\left(x_k\right)-\beta_k v_k^y\right\|^2\right]
 \\
& =&E_k\left[\left\|y_k-\beta_k D_y(x_k,y_k,z_k)-y^*\left(x_k\right)-\beta_k\left(v_k^y-D_y(x_k,y_k,z_k)\right)\right\|^2\right] \\
& \leq &\left(1+\frac{\delta_k}{2}\right)E_k\left[\left\|y_k-\beta_k D_y(x_k,y_k,z_k)-y^*\left(x_k\right)\right\|^2\right]
\\
&&+\left(1+\frac{2}{\delta_k}\right)E_k\left[\| \beta_k\left(v_k^y-D_y(x_k,y_k,z_k)\right) \|^2\right],
\end{eqnarray*}
Utilizing Lemma \ref{stronglysmooth},
we can thus establish the following inequality
\begin{align*}
    E_k\left[\left\|y_k-\beta_k D_y(x_k,y_k,z_k)-y^*\left(x_k\right)\right\|^2\right]
    =&E_k\left[\left\|y_k-y^*\left(x_k\right)\right\|^2\right]
    +E_k\left[\left\|\beta_k D_y(x_k,y_k,z_k)\right\|^2\right]
    \\&
    -2E_k\beta_k\langle D_y(x_k,y_k,z_k),y_k-y^*\left(x_k\right)\rangle
    \\\leq&
    \left(1-2\beta_k\frac{\mu L_1^g}{\mu+L_1^g}\right)E_k\left[\left\|y_k-y^*\left(x_k\right)\right\|^2\right]\\
    &+\left(\beta_k^2-2\beta_k\frac{1}{\mu+L_1^g}\right)E_k\left[\left\| D_y(x_k,y_k,z_k)\right\|^2\right].
\end{align*}

Plugging it into (\ref{y11}) and taking the total expectation, we have
\begin{align*}
\mathbb{E}\left[\left\|y_{k+1}-y^*\left(x_{k+1}\right)\right\|^2\right]
\leq&
\left(1+\delta_k\right)\left(1+\frac{\delta_k}{2}\right)\left(1-2\beta_k\frac{\mu L_1^g}{\mu+L_1^g}\right)\mathbb{E}\left[\left\|y_k-y^*\left(x_k\right)\right\|^2\right]\\
&+
\left(1+\delta_k\right)\left(1+\frac{\delta_k}{2}\right)\left(\beta_k^2-2\beta_k\frac{1}{\mu+L_1^g}\right)\mathbb{E}\left[\left\| D_y(x_k,y_k,z_k)\right\|^2\right]\\
&+\left(1+\delta_k\right)\left(1+\frac{2}{\delta_k}\right)\beta_k^2\mathbb{E}\left[\| v_k^y-D_y(x_k,y_k,z_k)\|^2\right]\\
&+ \left(1+\frac{1}{\delta_k }\right)L_{y^*}^2\alpha_k^2\mathbb{E}\left[\left\|v_k^x\right\|^2\right],
\end{align*}
We choose the parameter $\delta_k$ and the step size $\beta_k$ to satisfy
\begin{align*}
    \delta_k=\frac{\mu L_1^g}{\mu+L_1^g}\beta_k,
    \quad\beta_k\leq\min\left\{\frac{\mu+L_1^g}{\mu L_1^g},\frac{1}{\mu+L_1^g}\right\},
\end{align*}
Consequently, the lemma concerning $y$ is proven.

 \textbf{Inequality for $z$.}

Based on the definition of $z_{k+1}$ and utilizing Young's inequality twice, we obtain
\begin{eqnarray*}
	\left\|z_{k+1}-z^*\left(x_{k+1}\right)\right\|^2
	& \leq & \left(1+\delta_k'\right)\left\|z_{k+1}-z^*\left(x_k\right)\right\|^2
	+ \left(1+\frac{1}{\delta_k'}\right)L_{z^*}^2\alpha_k^2\left\|v_k^x\right\|^2.
\end{eqnarray*}
Using this nonexpansiveness of projection and applying Young's inequality, we have
\begin{align*}
	\left\|z_{k+1}-z^*\left(x_k\right)\right\|^2 
	&\leq \left\|z_k-\gamma_k v_k^z-z^*\left(x_k\right)\right\|^2
	\\
	& \leq\left(1+\frac{\mu\gamma_k}{4}\right)\left\|z_k-\gamma_k D_z(x_k,y_k,z_k)-z^*\left(x_k\right)\right\|^2
	+\left(1+\frac{4}{\mu\gamma_k}\right)\| \gamma_k\left(v_k^z-D_z(x_k,y_k,z_k)\right) \|^2. 
\end{align*}
For the first term, we use the identity
\(\nabla_{22}^2 g(x_k,y^*(x_k))z^*(x_k)-\nabla_2 f(x_k,y^*(x_k))=0\)
to obtain the decomposition
\begin{align*}
	&\left\|z_k-\gamma_k D_z(x_k,y_k,z_k)-z^*\left(x_k\right)\right\|^2\notag
	= \Big\|
	(I-\gamma_k\nabla_{22}^2g(x_k,y_{k}))(z_k-z^*(x_k))
	\notag\\
	&\quad\quad\quad+\gamma_k\big((\nabla_{22}^2g(x_k,y^*(x_k))-\nabla_{22}^2g(x_k,y_k))z^*(x_k)+\nabla_2f(x_k,y_k)-\nabla_2f(x_k,y^*(x_k))\big) \Big\|^2.\notag
\end{align*}
Applying Young's inequality yields
\begin{align*}
	\left\|z_k-\gamma_k D_z(x_k,y_k,z_k)-z^*\left(x_k\right)\right\|^2
	\leq &\left(1+\frac{\gamma_k\mu}{2}\right)\left\|(I-\gamma_k\nabla_{22}^2g(x_k,y_{k}))(z_k-z^*(x_k))\right\|^2\notag
	\\
	&+\left(2+\frac{4}{\gamma_k\mu}\right)\gamma_k^2 \|\nabla_{22}^2g(x_k,y^*(x_k))-\nabla_{22}^2g(x_k,y_k)\|^2\|z^*(x_k)\|^2\notag
	\\&+\left(2+\frac{4}{\gamma_k\mu}\right)\gamma_k^2 \|\nabla_2f(x_k,y_k)-\nabla_2f(x_k,y^*(x_k))\|^2. \notag
\end{align*}
Using Assumption \ref{assump LL}(2) and the step size condition \(\gamma_k\le 1/L_1^g\), we have
\[\bigl\|I-\gamma_k\nabla_{22}^2 g(x_k,y_k)\bigr\|^2\le (1-\gamma_k\mu)^2.\]
Moreover, by Assumptions~\ref{assump UL}(1) and~\ref{assump LL}(1), together with \(\|z^*(x_k)\|\le R\), we obtain
\begin{align*}
	&\left\|z_k-\gamma_k D_z(x_k,y_k,z_k)-z^*\left(x_k\right)\right\|^2\notag\\
	\leq &
	\left(1+\frac{\gamma_k\mu}{2}\right)(1-\gamma_k\mu)^2\left\|z_k-z^*(x_k)\right\|^2
	{+\left(2+\frac{4}{\gamma_k\mu}\right)\gamma_k^2c_1\|y_k-y^*(x_k)\|^2},
\end{align*}
where $c_1=\left(L_2^gR\right)^2+(L^f)^2$.

Rearranging the above inequalities and taking the total expectation yields
\begin{align*}
	\mathbb{E}\left[\left\|z_{k+1}-z^*\left(x_{k+1}\right)\right\|^2\right]
	\leq&\left(1+\frac{\mu\gamma_k}{4}\right)\left(1+\frac{\gamma_k\mu}{2}\right)^2(1-\gamma_k\mu)^2\mathbb{E}\left[\left\|z_k-z^*\left(x_k\right)\right\|^2\right]\\
	&+
	\left(1+\frac{\mu\gamma_k}{4}\right)\left(1+\frac{\gamma_k\mu}{2}\right)\left(2+\frac{4}{\gamma_k\mu}\right)\gamma_k^2c_1\mathbb{E}\left[\left \|y_k-y^*(x_k)\right\|^2\right]\\
	&+\left(1+\frac{\mu\gamma_k}{2}\right)\left(1+\frac{4}{\mu\gamma_k}\right)\gamma_k^2\mathbb{E}\left[\| v_k^z-D_z(x_k,y_k,z_k)\|^2\right]\\
	&+ \left(1+\frac{2}{\mu\gamma_k }\right)L_{z^*}^2\alpha_k^2\mathbb{E}\left[\left\|v_k^x\right\|^2\right]. 
\end{align*}
Using  the step size condition $\gamma_k\leq\frac{1}{L_1^g+\mu},$ we have
\begin{align*}
	\mathbb{E}\left[\left\|z_{k+1}-z^*\left(x_{k+1}\right)\right\|^2\right]
	\leq&
	\left(1-\frac{\gamma_k\mu}{4}\right)\mathbb{E}\left[\left\|z_k-z^*\left(x_k\right)\right\|^2\right]
	+
	\frac{18\left(\left(L_2^gR\right)^2+(L^f)^2\right)\gamma_k}{\mu}\mathbb{E}\left[\left \|y_k-y^*(x_k)\right\|^2\right]
	\\
	&+\frac{12\gamma_k}{\mu}\mathbb{E}\left[\| v_k^z-D_z(x_k,y_k,z_k)\|^2\right]
	+ \frac{4L_{z^*}^2\alpha_k^2}{\mu\gamma_k }\mathbb{E}\left[\left\|v_k^x\right\|^2\right].
\end{align*}
\end{proof}

\begin{lemma}\label{pagefinitemoment}
    Under the Assumption \ref{assump UL}, \ref{assump LL}, \ref{assumporacle} and \ref{assumptionijsmooth}, We have the following inequalities established:
\begin{align*}
   (1)\quad  \mathbb{E}\left[\left\|v_{k+1}^y-D_y\left(x_{k+1}, y_{k+1}, z_{k+1}\right)\right\|^2\right]
    \leq&(1-p)\mathbb{E}\left[\left\|v_{k}^{y}-D_y(x_k,y_k,z_k)\right\|^2\right]+\frac{(1-p)}{b}(L_1^g)^2\alpha_k^2\mathbb{E}\left[\left\|v_k^x\right\|^2\right]
         \\&+\frac{2(1-p)}{b}(L_1^g)^2\beta_k^2\mathbb{E}\left[\left\|v_k^y-D_y(x_k,y_k,z_k)\right\|^2\right]
         \\&+\frac{2(1-p)}{b}(L_1^g)^4\beta_k^2 \mathbb{E}\left[\| y_k-y^*\left(x_k\right) \|^2\right].
\end{align*}
\begin{align*}
     (2)\quad \mathbb{E}\left[\left\|v_{k+1}^x-D_x\left(x_{k+1}, y_{k+1}, z_{k+1}\right)\right\|^2\right]
     &{\leq}(1-p)\mathbb{E}\left[\left\|v_{k}^{x}-D_x(x_k,y_k,z_k)\right\|^2\right]
    \\& +\left(2(L^f)^2+4R^2(L_2^g)^2\right)\frac{(1-p)}{b}\alpha_k^2\mathbb{E}\left[\left\|v_k^x\right\|^2\right]
    \\&+\left(2(L^f)^2+4R^2(L_2^g)^2\right)\frac{(1-p)}{b}\beta_k^2\mathbb{E}\left[\left\|v_k^y-D_y(x_k,y_k,z_k)\right\|^2\right]
    \\&+\left(2(L^f)^2+4R^2(L_2^g)^2\right)\frac{(1-p)}{b}\beta_k^2\left(L^g_1\right)^2 \mathbb{E}\left[\| y_k-y^*\left(x_k\right) \|^2\right]
     \\&
     +\frac{4(1-p)}{b}(L^g_2)^2\gamma_k^2\mathbb{E}\left[\left\|v^z_{k}-D_z(x_k,y_k,z_k)\right\|^2 \right]
      \\&+\frac{4(1-p)}{b}(L^g_2)^2\gamma_k^2L_z^2\mathbb{E}\left[\left\|z_k-z^*\left(x_k\right)\right\|^2\right]
     \\&+\frac{4(1-p)}{b}(L^g_2)^2\gamma_k^2 L_z^2\mathbb{E}\left[\left\|y_k-y^*\left(x_k\right)\right\|^2\right]
\end{align*}
\begin{align*}
    (3)\quad \mathbb{E}\left[\left\|v_{k+1}^z-D_z\left(x_{k+1}, y_{k+1}, z_{k+1}\right)\right\|^2\right]
     \leq
     &(1-p+\frac{4(1-p)}{b}(L^g_2)^2\gamma_k^2)\mathbb{E}\left[\left\|v_{k}^{z}-D_z(x_k,y_k,z_k)\right\|^2\right]
     \\&
     +\left(2(L^f)^2+4R^2(L_2^g)^2\right)\frac{(1-p)}{b}\alpha_k^2\mathbb{E}\left[\left\|v_k^x\right\|^2\right]
    \\&+\left(2(L^f)^2+4R^2(L_2^g)^2\right)\frac{(1-p)}{b}\beta_k^2\mathbb{E}\left[\left\|v_k^y-D_y(x_k,y_k,z_k)\right\|^2\right]
    \\&+\left(2(L^f)^2+4R^2(L_2^g)^2\right)\frac{(1-p)}{b}\beta_k^2\left(L^g_1\right)^2 \mathbb{E}\left[\| y_k-y^*\left(x_k\right) \|^2\right]
     \\&
     +\frac{4(1-p)}{b}(L^g_2)^2\gamma_k^2L_z^2\mathbb{E}\left[\left\|z_k-z^*\left(x_k\right)\right\|^2\right]
     \\&+\frac{4(1-p)}{b}(L^g_2)^2\gamma_k^2 L_z^2\mathbb{E}\left[\left\|y_k-y^*\left(x_k\right)\right\|^2\right]
\end{align*}
\end{lemma}
\begin{proof}
    \textbf{Proof of (1).}

By the definition of $v_{k+1}^y$, we have
    \begin{align*}
         &\mathbb{E}\left[\left\|v_{k+1}^y-D_y\left(x_{k+1}, y_{k+1}, z_{k+1}\right)\right\|^2\right]
         \\=&p\mathbb{E}\left[\left\|\frac{1}{m} \sum_{j \in [m]} \nabla_2 G_j\left(x_{k+1},y_{k+1}\right)-D_y\left(x_{k+1}, y_{k+1}, z_{k+1}\right)\right\|^2\right]
         \\&+(1-p)\mathbb{E}\left[\left\|v_{k}^{y}+\frac{1}{b} \sum_{j \in J}\left[\nabla_2 G_j\left(x_{k+1},y_{k+1}\right)-\nabla_2 G_j\left(x_{k},y_{k}\right)\right]-D_y\left(x_{k+1}, y_{k+1}, z_{k+1}\right)\right\|^2\right]
         \\=&(1-p)\mathbb{E}\left[\left\|v_{k}^{y}+\frac{1}{b} \sum_{j \in J}\left[\nabla_2 G_j\left(x_{k+1},y_{k+1}\right)-\nabla_2 G_j\left(x_{k},y_{k}\right)\right]-D_y\left(x_{k+1}, y_{k+1}, z_{k+1}\right)\right\|^2\right]\\
         =&(1-p)\mathbb{E}\left[\left\|v_{k}^{y}-D_y(x_k,y_k,z_k)\right\|^2\right]\\
         &+(1-p)\mathbb{E}\left[\left\|\frac{1}{b} \sum_{j \in J}\left(\nabla_2 G_j\left(x_{k+1},y_{k+1}\right)-D_y\left(x_{k+1}, y_{k+1}, z_{k+1}\right)+D_y(x_k,y_k,z_k)-\nabla_2 G_j\left(x_{k},y_{k}\right)\right)\right\|^2\right]
         \\
         \leq&(1-p)\mathbb{E}\left[\left\|v_{k}^{y}-D_y(x_k,y_k,z_k)\right\|^2\right]
         +\frac{1-p}{b}\mathbb{E}\left[\left\|\nabla_2 G_j\left(x_{k+1},y_{k+1}\right)-\nabla_2 G_j\left(x_{k},y_{k}\right)\right\|^2\right],
         \end{align*}
where the last equation uses the fact that $$\mathbb{E}\left[\frac{1}{b} \sum_{j \in J}\nabla_2 G_j\left(x_{k+1},y_{k+1}\right)\right]=D_y\left(x_{k+1}, y_{k+1}, z_{k+1}\right),\quad\mathbb{E}\left[\frac{1}{b} \sum_{j \in J}\nabla_2 G_j\left(x_{k},y_{k}\right)\right]=D_y\left(x_{k}, y_{k}, z_{k}\right).$$
The final inequality arises due to $\mathbb{E}\|X-\mathbb{E}[X]\|^2\leq \mathbb{E}[X^2]$.
Additionally, utilizing Assumption \ref{assumptionijsmooth}, we obtain
         \begin{align*}
          &\mathbb{E}\left[\left\|v_{k+1}^y-D_y\left(x_{k+1}, y_{k+1}, z_{k+1}\right)\right\|^2\right]
         \\
        {\leq}&(1-p)\mathbb{E}\left[\left\|v_{k}^{y}-D_y(x_k,y_k,z_k)\right\|^2\right]+\frac{(1-p)}{b}(L_1^g)^2\left(\alpha_k^2\mathbb{E}\left[\left\|v_k^x\right\|^2\right]+\beta_k^2\mathbb{E}\left[\left\|v_k^y\right\|^2\right]\right)
          \\
         \leq&(1-p)\mathbb{E}\left[\left\|v_{k}^{y}-D_y(x_k,y_k,z_k)\right\|^2\right]+\frac{(1-p)}{b}(L_1^g)^2\alpha_k^2\mathbb{E}\left[\left\|v_k^x\right\|^2\right]
         \\&+\frac{2(1-p)}{b}(L_1^g)^2\beta_k^2\mathbb{E}\left[\left\|v_k^y-D_y(x_k,y_k,z_k)\right\|^2\right]
         +\frac{2(1-p)}{b}(L_1^g)^4\beta_k^2 \mathbb{E}\left[\| y_k-y^*\left(x_k\right) \|^2\right].
    \end{align*}

\textbf{Proof of (2) and (3).}
Similarly, from the definition of $v_{k+1}^x$, we have
\begin{align*}
    &\mathbb{E}\left[\left\|v_{k+1}^x-D_x\left(x_{k+1}, y_{k+1}, z_{k+1}\right)\right\|^2\right]\\
    =&(1-p)\mathbb{E}\left[\left\|v_{k}^{x}+\frac{1}{b} \sum_{i \in I}\left( \nabla_{1} F_i\left(x_{k+1},y_{k+1}\right)-\nabla_{1} F_i\left(x_{k},y_{k}\right)\right)\right.\right.\\&\left. \left.- \frac{1}{b} \sum_{j \in J}\left( \nabla^2_{12} G_j\left(x_{k+1},y_{k+1}\right)z_{k+1-} \nabla^2_{12} G_j\left(x_{k},y_{k}\right)z_{k}\right)-D_x\left(x_{k+1}, y_{k+1}, z_{k+1}\right)\right\|^2\right].
    \end{align*}
Based on the fact that ${E}_k\left[\frac{1}{b}\sum_{i \in I}\nabla_{1} F_i\left(x_{k+1},y_{k+1}\right)
-\frac{1}{b}
     \sum_{j \in J} \nabla^2_{12} G_j\left(x_{k+1},y_{k+1}\right)z_{k+1}\right]=D_x(x_{k+1},y_{k+1},z_{k+1}),$ and
     ${E}_k\left[\frac{1}{b}\sum_{i \in I}\nabla_{1} F_i\left(x_{k},y_{k}\right)-\frac{1}{b}
     \sum_{j \in J} \nabla^2_{12} G_j\left(x_{k},y_{k}\right)z_{k}\right]=D_x(x_{k},y_{k},z_{k}),$
     we deduce
    \begin{align*}
     &\mathbb{E}\left[\left\|v_{k+1}^x-D_x\left(x_{k+1}, y_{k+1}, z_{k+1}\right)\right\|^2\right]\\
     \leq&(1-p)\mathbb{E}\left[\left\|v_{k}^{x}-D_x(x_k,y_k,z_k)\right\|^2\right]
     +2(1-p)\mathbb{E}\left[\left\|\frac{1}{b}\sum_{i \in I}\left[\nabla_{1} F_i\left(x_{k+1},y_{k+1}\right)-\nabla_{1} F_i\left(x_{k},y_{k}\right)\right]\right\|^2\right]
     \\&+2(1-p)\mathbb{E}\left[\left\|\frac{1}{b}
     \sum_{j \in J} [\nabla^2_{12} G_j\left(x_{k+1},y_{k+1}\right)z_{k+1}- \nabla^2_{12} G_j\left(x_{k},y_{k}\right)z_{k}]\right\|^2\right]
      \\
    \leq&(1-p)\mathbb{E}\left[\left\|v_{k}^{x}-D_x(x_k,y_k,z_k)\right\|^2\right]
     +\frac{2(1-p)}{b}\mathbb{E}\left[\left\| \nabla_{1} F_i\left(x_{k+1},y_{k+1}\right)-\nabla_{1} F_i\left(x_{k},y_{k}\right)\right\|^2\right]
     \\&+\frac{4(1-p)}{b}\mathbb{E}\left[\left\| \nabla^2_{12} G_j\left(x_{k+1},y_{k+1}\right)z_{k+1}- \nabla^2_{12} G_j\left(x_{k},y_{k}\right)z_{k+1}\right\|^2\right]
     \\&+\frac{4(1-p)}{b}\mathbb{E}\left[\left\|\nabla^2_{12} G_j\left(x_{k},y_{k}\right)z_{k+1}- \nabla^2_{12} G_j\left(x_{k},y_{k}\right)z_{k}\right\|^2\right].
\end{align*}

Under Assumption \ref{assumptionijsmooth}, it further implies that
\begin{align*}
&\mathbb{E}\left[\left\|v_{k+1}^x-D_x\left(x_{k+1}, y_{k+1}, z_{k+1}\right)\right\|^2\right]
    \\
     \leq&(1-p)\mathbb{E}\left[\left\|v_{k}^{x}-D_x(x_k,y_k,z_k)\right\|^2\right]
    +\frac{2(1-p)}{b}(L^f)^2\left(\alpha_k^2\mathbb{E}\left[\left\|v_k^x\right\|^2\right]+\beta_k^2\mathbb{E}\left[\left\|v_k^y\right\|^2\right]\right)
     \\&+\frac{4(1-p)}{b}R^2(L_2^g)^2\left(\alpha_k^2\mathbb{E}\left[\left\|v_k^x\right\|^2\right]+\beta_k^2\mathbb{E}\left[\left\|v_k^y\right\|^2\right]\right)
     +\frac{4(1-p)}{b}(L^g_2)^2\gamma_k^2\mathbb{E}\left[\left\|v^z_{k}\right\|^2 \right]
       \\
    =&(1-p)\mathbb{E}\left[\left\|v_{k}^{x}-D_x(x_k,y_k,z_k)\right\|^2\right]
+\left(\frac{2(1-p)}{b}(L^f)^2+\frac{4(1-p)}{b}R^2(L_2^g)^2\right)\alpha_k^2\mathbb{E}\left[\left\|v_k^x\right\|^2\right]
    \\&+\left(\frac{2(1-p)}{b}(L^f)^2+\frac{4(1-p)}{b}R^2(L_2^g)^2\right)\beta_k^2\mathbb{E}\left[\left\|v_k^y\right\|^2\right]
     +\frac{4(1-p)}{b}(L^g_2)^2\gamma_k^2\mathbb{E}\left[\left\|v^z_{k}\right\|^2 \right]
      \\
     {\leq}&(1-p)\mathbb{E}\left[\left\|v_{k}^{x}-D_x(x_k,y_k,z_k)\right\|^2\right]
     +\left(\frac{2(1-p)}{b}(L^f)^2+\frac{4(1-p)}{b}R^2(L_2^g)^2\right)\alpha_k^2\mathbb{E}\left[\left\|v_k^x\right\|^2\right]
    \\&+\left(\frac{2(1-p)}{b}(L^f)^2+\frac{4(1-p)}{b}R^2(L_2^g)^2\right)\beta_k^2\mathbb{E}\left[\left\|v_k^y-D_y(x_k,y_k,z_k)\right\|^2\right]
    \\&+\left(\frac{2(1-p)}{b}(L^f)^2+\frac{4(1-p)}{b}R^2(L_2^g)^2\right)\beta_k^2\left(L^g_1\right)^2 \mathbb{E}\left[\| y_k-y^*\left(x_k\right) \|^2\right]
     \\&
     +\frac{4(1-p)}{b}(L^g_2)^2\gamma_k^2\mathbb{E}\left[\left\|v^z_{k}-D_z(x_k,y_k,z_k)\right\|^2 \right]
     +\frac{4(1-p)}{b}(L^g_2)^2\gamma_k^2L_z^2\mathbb{E}\left[\left\|z_k-z^*\left(x_k\right)\right\|^2\right]
     \\&+\frac{4(1-p)}{b}(L^g_2)^2\gamma_k^2 L_z^2\mathbb{E}\left[\left\|y_k-y^*\left(x_k\right)\right\|^2\right].
\end{align*}
The proof of (3) is analogous, and hence we omit the details here.
\end{proof}

\begin{theorem}\textbf{(Restatement of Theorem \ref{thpagefini})}

    Fix an iteration $K>1$ and assume that Assumption \ref{assump UL} to \ref{assump LL} and \ref{assumptionijsmooth} hold.
    Choose minibatch size $b<(n+m)$ and the probability $p\in (0,1]$.
Then there exist positive constants $c$, $c_{\beta}$, and $c_{\gamma}$, such that if
\begin{gather*}
    \alpha_k\leq  \frac{c}{1+\sqrt{\frac{1-p}{pb}}},\quad\beta_k\leq c_{\beta}\alpha_k,\quad  \gamma_k\leq c_{\gamma}\alpha_k,\quad
    \gamma_k\leq c_{\gamma\beta}\beta_k,
\end{gather*}
the iterates in SPABA satisfy
$$ \frac{1}{K}\sum_{k=0}^{K-1}  \mathbb{E}\left[\left\|\nabla H\left(x_k\right)\right\|^2\right]
=\mathcal{O}\left(\frac{1+\sqrt{\frac{1-p}{pb}}}{K}\right).$$
\end{theorem}
\begin{proof}
We consider the Lyapunov function
\begin{align*}
   L_k=&H_k+\mathbb{E}\left[\left\|y_k-y^*\left(x_k\right)\right\|^2\right]+\mathbb{E}\left[\left\|z_k-z^*\left(x_k\right)\right\|^2\right]\\
   &+\frac{\alpha_k}{p}\left(\mathbb{E}\left[\left\|v_{k}^{x}-D_x(x_k,y_k,z_k)\right\|^2\right]+\mathbb{E}\left[\left\|v_{k}^{y}-D_y(x_k,y_k,z_k)\right\|^2\right]+\mathbb{E}\left[\left\|v_{k}^{z}-D_z(x_k,y_k,z_k)\right\|^2\right]\right)
\end{align*}

\begin{align*}
   &L_{k+1}-L_k
   \\
   \leq &
   -\frac{\alpha_k}{2} \mathbb{E}\left[\left\|\nabla H\left(x_k\right)\right\|^2\right]
   \\
   &+\left(\alpha_k-\alpha_k\right)
   \mathbb{E}\left[\left\|D_x(x_k,y_k,z_k)-v_k^x\right\|^2\right]
   \\
    &+\left\{\frac{L^H\alpha_k^2}{2}-\frac{\alpha_k}{2}+ \frac{2(\mu+L_1^g)L_{y^*}^2\alpha_k^2}{\mu L_1^g\beta_k}+ \frac{4L_{z^*}^2\alpha_k^2}{\mu \gamma_k}+\frac{(1-p)}{bp}(L_1^g)^2\alpha_k^3\right.
    \\&\quad \left.+\left(4(L^f)^2+8R^2(L_2^g)^2\right)\frac{(1-p)}{pb}\alpha_k^3
    \right\}
\mathbb{E}\left[\left\|v_k^x\right\|^2\right]\\
     &+\left\{3{\alpha_k} \left(\left(L^f\right)^2+\left(L^g_2 R\right)^2\right)-\frac{\mu L_1^g\beta_k}{2(\mu+L_1^g)}+\frac{2(1-p)}{bp}(L_1^g)^4\alpha_k\beta_k^2
     +\frac{18((L_2^gR)^2+(L^f)^2)}{\mu}\gamma_k
     \right.
     \\&\left.\quad +\left(4(L^f)^2+8R^2(L_2^g)^2\right)\frac{(1-p)}{bp}\alpha_k\beta_k^2\left(L^g_1\right)^2+\frac{8(1-p)}{bp}(L^g_2)^2\alpha_k\gamma_k^2 L_z^2\right\}
\mathbb{E}\left[\left\|y_k-y^*\left(x_k\right)\right\|^2\right] \\
      &+\left\{3{\alpha_k} \left(L^g_1\right)^2-\frac{\mu \gamma_k}{4}+\frac{8(1-p)}{bp}(L^g_2)^2\alpha_k\gamma_k^2 L_z^2\right\}
\mathbb{E}\left[\left\|z_k-z^*\left(x_k\right)\right\|^2\right]\\
        &+\left\{6\frac{\mu+L_1^g}{\mu L_1^g}\beta_k-\alpha_k +\frac{2(1-p)}{bp}(L_1^g)^2\alpha_k\beta_k^2
        +\left(4(L^f)^2+8R^2(L_2^g)^2\right)\frac{(1-p)}{bp}\alpha_k\beta_k^2\right\}
\mathbb{E}\left[\| v_k^y-D_y(x_k,y_k,z_k)\|^2\right]\\
        &+\left\{\frac{12}{\mu }\gamma_k+\frac{4(1-p)}{bp}(L^g_2)^2\alpha_k\gamma_k^2-\alpha_k\right\}
\mathbb{E}\left[\| v_k^z-D_z(x_k,y_k,z_k)\|^2\right]
   \end{align*}

We choose the step sizes to be
$$\alpha_k=\min\left\{\frac{1}{4L^H},\frac{c_{\alpha}}{\sqrt{\frac{1-p}{bp}}}\right\},\quad\beta_k={c_{\beta}}\alpha_k,\quad\gamma_k=c_{\gamma}\alpha_k.$$
Furthermore, by analyzing the coefficients of each term in the aforementioned inequalities, we can determine the range of values for $c_{\alpha}$, $c_{\beta}$ and $c_{\gamma}$.

\textbf{Analysis of the Coefficient for} $\mathbb{E}\left[\left\|v_k^x\right\|^2\right]$

By assuming
$$\alpha_k\leq\min\left\{\frac{1}{4L^H},\frac{\mu L_1^g}{16(\mu+ L_1^g)L_{y^*}}\beta_k,\frac{\mu}{32L_{z^*}^2}\gamma_k\right\},\quad c_{\alpha}\leq\frac{1}{4\Delta_1},\quad \Delta_1=\sqrt{(L_1^g)^2+ 4(L^f)^2+8R^2(L_2^g)^2},$$
we can deduce
\begin{align*}
&\frac{L^H\alpha_k^2}{2}-\frac{\alpha_k}{2}+ \frac{2(\mu+L_1^g)L_{y^*}^2\alpha_k^2}{\mu L_1^g\beta_k}
+\frac{4L_{z^*}^2\alpha_k^2}{\mu \gamma_k}
+\frac{(1-p)}{bp}(L_1^g)^2\alpha_k^3+\left(4(L^f)^2+8R^2(L_2^g)^2\right)\frac{(1-p)}{pb}\alpha_k^3
\\
=&\frac{L^H\alpha_k^2}{2}-\frac{\alpha_k}{2}+ \frac{2(\mu+L_1^g)L_{y^*}^2\alpha_k^2}{\mu L_1^g\beta_k}
+\frac{4L_{z^*}^2\alpha_k^2}{\mu \gamma_k}
+
\frac{(1-p)}{pb}\alpha_k^3\left((L_1^g)^2+ 4(L^f)^2+8R^2(L_2^g)^2\right)
\\
\triangleq&
\frac{L^H\alpha_k^2}{2}-\frac{\alpha_k}{2}+ \frac{2(\mu+L_1^g)L_{y^*}^2\alpha_k^2}{\mu L_1^g\beta_k}
+\frac{4L_{z^*}^2\alpha_k^2}{\mu \gamma_k}
+
\frac{(1-p)}{pb}\alpha_k^3\Delta_1^2
   \\\leq&  0.
\end{align*}
\textbf{Analysis of the Coefficient for} $\mathbb{E}\left[\left\|y_k-y^*\left(x_k\right)\right\|^2\right] $

By assuming
$$\alpha_k\leq \frac{\mu L_1^g\beta_k}{18(\mu+L_1^g)\left(\left(L^f\right)^2+\left(L^g_2 R\right)^2\right)},\quad c_{\alpha}^2\leq\frac{\mu L_1^g c_{\beta}}{12(\mu+ L_1^g)\Delta_2},
\quad
c_{\gamma\beta}\leq \frac{\mu^2 L_1^g}{108(\mu+L_1^g)\left(\left(L^f\right)^2+\left(L^g_2 R\right)^2\right)}$$
we can deduce
\begin{align*}
&-\frac{\mu L_1^g\beta_k}{2(\mu+L_1^g)}+3{\alpha_k} \left(\left(L^f\right)^2+\left(L^g_2 R\right)^2\right)
+\frac{18((L_2^gR)^2+(L^f)^2)}{\mu}\gamma_k
\\
&+\frac{(1-p)}{bp}\alpha_k^3
\left(2(L_1^g)^2 c_{\beta_k}^4
+\left(4(L^f)^2+8R^2(L_2^g)^2\right)\left(L^g_1\right)^2 c_{\beta_k}^2+8(L_2^g)^2L_z^2\right)
\\ \triangleq &-\frac{\mu L_1^g\beta_k}{2(\mu+L_1^g)}+3{\alpha_k} \left(\left(L^f\right)^2+\left(L^g_2 R\right)^2\right)+\frac{(1-p)}{bp}\alpha_k^3\Delta_2
+\frac{18((L_2^gR)^2+(L^f)^2)}{\mu}\gamma_k
\\&\leq 0.
\end{align*}
\textbf{Analysis of the Coefficient for} $\mathbb{E}\left[\left\|z_k-z^*\left(x_k\right)\right\|^2\right]$

By assuming
$$\alpha_k\leq \frac{\mu \gamma_k}{24 \left(L^g_1\right)^2},\quad c_{\alpha}^2c_{\gamma}\leq\frac{\mu 
}{64L_z^2(L_2^g)^2,}$$
we can deduce
\begin{align*}
3{\alpha_k} \left(L^g_1\right)^2-\frac{\mu \gamma_k}{4}+\frac{8(1-p)}{bp}(L^g_2)^2\alpha_k\gamma_k^2 L_z^2\leq 0.
\end{align*}

\textbf{Analysis of the Coefficient for} $\mathbb{E}\left[\| v_k^y-D_y(x_k,y_k,z_k)\|^2\right]$

By assuming
$$\beta_k\leq \frac{\mu L_1^g}{12(\mu+L_1^g)}\alpha_k,\quad c_{\alpha}^2\leq\frac{1}{4\Delta_3},$$
we can deduce
\begin{align*}
&6\frac{\mu+L_1^g}{\mu L-_1^g}\beta_k-\alpha_k +\frac{2(1-p)}{bp}(L_1^g)^2\alpha_k\beta_k^2
        +\left(4(L^f)^2+8R^2(L_2^g)^2\right)\frac{(1-p)}{bp}\alpha_k\beta_k^2
        \\&=6\frac{\mu+L_1^g}{\mu L_1^g}\beta_k-\alpha_k
        +\frac{(1-p)}{bp}\alpha_k^3\Delta_3\leq 0,
\end{align*}
where $$\Delta_3=\left(2(L_1^g)^2+4(L^f)^2+8R^2(L_2^g)^2\right)c_{\beta}^2.$$

\textbf{Analysis of the Coefficient for} $\mathbb{E}\left[\| v_k^z-D_z(x_k,y_k,z_k)\|^2\right]$

By assuming
$$\gamma_k\leq \frac{\mu }{24}\alpha_k,\quad c_{\alpha}^2c_{\gamma}\leq\frac{1}{12(L^g_2)^2},$$
we can deduce
\begin{align*}
\frac{12}{\mu }\gamma_k
+\frac{4(1-p)}{bp}(L^g_2)^2\alpha_k\gamma_k^2-\alpha_k\leq 0.
\end{align*}

Thus, we have obtained the recursive inequality for this theorem

   Summing, taking the average, and rearranging, we obtain
   \begin{align*}
   \frac{1}{K}\sum_{k=0}^{K-1} \mathbb{E}\left[\left\|\nabla H\left(x_k\right)\right\|^2\right]\leq
  \frac{L_0}{K\alpha_k}.\end{align*}

From the above analysis,
the step size $\alpha_k$ should satisfy
$$\alpha_k=\min\left\{\frac{1}{c_{\alpha}'},\frac{c_{\alpha}}{\sqrt{\frac{1-p}{bp}}}\right\}$$
then we have $$\alpha_k\leq  \frac{c}{1+\sqrt{\frac{1-p}{pb}}},$$
for some constants $c$.
Therefore, we ultimately arrive at the conclusion that
$$ \frac{1}{K}\sum_{k=0}^{K-1}  \mathbb{E}\left[\left\|\nabla H\left(x_k\right)\right\|^2\right]
=\mathcal{O}\left(\frac{1+\sqrt{\frac{1-p}{pb}}}{K}\right)$$
\end{proof}

\begin{corollary}
    Suppose that Assumption \ref{assump UL} to Assumption \ref{assumptionijsmooth} hold. If we take $p=b/(n+m+b)$, and $b\leq\sqrt{n+m}$, then the sample complexity is $\mathcal{O}((n+m)^{1/2}\epsilon^{-1})$.
\end{corollary}
\begin{proof}
  In each iteration, it uses $p(n+m)+(1-p)b$ samples on expectation.
  Let $p=\frac{b}{n+m+b}$ and $b\leq({n+m})^{1/2}$.
  Thus, the total sample complexity is
  \begin{align*}
    K(p(n+m)+(1-p)b)= \mathcal{O}\left(\left(1+\frac{\sqrt{n+m}}{b}\right)\frac{2(n+m)b}{n+m+b}\epsilon^{-1}\right)
    =\mathcal{O}\left(({n+m})^{1/2}\epsilon^{-1}\right).
  \end{align*}
\end{proof}

\section{Proof of Theorem \ref{thpage1.5}}\label{detial:proofpage1.5}
Under the expected form setting, in the algorithm, we set
$n=m=\tau'$, which represents the mini-batch size.
\begin{lemma}\label{pagegeneralmonment}
    Under the Assumption \ref{assump UL}, \ref{assump LL}, \ref{assumporacle} and \ref{assumptionijsmooth}, We have the following inequalities established:
\begin{align*}
   (1)\quad  \mathbb{E}\left[\left\|v_{k+1}^y-D_y\left(x_{k+1}, y_{k+1}, z_{k+1}\right)\right\|^2\right]
    \leq&(1-p)\mathbb{E}\left[\left\|v_{k}^{y}-D_y(x_k,y_k,z_k)\right\|^2\right]+\frac{(1-p)}{b}(L_1^g)^2\alpha_k^2\mathbb{E}\left[\left\|v_k^x\right\|^2\right]
         \\&+\frac{2(1-p)}{b}(L_1^g)^2\beta_k^2\mathbb{E}\left[\left\|v_k^y-D_y(x_k,y_k,z_k)\right\|^2\right]
         \\&+\frac{2(1-p)}{b}(L_1^g)^4\beta_k^2 \mathbb{E}\left[\| y_k-y^*\left(x_k\right) \|^2\right]\\&+\frac{p\sigma_{g,1}^2}{\tau'}.
\end{align*}
\begin{align*}
     (2)\quad \mathbb{E}\left[\left\|v_{k+1}^x-D_x\left(x_{k+1}, y_{k+1}, z_{k+1}\right)\right\|^2\right]
     &{\leq}(1-p)\mathbb{E}\left[\left\|v_{k}^{x}-D_x(x_k,y_k,z_k)\right\|^2\right]
    \\& +\left(2(L^f)^2+4R^2(L_2^g)^2\right)\frac{(1-p)}{b}\alpha_k^2\mathbb{E}\left[\left\|v_k^x\right\|^2\right]
    \\&+\left(2(L^f)^2+4R^2(L_2^g)^2\right)\frac{(1-p)}{b}\beta_k^2\mathbb{E}\left[\left\|v_k^y-D_y(x_k,y_k,z_k)\right\|^2\right]
    \\&+\left(2(L^f)^2+4R^2(L_2^g)^2\right)\frac{(1-p)}{b}\beta_k^2\left(L^g_1\right)^2 \mathbb{E}\left[\| y_k-y^*\left(x_k\right) \|^2\right]
     \\&
     +\frac{4(1-p)}{b}(L^g_2)^2\gamma_k^2\mathbb{E}\left[\left\|v^z_{k}-D_z(x_k,y_k,z_k)\right\|^2 \right]
      \\&+\frac{4(1-p)}{b}(L^g_2)^2\gamma_k^2L_z^2\mathbb{E}\left[\left\|z_k-z^*\left(x_k\right)\right\|^2\right]
     \\&+\frac{4(1-p)}{b}(L^g_2)^2\gamma_k^2 L_z^2\mathbb{E}\left[\left\|y_k-y^*\left(x_k\right)\right\|^2\right]
     \\&
     +\frac{2p\sigma_f^2}{\tau'}+\frac{2p\sigma_{g,2}^2}{\tau'}.
\end{align*}
\begin{align*}
    (3)\quad \mathbb{E}\left[\left\|v_{k+1}^z-D_z\left(x_{k+1}, y_{k+1}, z_{k+1}\right)\right\|^2\right]
     \leq
     &(1-p+\frac{4(1-p)}{b}(L^g_2)^2\gamma_k^2)\mathbb{E}\left[\left\|v_{k}^{z}-D_z(x_k,y_k,z_k)\right\|^2\right]
     \\&
     +\left(2(L^f)^2+4R^2(L_2^g)^2\right)\frac{(1-p)}{b}\alpha_k^2\mathbb{E}\left[\left\|v_k^x\right\|^2\right]
    \\&+\left(2(L^f)^2+4R^2(L_2^g)^2\right)\frac{(1-p)}{b}\beta_k^2\mathbb{E}\left[\left\|v_k^y-D_y(x_k,y_k,z_k)\right\|^2\right]
    \\&+\left(2(L^f)^2+4R^2(L_2^g)^2\right)\frac{(1-p)}{b}\beta_k^2\left(L^g_1\right)^2 \mathbb{E}\left[\| y_k-y^*\left(x_k\right) \|^2\right]
     \\&
     +\frac{4(1-p)}{b}(L^g_2)^2\gamma_k^2L_z^2\mathbb{E}\left[\left\|z_k-z^*\left(x_k\right)\right\|^2\right]
     \\&+\frac{4(1-p)}{b}(L^g_2)^2\gamma_k^2 L_z^2\mathbb{E}\left[\left\|y_k-y^*\left(x_k\right)\right\|^2\right]
     \\& +\frac{2p\sigma_f^2}{\tau'}+\frac{2p\sigma_{g,2}^2}{\tau'}.
\end{align*}
\end{lemma}
\begin{proof}
    \textbf{Proof of (1).}

By the definition of $v_{k+1}^y$, we have
 \begin{align*}
         &\mathbb{E}\left[\left\|v_{k+1}^y-D_y\left(x_{k+1}, y_{k+1}, z_{k+1}\right)\right\|^2\right]
         \\=&p\mathbb{E}\left[\left\|\frac{1}{\tau'} \sum_{j \in [\tau']} \nabla_2 G\left(x_{k+1},y_{k+1};\zeta_j\right)-D_y\left(x_{k+1}, y_{k+1}, z_{k+1}\right)\right\|^2\right]
         \\&+(1-p)\mathbb{E}\left[\left\|v_{k}^{y}+\frac{1}{b} \sum_{\zeta_j \in J}\left[\nabla_2 G\left(x_{k+1},y_{k+1};\zeta_j\right)-\nabla_2 G\left(x_{k},y_{k};\zeta_j\right)\right]-D_y\left(x_{k+1}, y_{k+1}, z_{k+1}\right)\right\|^2\right]
         \end{align*}
For the first term,
based on Assumption \ref{assumporacle},
we have
\begin{align*}
    p\mathbb{E}\left[\left\|\frac{1}{\tau'} \sum_{j \in [\tau']} \nabla_2 G\left(x_{k+1},y_{k+1};\zeta_j\right)-D_y\left(x_{k+1}, y_{k+1}, z_{k+1}\right)\right\|^2\right]
    \leq\frac{p\sigma_{g,1}^2}{\tau'} .
\end{align*}
For the second term,
analogous to the proof of Lemma \ref{pagefinitemoment}, we have
\begin{align*}
    &(1-p)\mathbb{E}\left[\left\|v_{k}^{y}+\frac{1}{b} \sum_{\zeta_j \in J}\left[\nabla_2 G\left(x_{k+1},y_{k+1};\zeta_j\right)-\nabla_2 G\left(x_{k},y_{k};\zeta_j\right)\right]-D_y\left(x_{k+1}, y_{k+1}, z_{k+1}\right)\right\|^2\right]
    \\\leq&
    (1-p)\mathbb{E}\left[\left\|v_{k}^{y}-D_y(x_k,y_k,z_k)\right\|^2\right]+\frac{(1-p)}{b}(L_1^g)^2\alpha_k^2\mathbb{E}\left[\left\|v_k^x\right\|^2\right]
         \\&+\frac{2(1-p)}{b}(L_1^g)^2\beta_k^2\mathbb{E}\left[\left\|v_k^y-D_y(x_k,y_k,z_k)\right\|^2\right]
         +\frac{2(1-p)}{b}(L_1^g)^4\beta_k^2 \mathbb{E}\left[\| y_k-y^*\left(x_k\right) \|^2\right].
\end{align*}
In summary, (1) is proved.

\textbf{Proof of (2) and (3).}
Similarly, from the definition of $v_{k+1}^x$, we have
\begin{align*}
    &\mathbb{E}\left[\left\|v_{k+1}^x-D_x\left(x_{k+1}, y_{k+1}, z_{k+1}\right)\right\|^2\right]\\
    =&
    p\mathbb{E}\left[\left\|\frac{1}{\tau'} \sum_{i \in [\tau']}\nabla_1 F\left(x_{k+1},y_{k+1};\xi_i\right)
    -\frac{1}{\tau'} \sum_{j \in [\tau']} \nabla_{12}^2 G\left(x_{k+1},y_{k+1};\zeta_j\right)-D_y\left(x_{k+1}, y_{k+1}, z_{k+1}\right)\right\|^2\right]
    \\
    +&(1-p)\mathbb{E}\left[\left\|v_{k}^{x}+\frac{1}{b} \sum_{i \in I}\left( \nabla_{1} F_i\left(x_{k+1},y_{k+1}\right)-\nabla_{1} F_i\left(x_{k},y_{k}\right)\right)\right.\right.\\&\left. \left.- \frac{1}{b} \sum_{j \in J}\left( \nabla^2_{12} G_j\left(x_{k+1},y_{k+1}\right)z_{k+1-} \nabla^2_{12} G_j\left(x_{k},y_{k}\right)z_{k}\right)-D_x\left(x_{k+1}, y_{k+1}, z_{k+1}\right)\right\|^2\right].
    \end{align*}
For the first term,
based on Assumption \ref{assumporacle},
we have
\begin{align*}
  & p\mathbb{E}\left[\left\|\frac{1}{\tau'} \sum_{i \in [\tau']}\nabla_1 F\left(x_{k+1},y_{k+1};\xi_i\right)
    -\frac{1}{\tau'} \sum_{j \in [\tau']} \nabla_{12}^2 G\left(x_{k+1},y_{k+1};\zeta_j\right)-D_y\left(x_{k+1}, y_{k+1}, z_{k+1}\right)\right\|^2\right]
    \\
    &\leq \frac{2p\sigma_f^2}{\tau'}+\frac{2p\sigma_{g,2}^2}{\tau'}.
\end{align*}
For the second term,
analogous to the proof of Lemma \ref{pagefinitemoment}, we have
\begin{align*}
 &(1-p)\mathbb{E}\left[\left\|v_{k}^{x}+\frac{1}{b} \sum_{i \in I}\left( \nabla_{1} F_i\left(x_{k+1},y_{k+1}\right)-\nabla_{1} F_i\left(x_{k},y_{k}\right)\right)\right.\right.\\&\left. \left.- \frac{1}{b} \sum_{j \in J}\left( \nabla^2_{12} G_j\left(x_{k+1},y_{k+1}\right)z_{k+1-} \nabla^2_{12} G_j\left(x_{k},y_{k}\right)z_{k}\right)-D_x\left(x_{k+1}, y_{k+1}, z_{k+1}\right)\right\|^2\right]
    \\ {\leq}&(1-p)\mathbb{E}\left[\left\|v_{k}^{x}-D_x(x_k,y_k,z_k)\right\|^2\right]
     +\left(2(L^f)^2+4R^2(L_2^g)^2\right)\frac{(1-p)}{b}\alpha_k^2\mathbb{E}\left[\left\|v_k^x\right\|^2\right]
    \\&+\left(2(L^f)^2+4R^2(L_2^g)^2\right)\frac{(1-p)}{b}\beta_k^2\mathbb{E}\left[\left\|v_k^y-D_y(x_k,y_k,z_k)\right\|^2\right]
    \\&+\left(2(L^f)^2+4R^2(L_2^g)^2\right)\frac{(1-p)}{b}\beta_k^2\left(L^g_1\right)^2 \mathbb{E}\left[\| y_k-y^*\left(x_k\right) \|^2\right]
     \\&
     +\frac{4(1-p)}{b}(L^g_2)^2\gamma_k^2\mathbb{E}\left[\left\|v^z_{k}-D_z(x_k,y_k,z_k)\right\|^2 \right]
      \\&+\frac{4(1-p)}{b}(L^g_2)^2\gamma_k^2L_z^2\mathbb{E}\left[\left\|z_k-z^*\left(x_k\right)\right\|^2\right]
     +\frac{4(1-p)}{b}(L^g_2)^2\gamma_k^2 L_z^2\mathbb{E}\left[\left\|y_k-y^*\left(x_k\right)\right\|^2\right].
\end{align*}
Therefore, (2) is proven.
The proof of (3) is analogous, and hence we omit the details here.
\end{proof}

\begin{theorem}\textbf{(Restatement of Theorem \ref{thpage1.5})}
    Fix an iteration $K>1$ and assume that Assumption \ref{assump UL} to Assumption \ref{assumptionijsmooth} hold.
Choose minibatch size $\tau'$ and $b<\tau'$, the probability $p\in (0,1]$.
Then there exist positive constants $c$, $c_{\beta}$, and $c_{\gamma}$, such that if
\begin{gather*}
    \alpha_k\leq  \frac{c}{1+\sqrt{\frac{1-p}{pb}}},\quad\beta_k=c_{\beta}\alpha_k,\quad  \gamma_k=c_{\gamma}\alpha_k,
\end{gather*}
the iterates in SPABA satisfy
$$ \frac{1}{K}\sum_{k=0}^{K-1}  \mathbb{E}\left[\left\|\nabla H\left(x_k\right)\right\|^2\right]
=\mathcal{O}\left(\frac{1+\sqrt{\frac{1-p}{pb}}}{K}+\frac{1}{Kp\tau'}+\frac{\sigma}{\tau'}\right).$$
\end{theorem}
\begin{proof}
We consider the Lyapunov function
\begin{align*}
   L_k=&H_k+\mathbb{E}\left[\left\|y_k-y^*\left(x_k\right)\right\|^2\right]+\mathbb{E}\left[\left\|z_k-z^*\left(x_k\right)\right\|^2\right]\\
   &+\frac{\alpha_k}{p}\left(\mathbb{E}\left[\left\|v_{k}^{x}-D_x(x_k,y_k,z_k)\right\|^2\right]+\mathbb{E}\left[\left\|v_{k}^{y}-D_y(x_k,y_k,z_k)\right\|^2\right]+\mathbb{E}\left[\left\|v_{k}^{z}-D_z(x_k,y_k,z_k)\right\|^2\right]\right)
\end{align*}

\begin{align*}
   &L_{k+1}-L_k
   \\
   \leq &
   -\frac{\alpha_k}{2} \mathbb{E}\left[\left\|\nabla H\left(x_k\right)\right\|^2\right]
   \\
   &+\left(\alpha_k-\alpha_k\right)
   \mathbb{E}\left[\left\|D_x(x_k,y_k,z_k)-v_k^x\right\|^2\right]
   \\
    &+\left\{
    \frac{L^H\alpha_k^2}{2}-\frac{\alpha_k}{2}+ \frac{2(\mu+L_1^g)L_{y^*}^2\alpha_k^2}{\mu L_1^g\beta_k}
    + \frac{4L_{z^*}^2\alpha_k^2}{\mu \gamma_k}
    +\frac{(1-p)}{bp}(L_1^g)^2\alpha_k^3
   +\left(4(L^f)^2+8R^2(L_2^g)^2\right)\frac{(1-p)}{pb}\alpha_k^3
    \right\}
\mathbb{E}\left[\left\|v_k^x\right\|^2\right]\\
     &+\left\{3{\alpha_k} \left(\left(L^f\right)^2+\left(L^g_2 R\right)^2\right)-\frac{\mu L_1^g\beta_k}{2(\mu+L_1^g)}+\frac{2(1-p)}{bp}(L_1^g)^4\alpha_k\beta_k^2
     +\frac{18((L_2^gR)^2+(L^f)^2)}{\mu}\gamma_k
      \right.
     \\&\left.\quad +\left(4(L^f)^2+8R^2(L_2^g)^2\right)\frac{(1-p)}{bp}\alpha_k\beta_k^2\left(L^g_1\right)^2+\frac{8(1-p)}{bp}(L^g_2)^2\alpha_k\gamma_k^2 L_z^2\right\}
\mathbb{E}\left[\left\|y_k-y^*\left(x_k\right)\right\|^2\right] \\
      &+\left\{3{\alpha_k} \left(L^g_1\right)^2
      -\frac{\mu \gamma_k}{4}+\frac{8(1-p)}{bp}(L^g_2)^2\alpha_k\gamma_k^2 L_z^2\right\}
\mathbb{E}\left[\left\|z_k-z^*\left(x_k\right)\right\|^2\right]\\
        &+\left\{6\frac{\mu+L_1^g}{\mu L_1^g}\beta_k^2-\alpha_k +\frac{2(1-p)}{bp}(L_1^g)^2\alpha_k\beta_k^2
        +\left(4(L^f)^2+8R^2(L_2^g)^2\right)\frac{(1-p)}{bp}\alpha_k\beta_k^2\right\}
\mathbb{E}\left[\| v_k^y-D_y(x_k,y_k,z_k)\|^2\right]\\
        &+\left\{\frac{12}{\mu }\gamma_k+\frac{4(1-p)}{bp}(L^g_2)^2\alpha_k\gamma_k^2-\alpha_k\right\}
\mathbb{E}\left[\| v_k^z-D_z(x_k,y_k,z_k)\|^2\right]
\\
&+\frac{\sigma_{g,1}^2\alpha_k}{\tau'}+\frac{4\sigma_f^2\alpha_k}{\tau'}+\frac{4\sigma_{g,2}\alpha_k}{\tau'}.
   \end{align*}

Following the proof process of Theorem \ref{thpagefini}, we have obtained the recursive inequality for this theorem
\begin{align*}
   \frac{\alpha_k}{2} \mathbb{E}\left[\left\|\nabla H\left(x_k\right)\right\|^2\right]\leq L_k -L_{k+1}+\frac{\sigma_{g,1}^2\alpha_k}{\tau'}+\frac{4\sigma_f^2\alpha_k}{\tau'}+\frac{4\sigma_{g,2}^2\alpha_k}{\tau'}\leq L_k -L_{k+1}+\frac{\sigma}{\tau'}\alpha_k,
   \end{align*}
   where $$\sigma=\sigma_{g,1}^2+4\sigma_f^2+4\sigma_{g,2}^2,\quad\alpha_k\leq  \frac{c}{1+\sqrt{\frac{1-p}{pb}}}.$$

   Summing, taking the average, and rearranging, we obtain
   \begin{align*}
   \frac{1}{K}\sum_{k=0}^{K-1} \mathbb{E}\left[\left\|\nabla H\left(x_k\right)\right\|^2\right]&\leq
  \frac{2L_0}{K\alpha_k}+\frac{2\sigma}{\tau'}= \frac{2L_0'}{K\alpha_k}+\frac{2}{Kp\tau'}+\frac{2\sigma}{\tau'}
  ,\end{align*}
where the last equation is based on the fact that
$$L_0=H(x_0)+\mathbb{E}[\|y_0-y^*(y_0)\|^2]+\mathbb{E}[\|z_0-z^*(y_0)\|^2]+\frac{\alpha_k{\sigma}}{p\tau'}\triangleq L_0'+\frac{\alpha_k{\sigma}}{p\tau'}.$$

Therefore, we ultimately arrive at the conclusion that
$$ \frac{1}{K}\sum_{k=0}^{K-1}  \mathbb{E}\left[\left\|\nabla H\left(x_k\right)\right\|^2\right]
=\mathcal{O}\left(\frac{1+\sqrt{\frac{1-p}{pb}}}{K}+\frac{1}{Kp\tau'}+\frac{\sigma}{\tau'}\right).$$
\end{proof}

\begin{corollary}
    Suppose that Assumption \ref{assump UL} to Assumption \ref{assumptionijsmooth} hold. If we take $p=b/(n+m+b)$,  $\tau'=\mathcal{O}(\epsilon^{-1})$ and $b\leq\sqrt{\tau'}$, then the sample complexity is $\mathcal{O}(\epsilon^{-1.5})$.
\end{corollary}
\begin{proof}
  In each iteration, it uses $p(n+m)+(1-p)b$ samples on expectation.
  Let $p=\frac{b}{n+m+b}$, $\tau'=\mathcal{O}(\epsilon^{-1})$ and $b\leq\sqrt{\tau'}$.
  Thus, the total sample complexity is
  \begin{align*}
    K(p\tau'+(1-p)b)= \mathcal{O}\left(\epsilon^{-1}\left(1+\frac{\sqrt{\tau'}}{b}+\frac{\tau'+b}{\tau' b}\right)\frac{2\tau'b}{\tau'+b}\right)=\mathcal{O}\left(\sqrt{\tau'}\epsilon^{-1}\right)=\mathcal{O}\left(\epsilon^{-1.5}\right).
  \end{align*}
\end{proof}

\section{Proof of Theorem \ref{thstorm}}\label{detial:proofs trom}
\begin{lemma}\label{eg31}
    Suppose Assumptions \ref{assump UL}, \ref{assump LL} and \ref{assumporacle} hold, then we have
\begin{align*}
 \mathbb{E}\left[\left\|D_{k+1}^x-D_k^x\right\|^2\right]
 \leq&
\Tilde{\Delta}\alpha_k^2\mathbb{E}[\|v_k^x\|^2]
+\Tilde{\Delta}\beta_k^2\mathbb{E}[\|v_k^y\|^2]
+4\left(L_1^g\right)^2\gamma_k^2\mathbb{E}[\|v_k^z\|^2],\\
\mathbb{E}\left[\left\|D_{k+1}^z-D_k^z\right\|^2\right]
 \leq&
\Tilde{\Delta}\alpha_k^2\mathbb{E}[\|v_k^x\|^2]
+\Tilde{\Delta}\beta_k^2\mathbb{E}[\|v_k^y\|^2]
+4\left(L_1^g\right)^2\gamma_k^2\mathbb{E}[\|v_k^z\|^2],
\end{align*}
where $\Tilde{\Delta}=2\left(L^f\right)^2+4R^2\left(L_2^g\right)^2$.
\end{lemma}
\begin{proof}
Taking the expectation conditionally on $x_k$, $y_k$, $z_k$ yields
\begin{align*}
 E_k\left[\left\|D_{k+1}^x-D_k^x\right\|^2\right]
= & E_k\left[\left\|\nabla_1 F\left(x_{k+1}, y_{k+1};\xi\right)-\nabla_1 F\left(x_k, y_k;\xi\right)\right.
\right.
\\&\left.\left.-\nabla_{12}^2 G\left(x_{k+1}, y_{k+1};\zeta\right) z_{k+1}-\nabla_{12}^2 G\left(x_k, y_k;\zeta\right) z_k\right\|^2\right] \\
\leq & 2 E_k\left[\left\|\nabla_1  F\left(x_{k+1}, y_{k+1};\xi\right)-\nabla_1 F\left(x_k, y_k;\xi\right)\right\|^2\right]
\\&+4 E_k\left[\left\|\left(\nabla_{12}^2 G\left(x_{k+1}, y_{k+1};\zeta\right) -\nabla_{12}^2 G\left(x_k, y_k;\zeta\right) \right)z_{k+1}\right\|^2\right] \\
& +4 E_k\left[\left\|\nabla_{12}^2 G\left(x_k, y_k;\zeta\right)\left(z_{k+1}-z_k\right)\right\|^2\right],
\end{align*}
where the inequality is derived using the Cauchy-Schwarz inequality.
For the term $\left\|z_{k+1}-z_{k}\right\|^2$, based on the definition of $z_{k+1}$ and the contractivity of projection, we have
\begin{align*}
    \left\|z_{k+1}-z_{k}\right\|^2 = \left\|\text{Proj}_{\mathbb{B}(R)}\left(z_k-\gamma_k v^z_k\right)-z_{k}\right\|^2
     &= \left\|\text{Proj}_{\mathbb{B}(R)}\left(z_k-\gamma_k v^z_k\right)-\text{Proj}_{\mathbb{B}(R)}(z_{k})\right\|^2\\
     &\leq \left\|z_k-\gamma_k v^z_k-z_{k}\right\|^2
     =\gamma_k^2\|v_k^z\|^2.
\end{align*}
Thus, substituting into the above equation, we obtain
\begin{align*}
 E_k\left[\left\|D_{k+1}^x-D_k^x\right\|^2\right]
\leq &
2\left(L^f\right)^2\left(\alpha_k^2E_k[\|v_k^x\|^2]+\beta_k^2E_k[\|v_k^y\|^2]\right)\\
&+4R^2\left(L_2^g\right)^2\left(\alpha_k^2E_k[\|v_k^x\|^2]+\beta_k^2E_k[\|v_k^y\|^2]\right)
+4\left(L_1^g\right)^2\gamma_k^2E_k[\|v_k^z\|^2]
\\=&
\left(2\left(L^f\right)^2+4R^2\left(L_2^g\right)^2\right)\alpha_k^2E_k[\|v_k^x\|^2]
\\&+\left(2\left(L^f\right)^2+4R^2\left(L_2^g\right)^2\right)\beta_k^2E_k[\|v_k^y\|^2]
+4\left(L_1^g\right)^2\gamma_k^2E_k[\|v_k^z\|^2],
\end{align*}
where the inequality results from the boundedness generated by projecting $z_k$, as well as Assumptions \ref{assump UL}, \ref{assump LL} and \ref{assumporacle}.
Finally, by taking the total expectation, the lemma is proven.

Similarly, we can derive the inequality concerning $\mathbb{E}\left[\left\|D_{k+1}^z-D_k^z\right\|^2\right]$.
\end{proof}

\begin{lemma}\label{eg33}
Under the Assumption \ref{assump UL}, \ref{assump LL}, \ref{assumporacle} and \ref{assumptionijsmooth}, We have the following inequalities established:
\begin{align*}
   (1)\quad& \mathbb{E}\left[\left\|v_{k+1}^y-D_y\left(x_{k+1}, y_{k+1}, z_{k+1}\right)\right\|^2\right]\\
    \leq& \left(\left(1-\rho_k^y\right)^2+4\left(1-\rho_k^y\right)^2\left(L_1^g\right)^2\beta_k^2\right) \mathbb{E}\left[\left\|v_k^y-D_y\left(x_k, y_k, z_k\right)\right\|^2\right]
    \\&+2\left(1-\rho_k^y\right)^2\left(L_1^g\right)^2\alpha_k^2 \mathbb{E}\left[\left\|v_k^x\right\|^2\right]
    +4\left(1-\rho_k^y\right)^2\left(L_1^g\right)^4\beta_k^2\mathbb{E}\left[\left\|y_k-y^*(x_k)\right\|^2\right]
    +2\left(\rho_k^y\right)^2 \sigma_{g,1}^2,
\end{align*}
\begin{align*}
(2)\quad	&E_k\left[\left\|v_{k+1}^x-D_x\left(x_{k+1}, y_{k+1}, z_{k+1}\right)\right\|^2\right] \\
	\leq&
	\left(1-\rho_k^x\right)^2 E_k\left[\left\|v_k^x-D_x\left(x_k, y_k, z_k\right)\right\|^2\right]  
	+2\left(1-\rho_k^x\right)^2 \Tilde{\Delta}\alpha_k^2E_k[\|v_k^x \|^2]
	\\&
	+4\left(1-\rho_k^x\right)^2\Tilde{\Delta}\beta_k^2 E_k[\|v_k^y-D_y(x_k,y_k,z_k) \|^2]
	\\&
	+8\left(1-\rho_k^x\right)^2\gamma_k^2(L_1^g)^2E_k[\|v_k^z -D_z(x_k,y_k,z_k)\|^2] 
	\\&
	+\left(8\left(1-\rho_k^x\right)^2\gamma_k^2(L_1^g)^2L_z^2+4\left(1-\rho_k^x\right)^2(L_1^g)^2\Tilde{\Delta}\beta_k^2\right)
	 E_k[\|y_k-y^*(x_k)\|^2]
	\\&
	+8\left(1-\rho_k^x\right)^2\gamma_k^2(L_1^g)^2L_z^2 E_k[\|z_k-z^*(x_k)\|^2]\\&
	+2\left(\rho_k^x\right)^2 \left(\sigma_{f}+R^2\sigma_{g,2}\right)  ,
\end{align*}
\begin{align*}
	(3)\quad	&E_k\left[\left\|v_{k+1}^z-D_z\left(x_{k+1}, y_{k+1}, z_{k+1}\right)\right\|^2\right] \\
	\leq&
	\left(\left(1-\rho_k^z\right)^2 +8\left(1-\rho_k^z\right)^2\gamma_k^2(L_1^g)^2 \right)
	E_k\left[\left\|v_k^z-D_z\left(x_k, y_k, z_k\right)\right\|^2\right] 
	 \\&
	+2\left(1-\rho_k^z\right)^2 \Tilde{\Delta}\alpha_k^2E_k[\|v_k^x \|^2]
	\\&
	+4\left(1-\rho_k^z\right)^2\Tilde{\Delta}\beta_k^2 E_k[\|v_k^y-D_y(x_k,y_k,z_k) \|^2]
	\\&
	+\left(8\left(1-\rho_k^z\right)^2\gamma_k^2(L_1^g)^2L_z^2
	+4\left(1-\rho_k^z\right)^2(L_1^g)^2\Tilde{\Delta}\beta_k^2\right)
	 E_k[\|y_k-y^*(x_k)\|^2]
	\\&
	+8\left(1-\rho_k^z\right)^2\gamma_k^2(L_1^g)^2L_z^2 E_k[\|z_k-z^*(x_k)\|^2]
	\\&+2\left(\rho_k^z\right)^2 \left(\sigma_{f}+R^2\sigma_{g,2}\right)  .
\end{align*}
\end{lemma}
\begin{proof}

\textbf{proof of (1)}

By the definition of $v_{k+1}^y$, we have
\begin{align*}
\left\|v_{k+1}^y-D_y\left(x_{k+1}, y_{k+1}, z_{k+1}\right)\right\|^2
= &\left\|D_{k+1}^y+\left(1-\rho_k^y\right)\left(v_k^y-D_k^y\right)-D_y\left(x_{k+1}, y_{k+1}, z_{k+1}\right)\right\|^2 \\
=& \left\|\left(1-\rho_k^y\right)\left(v_k^y-D_y\left(x_k, y_k, z_k\right)\right)+\rho_k^y\left(D_{k+1}^y-D_y\left(x_{k+1}, y_{k+1}, z_{k+1}\right)\right)\right. \\
&\left.+\left(1-\rho_k^y\right)\left(D_{k+1}^y-D_y\left(x_{k+1}, y_{k+1}, z_{k+1}\right)-D_k^y+D_y\left(x_k, y_k, z_k\right)\right) \right\|^2 ,
\end{align*}

Taking the expectation conditionally on $x_k$, $y_k$, $z_k$, and utilizing that
$D_{k+1}^y$ and $D_{k}^y$ are unbiased estimates of $D_y(x_{k+1},y_{k+1},z_{k+1})$ and $D_y(x_{k},y_{k},z_{k})$ respectively, yields
\begin{align*}
&E_k\left[\left\|v_{k+1}^y-D_y\left(x_{k+1}, y_{k+1}, z_{k+1}\right)\right\|^2\right]
\\\leq&\left(1-\rho_k^y\right)^2 E_k\left[\left\|v_k^y-D_y\left(x_k, y_k, z_k\right)\right\|^2\right]  +2\left(\rho_k^y\right)^2 E_k\left[\left\|D_{k+1}^y-D_y\left(x_{k+1}, y_{k+1}, z_{k+1}\right)\right\|^2\right]  \\
&+2\left(1-\rho_k^y\right)^2 E_k\left[\left\|D_{k+1}^y-D_y\left(x_{k+1}, y_{k+1}, z_{k+1}\right)-D_k^y+D_y\left(x_k, y_k, z_k\right)\right\|^2\right],
\end{align*}
For the second term, by setting  $|\mathcal{S}_2|=1$ in Lemma \ref{eg1}, we obtain
$$E_k\left[\left\|D_{k+1}^y-D_y\left(x_{k+1}, y_{k+1}, z_{k+1}\right)\right\|^2\right]\leq \sigma_{g,1}^2,$$
For the third term, from Assumption \ref{assumptionijsmooth}, we obtain
\begin{align*}
&E_k\left[\left\|D_{k+1}^y-D_y\left(x_{k+1}, y_{k+1}, z_{k+1}\right)-D_k^y+D_y\left(x_k, y_k, z_k\right)\right\|^2\right]\\
\leq&E_k\left[\left\|\nabla_2 G(x_{k+1},y_{k+1};\zeta)-\nabla_2 G(x_{k},y_{k};\zeta)\right\|^2\right]
\leq \left(L_1^g\right)^2\left(\alpha_k^2 E_k\left[\left\|v_k^x\right\|^2\right]+\beta_k^2E_k\left[\left\|v_k^y\right\|^2\right]\right),
\end{align*}
Taking the total expectation, ultimately we can derive
\begin{align*}
    \mathbb{E}\left[\left\|v_{k+1}^y-D_y\left(x_{k+1}, y_{k+1}, z_{k+1}\right)\right\|^2\right]
    \leq& \left(1-\rho_k^y\right)^2 \mathbb{E}\left[\left\|v_k^y-D_y\left(x_k, y_k, z_k\right)\right\|^2\right]
    +2\left(\rho_k^y\right)^2 \sigma_{g,1}^2\\
    &+2\left(1-\rho_k^y\right)^2\left(L_1^g\right)^2\alpha_k^2 \mathbb{E}\left[\left\|v_k^x\right\|^2\right]
    +2\left(1-\rho_k^y\right)^2\left(L_1^g\right)^2\beta_k^2\mathbb{E}\left[\left\|v_k^y\right\|^2\right]
    \\
    \leq& \left(\left(1-\rho_k^y\right)^2+4\left(1-\rho_k^y\right)^2\left(L_1^g\right)^2\beta_k^2\right) \mathbb{E}\left[\left\|v_k^y-D_y\left(x_k, y_k, z_k\right)\right\|^2\right]
    \\
    &+2\left(1-\rho_k^y\right)^2\left(L_1^g\right)^2\alpha_k^2 \mathbb{E}\left[\left\|v_k^x\right\|^2\right]
    \\
    &+4\left(1-\rho_k^y\right)^2\left(L_1^g\right)^4\beta_k^2\mathbb{E}\left[\left\|y_k-y*(x_k)\right\|^2\right]
    \\
    &+2\left(\rho_k^y\right)^2 \sigma_{g,1}^2.
\end{align*}
Thus, the (1) is proven.

\textbf{proof of (2) and (3)}

Based on the definition of $v_{k+1}^x$ and the fact that $D_{k+1}^x$ and $D_{k}^x$ are unbiased estimates of $D_x(x_{k+1},y_{k+1},z_{k+1})$ and $D_x(x_{k},y_{k},z_{k})$ respectively, we have
\begin{align*}
&E_k\left[\left\|v_{k+1}^x-D_x\left(x_{k+1}, y_{k+1}, z_{k+1}\right)\right\|^2\right] \\
= &E_k\left\|D_{k+1}^x+\left(1-\rho_k^x\right)\left(v_k^x-D_k^x\right)-D_x\left(x_{k+1}, y_{k+1}, z_{k+1}\right)\right\|^2 \\
=&E_k \left[\left\|\left(1-\rho_k^x\right)\left(v_k^x-D_x\left(x_k, y_k, z_k\right)\right)+\rho_k^x\left(D_{k+1}^x-D_x\left(x_{k+1}, y_{k+1}, z_{k+1}\right)\right)\right.\right. \\
&\left.\left.+\left(1-\rho_k^x\right)\left(D_{k+1}^x-D_x\left(x_{k+1}, y_{k+1}, z_{k+1}\right)-D_k^x+D_x\left(x_k, y_k, z_k\right)\right) \right\|^2 \right]\\
\leq&\left(1-\rho_k^x\right)^2 E_k\left[\left\|v_k^x-D_x\left(x_k, y_k, z_k\right)\right\|^2\right]  +2\left(\rho_k^x\right)^2 E_k\left[\left\|D_{k+1}^x-D_x\left(x_{k+1}, y_{k+1}, z_{k+1} \right)  \right\|^2\right]  \\
&+2\left(1-\rho_k^x\right)^2 E_k\left[\left\|D_{k+1}^x-D_k^x -D_x\left(x_{k+1}, y_{k+1}, z_{k+1}\right) + D_x\left(x_{k}, y_{k}, z_{k}\right)\right\|^2\right].
\end{align*}
For the term $E_k\left[\left\|D_{k+1}^x-D_x\left(x_{k+1}, y_{k+1}, z_{k+1} \right)  \right\|^2\right]$, we have
\begin{align*}
&E_k\left[\left\|D_{k+1}^x-D_x\left(x_{k+1}, y_{k+1}, z_{k+1} \right)  \right\|^2\right]
\\=&E_k\left[\left\|\nabla_1 F(x_{k+1},y_{k+1};\xi)-\nabla_1 F(x_{k+1},y_{k+1}) -\nabla_{12}^2G(x_{k+1},y_{k+1};\zeta)z_{k+1} +\nabla_{12}^2G(x_{k+1},y_{k+1})z_{k+1}  \right\|^2\right]
\\
=&
E_k\left[\left\|\nabla_1 F(x_{k+1},y_{k+1};\xi)-\nabla_1 F(x_{k+1},y_{k+1})\right\|^2\right] +E_k\left[\left\|-\nabla_{12}^2G(x_{k+1},y_{k+1};\zeta)z_{k+1} +\nabla_{12}^2G(x_{k+1},y_{k+1})z_{k+1}  \right\|^2\right]
\\\leq&
\sigma_{f}+R^2\sigma_{g,2}.
\end{align*}
For the term $E_k\left[\left\|D_{k+1}^x-D_k^x -D_x\left(x_{k+1}, y_{k+1}, z_{k+1}\right) + D_x\left(x_{k}, y_{k}, z_{k}\right)\right\|^2\right]$, we have
\begin{align*}
	&E_k\left[\left\|D_{k+1}^x-D_k^x -D_x\left(x_{k+1}, y_{k+1}, z_{k+1}\right) + D_x\left(x_{k}, y_{k}, z_{k}\right)\right\|^2\right]
	\\\leq&
	E_k\left[\left\|\nabla_1 F(x_{k+1},y_{k+1};\xi)-\nabla_1 F(x_k,y_k;\xi)\right\|^2\right]
	+E_k\left[\left\|\nabla_{12}^2 G(x_{k+1},y_{k+1};\zeta)z_{k+1}-\nabla_{12}^2 G(x_k,y_k;\zeta)z_{k} \right\|^2\right]
\\\leq&
E_k\left[\left\|\nabla_1 F(x_{k+1},y_{k+1};\xi)-\nabla_1 F(x_k,y_k;\xi)\right\|^2\right]
+2E_k\left[\left\|\nabla_{12}^2 G(x_{k+1},y_{k+1};\zeta)(z_{k+1}-z_{k}) \right\|^2\right]
\\&+2E_k\left[\left\|(\nabla_{12}^2 G(x_{k+1},y_{k+1};\zeta)-\nabla_{12}^2 G(x_k,y_k;\zeta))z_{k} \right\|^2\right]
\\\leq&
\Tilde{\Delta}\alpha_k^2E_k\left[\left\|v_k^x \right\|^2\right]
+\Tilde{\Delta}\beta_k^2E_k\left[\left\|v_k^y \right\|^2\right]
+2(L_1^g)^2\gamma_k^2E_k\left[\left\|v_k^z \right\|^2\right],
\end{align*}
where $\Tilde{\Delta}=(L^f)^2+2(L_2^g)^2R^2$.
Based on Lemma \ref{dL}, we have
\begin{align*}
	&E_k\left[\left\|v_{k+1}^x-D_x\left(x_{k+1}, y_{k+1}, z_{k+1}\right)\right\|^2\right] \\
	\leq&\left(1-\rho_k^x\right)^2 E_k\left[\left\|v_k^x-D_x\left(x_k, y_k, z_k\right)\right\|^2\right]  +2\left(\rho_k^x\right)^2 \left(\sigma_{f}+R^2\sigma_{g,2}\right)  \\
	&+2\left(1-\rho_k^x\right)^2 \left(\Tilde{\Delta}\alpha_k^2E_k[\|v_k^x \|^2]
	+\Tilde{\Delta}\beta_k^2E_k[\|v_k^y \|^2]
	+2(L_1^g)^2\gamma_k^2E_k[\|v_k^z \|^2] \right)\\
	\leq&
	\left(1-\rho_k^x\right)^2 E_k\left[\left\|v_k^x-D_x\left(x_k, y_k, z_k\right)\right\|^2\right]  +2\left(\rho_k^x\right)^2 \left(\sigma_{f}+R^2\sigma_{g,2}\right)  
	+2\left(1-\rho_k^x\right)^2 \Tilde{\Delta}\alpha_k^2E_k[\|v_k^x \|^2]
	\\&
	+4\left(1-\rho_k^x\right)^2\Tilde{\Delta}\beta_k^2 E_k[\|v_k^y-D_y(x_k,y_k,z_k) \|^2]
	+4\left(1-\rho_k^x\right)^2(L_1^g)^2\Tilde{\Delta}\beta_k^2 E_k[\|y_k-y^*(x_k) \|^2]
	\\&
	+8\left(1-\rho_k^x\right)^2\gamma_k^2(L_1^g)^2E_k[\|v_k^z -D_z(x_k,y_k,z_k)\|^2] 
	+8\left(1-\rho_k^x\right)^2\gamma_k^2(L_1^g)^2L_z^2 E_k[\|y_k-y^*(x_k)\|^2]
	\\&
	+8\left(1-\rho_k^x\right)^2\gamma_k^2(L_1^g)^2L_z^2 E_k[\|z_k-z^*(x_k)\|^2],
\end{align*}
By taking the total expectation and rearranging the above expression, the (2) is proven.

We can similarly prove (3).
\end{proof}

\begin{theorem}\textbf{(Restatement of Theorem \ref{thstorm})}

    Fix an iteration $K>1$ and assume that Assumption \ref{assump UL} to Assumption \ref{assumptionijsmooth} hold.
Then there exist positive constants $\eta$, $c_{\beta}$, $c_{\gamma}$, $c_x$, $c_y$ and $c_z$ such that if
\begin{gather*}
    \alpha_k=\frac{1}{(\eta+k)^{1/3}},\quad \beta_k=c_{\beta}\alpha_k,\quad \gamma_k=c_{\gamma}\alpha_k;\quad
    \rho_k^x=c_x \alpha_k^2,\quad  \rho_k^y=c_y \alpha_k^2,\quad \rho_k^z=c_z \alpha_k^2,
\end{gather*}
the iterates in SRMBA satisfy
\begin{align*}
\frac{1}{K}\sum_{k=0}^{K-1} \mathbb{E}\left[\left\|\nabla H\left(x_k\right)\right\|^2\right]
=\mathcal{O}\left(\frac{1}{K^{2/3} }+\frac{\log(K-1)\sigma}{K^{2/3}}\right).
\end{align*}
\end{theorem}
\begin{proof}
    We consider the Lyapunov function
\begin{align*}
L_k=&\mathbb{E}\left[H\left(x_k\right)\right]
+A\mathbb{E}\left[\left\|y_k-y^*\left(x_k\right)\right\|^2\right]
+B\mathbb{E}\left[\left\|z_k-z^*\left(x_k\right)\right\|^2\right]
\\
&+\frac{1}{C_{k-1}}\mathbb{E}[\| v^y_k-D_y(x_k,y_k,z_k)\|^2]
+\frac{1}{D_{k-1}}\mathbb{E}[\|v^x_k-D_x(x_k,y_k,z_k)\|^2]
+\frac{1}{F_{k-1}}\mathbb{E}[\| v^z_k-D_z(x_k,y_k,z_k)\|^2].
\end{align*}
Using Lemma \ref{DH} Lemma \ref{Hlemma}, Lemma \ref{y,z}, Lemma \ref{vh} and Lemma \ref{S}, we get

\begin{align*}
 L_{k+1}-L_k=&\mathbb{E}\left[H\left(x_{k+1}\right)\right]-E\left[H\left(x_k\right)\right]
+A\left(\mathbb{E}\left[\left\|y_{k+1}-y^*\left(x_{k+1}\right)\right\|^2\right]-\mathbb{E}\left[\left\|y_k-y^*\left(x_k\right)\right\|^2\right]\right) \\
&+B\left(\mathbb{E}\left[\left\|z_{k+1}-z^*\left(x_{k+1}\right)\right\|^2\right]-\mathbb{E}\left[\left\|z_k-z^*\left(x_k\right)\right\|^2\right]\right)
\\&+\frac{1}{C_{k}}\mathbb{E}\left[\left\| v^y_{k+1}-D_y(x_{k+1},y_{k+1},z_{k+1})\right\|^2\right]
-\frac{1}{C_{k-1}}\mathbb{E}\left[\left\|v^y_k-D_y(x_k,y_k,z_k)\right\|^2\right]
\\&+\frac{1}{D_{k}}\mathbb{E}\left[\left\| v^x_{k+1}-D_x(x_{k+1},y_{k+1},z_{k+1})\right\|^2\right]
-\frac{1}{D_{k-1}}\mathbb{E}\left[\left\| v^x_k-D_x(x_k,y_k,z_k)\right\|^2\right]\\
&+\frac{1}{F_{k}}\mathbb{E}\left[\left\| v^z_{k+1}-D_z(x_{k+1},y_{k+1},z_{k+1})\right\|^2\right]
-\frac{1}{F_{k-1}}\mathbb{E}\left[\left\| v^z_k-D_z(x_k,y_k,z_k)\right\|^2\right]\\
\end{align*}

By incorporating Corollary \ref{Hlemmabaised}, Lemma \ref{eg30}, and Lemma \ref{eg33} into the above inequality, we can derive
\begin{align*}
& L_{k+1}-L_k
 \\
 \leq&
 -\frac{\alpha_k}{2} \mathbb{E}\left[\left\|\nabla H\left(x_k\right)\right\|^2\right]
\\
 &
    +\Big( 3\left(\left(L^f\right)^2+\left(L^g_2 R\right)^2\right)\alpha_k-\frac{\mu L_1^g\beta_k}{2(\mu+L_1^g)}A
    +
    \frac{18((L_2^gR)^2+(L^f)^2)\gamma_k}{\mu}B
    +\frac{4(L_1^g)^2\beta_k^2(1-\rho^y_k)^2}{C_k}
    \\&
    +\big(\frac{(1-\rho^x_k)^2}{D_k} + \frac{(1-\rho^z_k)^2}{F_k}\big)(8L_z^2(L_1^g)^2\gamma_k^2+4(L_1^g)^2\Tilde{\Delta}\beta_k^2)
    \Big)\times
       \mathbb{E}\left[\left\|y_k-y^*\left(x_k\right)\right\|^2\right]
 \\
    &+\left(-\frac{\mu\gamma_k}{4}B+8\gamma_k^2(L_1^g)^2L_z^2\left(\frac{(1-\rho^x_k)^2}{D_k} + \frac{(1-\rho^z_k)^2}{F_k}\right)+3\left(L^g_1\right)^2\alpha_k\right)
        \mathbb{E}\left[\left\|z_k-z^*\left(x_k\right)\right\|^2\right]
 \\
    &+\left({\alpha_k} +\frac{\left(1-\rho_k^x\right)^2 }{D_k}-\frac{1}{D_{k-1}}\right)
        \mathbb{E}\left[\left\|v_k^x-D_x(x_k,y_k,z_k)\right\|^2\right]
 \\
    &+\left(6\frac{\mu+L_1^g}{\mu L_1^g}\beta_k A +\frac{\left(1-\rho_k^y\right)^2+4\left(1-\rho_k^y\right)^2\left(L_1^g\right)^2\beta_k^2}{C_k}-\frac{1}{C_{k-1}}+\frac{4\left(1-\rho_k^x\right)^2\Tilde{\Delta}\beta_k^2}{D_k}+\frac{4\left(1-\rho_k^z\right)^2\Tilde{\Delta}\beta_k^2}{F_{k}}\right)\times\\
       &\quad\quad
       \mathbb{E}\left[\left\| v_k^y-D_y(x_k,y_k,z_k)\right\|^2\right]
\\
    &+\left(\frac{12B\gamma_k}{\mu}  +\frac{8\left(1-\rho_k^x\right)^2\left(L_1^g\right)^2\gamma_k^2}{D_k}
    +\frac{\left(1-\rho_k^z\right)^2
    +8\left(1-\rho_k^z\right)^2\left(L_1^g\right)^2\gamma_k^2}{F_k}-\frac{1}{F_{k-1}}\right) \mathbb{E}\left[\| v_k^z-D_z(x_k,y_k,z_k)\|^2\right]
\\
    &+ \left(\frac{L^H\alpha_k^2}{2}-\frac{\alpha_k}{2}+\frac{2(\mu+L_1^g)L_{y^*}^2\alpha_k^2}{\mu L_1^g\beta_k} A
    +\frac{4L_{z^*}^2\alpha_k^2}{\mu \gamma_k}B\right.\\
    &\left.\quad+\frac{2\left(1-\rho_k^y\right)^2\left(L_1^g\right)^2\alpha_k^2 }{C_k}+\frac{2\left(1-\rho_k^x\right)^2 \Tilde{\Delta}\alpha_k^2}{D_k}+\frac{2\left(1-\rho_k^z\right)^2 \Tilde{\Delta}\alpha_k^2}{F_k}\right)
        \mathbb{E}\left[\left\|v_k^x\right\|^2\right]
\\
    &+\frac{2\left(\rho_k^y\right)^2 \sigma_{g,1}^2}{C_k}
    +\frac{2\left(\rho_k^x\right)^2\left(2 \sigma_{g,2}^2R^2+\sigma_f^2\right)}{D_k}
+\frac{2\left(\rho_k^z\right)^2\left(2 \sigma_{g,2}^2R^2+\sigma_f^2\right)}{F_k},
\end{align*}
We select the coefficients of the Lyapunov function and the step sizes of the algorithm as follows
\begin{gather*}
    \beta_k=\frac{1}{(\eta+k)^{1/3}},\quad \beta_k\leq C_{\beta},\quad \alpha_k=c_{\alpha}\beta_k,\quad \gamma_k=c_{\gamma}\beta_k;\quad
    \rho_k^x=c_x \beta_k^2,\quad  \rho_k^y=c_y \beta_k^2,\quad \rho_k^z=c_z \beta_k^2;\\
    A=\frac{\mu L_1^g}{12(\mu+L_1^g)},\quad B=1,\quad C_{k}=\frac{ \beta_k}{\phi_c},\quad \quad D_{k}=\frac{ \beta_k}{\phi_d},\quad \quad F_{k}=\frac{ \beta_k}{\phi_f}.
\end{gather*}

Based on the definition of $\beta_k$ and the choice of $C_k$
\begin{align*}
	\frac{\left(1-\rho_k^y\right)^2 }{C_k}-\frac{1}{C_{k-1}}
	&\leq
	\frac{1}{C_k}-\frac{1}{C_{k-1}}-\frac{\rho_k^y }{C_k}
	=\phi_c\left(\frac{1}{\beta_k}-\frac{1}{\beta_{k-1}}-c_y\beta_k\right)\\
	&=\phi_c\left((\eta+k)^{1/3}-(\eta+k-1)^{1/3}-c_y\beta_k\right)\\
	&\leq \phi_c\left(\frac{2^{2/3}}{3(\eta+k)^{2/3}}-c_y\beta_k\right)\\
	&\leq \phi_c\left(\frac{2^{2/3}}{3}\beta_k^2-c_y\beta_k\right)\\
	&
	\leq -\beta_k,
\end{align*}
where the second inequality follows from $(x+y)^{1/3}-x^{1/3}\leq y/(3x^{2/3})$ and $\eta\geq 2$,
the third inequality is based on the definition of $\alpha_k$,
and the final inequality results from our choice of $\beta_k\leq C_\beta$ and $c_y=\frac{1}{\phi_c}+{C_{\beta}}$.
Similarly, we have
\begin{align*}
	\frac{\left(1-\rho_k^x\right)^2 }{D_k}-\frac{1}{D_{k-1}}\leq  -\beta_k,\quad
		\frac{\left(1-\rho_k^z\right)^2 }{F_k}-\frac{1}{F_{k-1}}\leq  -\beta_k
\end{align*}

\textbf{Analysis of the Coefficient for} $\mathbb{E}\left[\left\|y_k-y^*\left(x_k\right)\right\|^2\right] $

Due to the assumption that
\begin{align*}
     &\alpha_k\leq \frac{3A^2}{\left(L^f\right)^2+\left(L^g_2 R\right)^2}\beta_k,
     \quad
     \gamma_k\leq \frac{\mu A^2}{2\left(L^f\right)^2+\left(L^g_2 R\right)^2}\beta_k,\\&
     \phi_c\leq \frac{9A^2}{4(L_1^g)^2},
     \quad
     \phi_dc_\gamma^2\leq\frac{9A^2}{8L_z^2(L_1^g)^2},
     \quad
     \phi_fc_\gamma^2\leq\frac{9A^2}{8L_z^2(L_1^g)^2},\\
     &\phi_d\leq \frac{9A^2}{4(L_1^g)^2\Tilde{\Delta}},\quad
     \phi_d\leq \frac{9A^2}{4(L_1^g)^2\Tilde{\Delta}},
\end{align*}
it follows that
\begin{align*}
   & 3\left(\left(L^f\right)^2+\left(L^g_2 R\right)^2\right)\alpha_k-\frac{\mu L_1^g\beta_k}{2(\mu+L_1^g)}A
    +
    \frac{18((L_2^gR)^2+(L^f)^2)\gamma_k}{\mu}B
    +\frac{4(L_1^g)^2\beta_k^2(1-\rho^y_k)^2}{C_k}
    \\&
    +\big(\frac{(1-\rho^x_k)^2}{D_k} + \frac{(1-\rho^z_k)^2}{F_k}\big)(8L_z^2(L_1^g)^2\gamma_k^2+4(L_1^g)^2\Tilde{\Delta}\beta_k^2)
    \\&=3\left(\left(L^f\right)^2+\left(L^g_2 R\right)^2\right)\alpha_k-72A^2\beta_k
    +
    \frac{18((L_2^gR)^2+(L^f)^2)\gamma_k}{\mu}B
    +\frac{4(L_1^g)^2\beta_k^2(1-\rho^y_k)^2}{C_k}
    \\&
    +\big(\frac{(1-\rho^x_k)^2}{D_k} + \frac{(1-\rho^z_k)^2}{F_k}\big)(8L_z^2(L_1^g)^2\gamma_k^2+4(L_1^g)^2\Tilde{\Delta}\beta_k^2)
    \\&\leq 0.
\end{align*}
\textbf{Analysis of the Coefficient for} $ \mathbb{E}\left[\left\|z_k-z^*\left(x_k\right)\right\|^2\right]$

By assuming
\begin{align*}
\alpha_k\leq \frac{\mu}{36(L_1^g)^2}\gamma_k,\quad
c_\gamma\phi_d\leq \frac{\mu}{96(L_1^g)^2L_z62},\quad
c_\gamma\phi_f\leq \frac{\mu}{96(L_1^g)^2L_z62},
\end{align*}
we have
\begin{align*}
   & -\frac{\mu\gamma_k}{4}B+8\gamma_k^2(L_1^g)^2L_z^2\left(\frac{(1-\rho^x_k)^2}{D_k} + \frac{(1-\rho^z_k)^2}{F_k}\right)+3\left(L^g_1\right)^2\alpha_k
    \leq 0.
\end{align*}
\textbf{Analysis of the Coefficient for} $\mathbb{E}\left[\left\|v_k^x-D_x(x_k,y_k,z_k)\right\|^2\right]$

By assuming $\alpha_k\le q\beta_k,$ we have
\begin{align*}
    {\alpha_k} +\frac{\left(1-\rho_k^x\right)^2 }{D_k}-\frac{1}{D_{k-1}}\leq
    {\alpha_k}-\beta_k\leq  0.
\end{align*}

\textbf{Analysis of the Coefficient for} $\mathbb{E}\left[\left\| v_k^y-D_y(x_k,y_k,z_k)\right\|^2\right]$

Furthermore, by assuming
\begin{align*}
    \phi_c\leq \frac{1}{24(L_1^g)^2},\quad
    \phi_d\leq \frac{1}{24\Tilde{\Delta}},\quad
    \phi_f\leq \frac{1}{24\Tilde{\Delta}},
\end{align*}
we have
\begin{align*}
  & 6\frac{\mu+L_1^g}{\mu L_1^g}\beta_k A +\frac{\left(1-\rho_k^y\right)^2+4\left(1-\rho_k^y\right)^2\left(L_1^g\right)^2\beta_k^2}{C_k}-\frac{1}{C_{k-1}}+\frac{4\left(1-\rho_k^x\right)^2\Tilde{\Delta}\beta_k^2}{D_k}+\frac{4\left(1-\rho_k^z\right)^2\Tilde{\Delta}\beta_k^2}{F_{k}}
  \\&
  \leq -\frac{\beta_k}{2} +\frac{4\left(1-\rho_k^y\right)^2\left(L_1^g\right)^2\beta_k^2}{C_k}+\frac{4\left(1-\rho_k^x\right)^2\Tilde{\Delta}\beta_k^2}{D_k}+\frac{4\left(1-\rho_k^z\right)^2\Tilde{\Delta}\beta_k^2}{F_{k}}
  \\&\leq 0.
\end{align*}

\textbf{Analysis of the Coefficient for} $\mathbb{E}\left[\| v_k^z-D_z(x_k,y_k,z_k)\|^2\right]$

Similar to the analysis of $\mathbb{E}\left[\left\| v_k^y-D_y(x_k,y_k,z_k)\right\|^2\right]$,
by assuming
\begin{align*}
 \gamma_k\leq\frac{\mu }{36}\beta_k,\quad
    \phi_d c_\gamma^2\leq \frac{1}{24(L_1^g)^2},\quad
    \phi_f c_\gamma^2\leq \frac{1}{24(L_1^g)^2},
\end{align*}
we have
\begin{align*}
    &\frac{12B\gamma_k}{\mu}  +\frac{8\left(1-\rho_k^x\right)^2\left(L_1^g\right)^2\gamma_k^2}{D_k}
    +\frac{\left(1-\rho_k^z\right)^2
    	+8\left(1-\rho_k^z\right)^2\left(L_1^g\right)^2\gamma_k^2}{F_k}-\frac{1}{F_{k-1}}\\
    	\leq&\frac{12\gamma_k}{\mu}  +\frac{8\left(1-\rho_k^x\right)^2\left(L_1^g\right)^2\gamma_k^2}{D_k}
    	+\frac{8\left(1-\rho_k^z\right)^2\left(L_1^g\right)^2\gamma_k^2}{F_k}-\beta_k\\
    	\leq& 0.
\end{align*}

\textbf{Analysis of the Coefficient for} $\mathbb{E}\left[\left\|v_k^x\right\|^2\right]$

By assuming that the step sizes and the parameters of the Lyapunov function satisfy
\begin{align*}
    \alpha_k\leq \min\left\{\frac{1}{6L^H},
    \,\frac{1}{L_{y^*}}\beta_k,
    \,\frac{1}{24L_{z^*}}\gamma_k\right\},\quad
    \phi_c c_\alpha\leq \frac{1}{12(L_1^g)^2},\quad\phi_d c_\alpha\leq \frac{1}{12\Tilde{\Delta}},\quad\phi_f c_\alpha\leq \frac{1}{12\Tilde{\Delta}},
\end{align*}
we have
\begin{align*}
    &\frac{L^H\alpha_k^2}{2}-\frac{\alpha_k}{2}+\frac{2(\mu+L_1^g)L_{y^*}^2\alpha_k^2}{\mu L_1^g\beta_k} A
    +\frac{4L_{z^*}^2\alpha_k^2}{\mu \gamma_k}B\\
    &+\frac{2\left(1-\rho_k^y\right)^2\left(L_1^g\right)^2\alpha_k^2 }{C_k}+\frac{2\left(1-\rho_k^x\right)^2 \Tilde{\Delta}\alpha_k^2}{D_k}+\frac{2\left(1-\rho_k^z\right)^2 \Tilde{\Delta}\alpha_k^2}{F_k}\leq 0.
\end{align*}
Combining the above analysis of the coefficients of each term, we can obtain the simplified inequality as follows
\begin{align*}
 L_{k+1}-L_k
 \leq&
 -\frac{\alpha_k}{2} \mathbb{E}\left[\left\|\nabla H\left(x_k\right)\right\|^2\right]+\frac{2\left(\rho_k^y\right)^2 \sigma_{g,1}^2}{C_k}
    +\frac{2\left(\rho_k^x\right)^2\left(2 \sigma_{g,2}^2R^2+\sigma_f^2\right)}{D_k}
+\frac{2\left(\rho_k^z\right)^2\left(2 \sigma_{g,2}^2R^2+\sigma_f^2\right)}{F_k},
\end{align*}
Summing, taking the average, and rearranging, we obtain
\begin{align*}
\frac{1}{K}\sum_{k=0}^{K-1}\frac{\alpha_k}{2} \mathbb{E}\left[\left\|\nabla H\left(x_k\right)\right\|^2\right]
= \mathcal{O}\left(\frac{L_0}{K}+\frac{1}{K}\sum_{k=0}^{K-1}
\alpha_k^3
\right)
=\mathcal{O}\left(\frac{1}{K^{2/3} }+\frac{\log(K-1)}{K^{2/3}}\right).
\end{align*}
\end{proof}

\section{MA-SOBA-q: Vanilla minibatch SGD + Standard Momentum} \label{detial:proofsoba}
For the expectation form setting, we introduce MA-SOBA-q, which employs mini-batch stochastic estimation in its estimation module and selects  Moving-average in the acceleration module that reference the update direction from the previous iteration.

To illustrate further, at each iteration $k$, we draw two random set $\mathcal{S}_1$ and $\mathcal{S}_2$ with a fixed mini-batch size of $S$, for the functions $f$ and $g$ respectively, to perform a stochastic estimation of $D_{\bullet}$.
$\gamma_k$, $\beta_k$ and $\alpha_k$ are the step sizes and
$\rho_k$ is the moving average parameter.
A trade-off between the step size and batch size has been made, with more detailed descriptions to be provided in Theorem \ref{th1}.
Furthermore, we introduce historical information $v_{k-1}^x$ and $u_{k-1}$, and employ the moving average technique for acceleration, specifically by forming a convex combination of $v_{k-1}^x$ and $D_{k-1}^x$.

\begin{theorem}\textbf{(Expection form problem (\ref{generalpro}))}\label{th1}

Fix an iteration $K>1$ and assume that Assumption \ref{assump UL} to \ref{assumporacle} hold.
The mini-batch size $S$ is chosen to be ${K^q}$.
Then there exist positive constants $c_{\alpha}$, $c_{\beta}$, $c_{\gamma}$ and $c_{\rho}$
such that if
\begin{gather*}\alpha_k=c_{\alpha}K^{-p},\quad\beta_k=c_{\beta}K^{-p},\\\gamma_k=c_{\gamma}K^{-p},\quad\rho_k=c_{\rho}K^{-p},\end{gather*}
the iterates in MA-SOBA-q satisfy
\begin{eqnarray*}
\frac{1}{K} \sum_{k=1}^K \mathbb{E}\left[\left\|\nabla H\left(x_k\right)\right\|^2\right]=\mathcal{O}\left(\frac{1}{K^{1-p}}+\frac{1}{K^{p+q}}\right).
\end{eqnarray*}
\end{theorem}
\begin{remark}
In Theorem \ref{th1}, we discuss the trade-off between step sizes and mini-batch sizes, their exponents need to satisfy $q+2p=1$, ensuring that the sampling complexity is $\mathcal{O}(\epsilon^{-2})$.
\end{remark}

\begin{algorithm}[h]
  \caption{MA-SOBA-q}
\begin{algorithmic}[1]
   \STATE {\bfseries Input:} Initializations $(x_{-1},y_{-1},z_{-1})$, $(x_{0},y_{0},z_{0})$, and $v_{-1}^x$, number of total iterations $K$, step size $\{\alpha_k, \beta_k,\gamma_k\}$, momentum parameter$\rho_k$;

   \FOR{$k=0$ {\bfseries to} $K-1$}
   \STATE Sample $\mathcal{S}_1$ for $f$ and $\mathcal{S}_2$ for $g$;
    \STATE
        $v^x_k=(1-\rho_{k-1})v_{k-1}^x+\rho_{k-1} D^x_{k-1}$;
   \STATE
       $x_{k+1}= x_{k}-\alpha_k v_k^x$;
    \STATE
        $D^x_{k}=\nabla_1 F(x_{k},y_{k};\mathcal{S}_1)-\nabla_{12}^2G(x_{k},y_{k};\mathcal{S}_2)z_{k}$;
   \STATE
        $v^y_k=\nabla_2G(x_{k},y_{k};\mathcal{S}_2)$;
   \STATE
        $y_{k+1} = y_{k}-\beta_k v_k^y$;
   \STATE
        $ v^z_k=\nabla_{22}^2G(x_{k},y_{k};\mathcal{S}_2)z_{k}- \nabla_2 F(x_{k},y_{k};\mathcal{S}_1)$;
   \STATE
       $z_{k+1} = z_{k}-\gamma_k v_k^z$.

   \ENDFOR
\end{algorithmic}
\end{algorithm}

\begin{lemma}\label{eg1}
Under the Assumption \ref{assump UL} to \ref{assumporacle}, we have
\begin{align*}
    \mathbb{E}\left[\left\|D_y(x_k,y_k,z_k)-D_k^y\right\|^2\right]
    &\leq \frac{ \sigma_{g,1}^2  }{|\mathcal{S}_2|}, \\
    \mathbb{E}\left[\left\|D_z(x_k,y_k,z_k)-D_k^z\right\|^2\right]
     &\leq \frac{2 \sigma_{g,2}^2  }{|\mathcal{S}_2|}(R^2+\mathbb{E}\left[\left\|z_k-z^*(x_k)\right\|^2\right])
    +\frac{\sigma_{f}^2}{|\mathcal{S}_1|},\\
    \mathbb{E}\left[\left\|D_x(x_k,y_k,z_k)-D_k^x\right\|^2\right]
     &\leq \frac{2 \sigma_{g,2}^2  }{|\mathcal{S}_2|}(R^2+\mathbb{E}\left[\left\|z_k-z^*(x_k)\right\|^2\right])
    +\frac{\sigma_{f}^2}{|\mathcal{S}_1|}.
\end{align*}
\end{lemma}
\begin{proof}
Based on the definition of $D_k^y$ and Assumption \ref{assumporacle}, we have
\begin{eqnarray*}
\mathbb{E}\left[\left\|D_y(x_k,y_k,z_k)-D_k^y\right\|^2\right]
=\mathbb{E}[\|D_y(x_k,y_k,z_k)-\frac{1}{|\mathcal{S}_2|}\sum_{\zeta\in \mathcal{S}_2}\nabla_2G_j(x_k,y_k;\zeta )\|^2]
    \leq \frac{ \sigma_{g,1}^2  }{|\mathcal{S}_2|}.
\end{eqnarray*}
\begin{eqnarray*}
\mathbb{E}\left[\left\|D_z(x_k,y_k,z_k)-D_k^z\right\|^2\right]
&=&\mathbb{E}\left[\left\| \nabla_{22}^2g(x_k,y_k)z_k- \nabla_2 f(x_k,y_k)-\nabla_{22}^2G(x_{k},y_{k};\mathcal{S}_2)z_{k}+ \nabla_2 F(x_{k},y_{k};\mathcal{S}_1)\right\|^2\right]
\\&=&
\mathbb{E}\left[\left\| \nabla_{22}^2g(x_k,y_k)z_k-\nabla_{22}^2G(x_{k},y_{k};\mathcal{S}_2)z_{k}\right\|^2\right]+
\mathbb{E}\left[\left\|\nabla_2 f(x_k,y_k)-\nabla_2 F(x_{k},y_{k};\mathcal{S}_1)\right\|^2\right]
\\&\leq&
2\mathbb{E}\left[\left\| \nabla_{22}^2g(x_k,y_k)-\nabla_{22}^2G(x_{k},y_{k};\mathcal{S}_2)\right\|^2\right](\mathbb{E}[\|z_k-z^*(x_k)\|^2]+\mathbb{E}[\|z^*(x_k)\|^2])
\\&&+
\mathbb{E}\left[\left\|\nabla_2 f(x_k,y_k)-\nabla_2 F(x_{k},y_{k};\mathcal{S}_1)\right\|^2\right] 
\\&\leq&
\frac{2 \sigma_{g,2}^2  }{|\mathcal{S}_2|}(R^2+\mathbb{E}\left[\left\|z_k-z^*(x_k)\right\|^2\right])+\frac{\sigma_{f}^2}{|\mathcal{S}_1|},
\end{eqnarray*}
where the second equation holds because mini-batch estimation is unbiased. The last inequation is due to Assumption \ref{assumporacle} and Lemma \ref{R}.
Similarly, we can obtain inequalities regarding $D_k^x$.
\end{proof}

\begin{theorem}\textbf{(Restatement of Theorem \ref{th1})}

Fix an iteration $K>1$ and assume that Assumption \ref{assump UL} to \ref{assumporacle} hold. Let $|\mathcal{S}_1|=|\mathcal{S}_2|={K^q}$.
The step sizes $\alpha_k$, $\beta_k$, $\gamma_k$, and $\rho_k$ have the same order of $\Theta(K^{-p})$, $p>0$, and satisfy
\begin{align*}
    \alpha_k\leq&\min\{\frac{1}{2L^H},\frac{1}{16L^2_{y^*}}\beta_k,\frac{1}{64(L^H)^2}\rho_k,2\rho_k\},\quad
    \beta_k\leq\min\{\frac{4}{\mu},\frac{\mu}{16\Delta}\beta_k\},\\
    \rho_k\leq&\min\{\frac{\mu^2}{24\Delta}\beta_k,\frac{\mu^2}{24(L_1^g)^2}\gamma_k,1\},\quad
    \rho_k^2\leq\frac{\mu^2}{8}\gamma_k.
\end{align*}
Then the iterates in MA-SOBA-q satisfy
\begin{eqnarray*}
\frac{1}{K} \sum_{k=1}^K \mathbb{E}\left[\left\|\nabla H\left(x_k\right)\right\|^2\right]=\mathcal{O}(\frac{1}{K^{1-p}}+\frac{1}{K^{p+q}}).
\end{eqnarray*}
\end{theorem}
\begin{proof}
Consider the Lyapunov function in the form of
\begin{eqnarray}\label{L}
L_k=\mathbb{E}\left[H\left(x_k\right)\right]
+A\mathbb{E}\left[\left\|y_k-y^*\left(x_k\right)\right\|^2\right]
+B\mathbb{E}\left[\left\|z_k-z^*\left(x_k\right)\right\|^2\right]
+C\mathbb{E}[\|\nabla H(x_k)-v^x_k\|^2].
\end{eqnarray}
In Lemma \ref{y,z}, we provide the descent lemma for the second and third terms in (\ref{L}), and Lemma \ref{Hlemma} provides the descent lemma for the first term.
Related to the last term, refer to Lemma \ref{vh}.
Therefore, we have
\begin{align*}
 L_{k+1}-L_k=&\mathbb{E}\left[H\left(x_{k+1}\right)\right]-E\left[H\left(x_k\right)\right]
+A\left(\mathbb{E}\left[\left\|y_{k+1}-y^*\left(x_{k+1}\right)\right\|^2\right]-\mathbb{E}\left[\left\|y_k-y^*\left(x_k\right)\right\|^2\right]\right) \\
&+B\left(\mathbb{E}\left[\left\|z_{k+1}-z^*\left(x_{k+1}\right)\right\|^2\right]-\mathbb{E}\left[\left\|z_k-z^*\left(x_k\right)\right\|^2\right]\right)
\\&+C\left(\mathbb{E}[\|\nabla H(x_{k+1})-v^x_{k+1}\|^2]-\mathbb{E}[\|\nabla H(x_k)-v^x_k\|^2]\right)\\
&\leq
-\frac{\alpha_k}{2} \mathbb{E}\left[\left\|\nabla H\left(x_k\right)\right\|^2\right]
\\&
+\left(\frac{L^H\alpha_k^2}{2}-\frac{\alpha_k}{2}+A\frac{2 L_{y^*}^2 \alpha_k^2}{\beta_k \mu}+B\frac{3 L_{z^*}^2 \alpha_k^2}{\gamma_k \mu}
+ C\frac{2\left(L^H\right)^2 \alpha_k^2}{\rho_k}\right)
\mathbb{E}\left[\left\|v_k^x\right\|^2\right]
\\&+\left(-A\beta_k \mu +8\Delta B\gamma_k+6C\Delta\rho_k\right)
\mathbb{E}\left[\left\| y_k-y^* (x_k)\right\|^2\right]
\\&
+\left(-B\gamma_k\mu + 2B\gamma_k^2\frac{2\sigma_{g,2}^2}{|\mathcal{S}_2|}+6C\left(L_1^g\right)^2\rho_k+\rho_k^2C\frac{2\sigma_{g,2}^2}{|\mathcal{S}_2|}\right)
\mathbb{E}\left[\left\|z_k-z^*\left(x_k\right)\right\|^2\right]
\\&
+\left(\frac{\alpha_k}{2}-C\rho_k\right)
\mathbb{E}\left[\left\|\nabla H\left(x_k\right)-v_k^x\right\|^2\right]
\\&
+2A\beta_k^2\frac{\sigma_{g,1}^2}{|\mathcal{S}_2|}
+2B\gamma_k^2\left(\frac{2\sigma_{g,2}^2}{|\mathcal{S}_2|}R^2+\frac{2\sigma_{f}^2}{|\mathcal{S}_1|}\right)
+2C\rho_k^2\left(\frac{2\sigma_{g,2}^2}{|\mathcal{S}_2|}R^2+\frac{2\sigma_{f}^2}{|\mathcal{S}_1|}\right),
\end{align*}
where the inequality holds by utilizing Lemma \ref{DH} and Lemma \ref{eg1}.
Furthermore, we have
\begin{align}\label{eg1.1}
 L_{k+1}-L_k
&\leq
-\frac{\alpha_k}{2} \mathbb{E}\left[\left\|\nabla H\left(x_k\right)\right\|^2\right]
+2A\beta_k^2\frac{\sigma_{g,1}^2}{|\mathcal{S}_2|}
+2B\gamma_k^2\left(\frac{2\sigma_{g,2}^2}{|\mathcal{S}_2|}R^2+\frac{2\sigma_{f}^2}{|\mathcal{S}_1|}\right)
+2C\rho_k^2\left(\frac{2\sigma_{g,2}^2}{|\mathcal{S}_2|}R^2+\frac{2\sigma_{f}^2}{|\mathcal{S}_1|}\right),
\end{align}
if the following system of inequalities holds
\begin{align}\label{step}
\left\{\begin{aligned}
\frac{L^H\alpha_k^2}{2}-\frac{\alpha_k}{2}+A\frac{2 L_{y^*}^2 \alpha_k^2}{\beta_k \mu}+B\frac{3 L_{z^*}^2 \alpha_k^2}{\gamma_k \mu}
+ C\frac{2\left(L^H\right)^2 \alpha_k^2}{\rho_k}\leq 0, \\
-A\beta_k \mu +8\Delta B\gamma_k+6C\Delta\rho_k \leq 0, \\
-B\gamma_k\mu + 2B\gamma_k^2\frac{2\sigma_{g,2}^2}{|\mathcal{S}_2|}+6C\left(L_1^g\right)^2\rho_k+\rho_k^2C\frac{2\sigma_{g,2}^2}{|\mathcal{S}_2|} \leq 0, \\
\frac{\alpha_k}{2}-C\rho_k \leq 0.
\end{aligned}\right.
\end{align}
We choose the coefficients of the Lyapunov function to be $A=B=\mu$, $C=1$.
In fact, it is only necessary to require that $\alpha_k$, $\beta_k$, $\gamma_k$
, and $\rho_k$ have the same order of $\Theta(K^{-p})$, $p>0$, and satisfy
\begin{align*}
    \alpha_k\leq&\min\{\frac{1}{2L^H},\frac{1}{16L^2_{y^*}}\beta_k,\frac{1}{64(L^H)^2}\rho_k,2\rho_k\},\,
    \beta_k\leq\min\{\frac{4}{\mu},\frac{\mu}{16\Delta}\beta_k\},\,
    \rho_k\leq&\min\{\frac{\mu^2}{24\Delta}\beta_k,\frac{\mu^2}{24(L_1^g)^2}\gamma_k,1\},\,
    \rho_k^2\leq\frac{\mu^2}{8}\gamma_k,
\end{align*}
then (\ref{step}) holds.
Rearranging (\ref{eg1.1}), we have
\begin{align*}
\frac{\alpha_k}{2}\mathbb{E}\left[\left\|\nabla H\left(x_k\right)\right\|^2\right]
&\leq L_k-L_{k+1}
+2\mu\beta_k^2\frac{\sigma_{g,1}^2}{|\mathcal{S}_2|}
+2\mu\gamma_k^2\left(\frac{2\sigma_{g,2}^2}{|\mathcal{S}_2|}R^2+\frac{2\sigma_{f}^2}{|\mathcal{S}_1|}\right)
+2\rho_k^2\left(\frac{2\sigma_{g,2}^2}{|\mathcal{S}_2|}R^2+\frac{2\sigma_{f}^2}{|\mathcal{S}_1|}\right)
\end{align*}
Summing and telescoping  yields
\begin{align*}
\frac{1}{K}\sum_{k=0}^{K-1}\alpha_k\mathbb{E}\left[\left\|\nabla H\left(x_k\right)\right\|^2\right]
&\leq\frac{2L_0}{K}
+\frac{1}{K}\sum_{k=0}^{K-1}\left(2\mu\beta_k^2\frac{\sigma_{g,1}^2}{|\mathcal{S}_2|}
+2\mu\gamma_k^2\left(\frac{2\sigma_{g,2}^2}{|\mathcal{S}_2|}R^2+\frac{2\sigma_{f}^2}{|\mathcal{S}_1|}\right)
+2\rho_k^2\left(\frac{2\sigma_{g,2}^2}{|\mathcal{S}_2|}R^2+\frac{2\sigma_{f}^2}{|\mathcal{S}_1|}\right)\right),
\end{align*}
let $|\mathcal{S}_1|=|\mathcal{S}_2|=K^q$, then we have
\begin{align*}
\frac{1}{K}\sum_{k=0}^{K-1}\mathbb{E}\left[\left\|\nabla H\left(x_k\right)\right\|^2\right]
&= \mathcal{O}(\frac{1}{K^{1-p}}+\frac{1}{K^{p+q}})
.
\end{align*}
\end{proof}
\onecolumn

\end{document}